\documentclass[11pt, dvipdfmx]{article}
\usepackage{amssymb}
\usepackage{amsmath}
\usepackage{amsthm}
\usepackage{mathrsfs}
\usepackage{mathtools}
\usepackage{color}
\usepackage{hyperref}
\usepackage{xcolor}
\usepackage{cases}
\usepackage[normalem]{ulem}

\usepackage[scale=0.80]{geometry}

\allowdisplaybreaks[1]

\theoremstyle{plain}
\newtheorem{theorem}{Theorem}[section]
\newtheorem{lemma}[theorem]{Lemma}
\newtheorem{corollary}[theorem]{Corollary}
\newtheorem{proposition}[theorem]{Proposition}

\theoremstyle{definition}
\newtheorem{remark}[theorem]{Remark}
\newtheorem{definition}[theorem]{Definition}

\newtheorem{condition}[theorem]{Condition}

\newcommand{\floor}[1]{\lfloor #1 \rfloor}

\newcommand{\ep}{\epsilon}
\newcommand{\vep}{\varepsilon}
\newcommand{\RR}{{\mathbb R}}
\newcommand{\BB}{{\mathbb B}}
\newcommand{\RPB}{\boldsymbol{B}}
\newcommand{\XX}{{\mathbb X}}
\newcommand{\RPX}{\boldsymbol{X}}
\newcommand{\la}{\lambda}
\newcommand{\tJm}{\tilde{J}^m}
\newcommand{\tJmr}{\tilde{J}^{m,\rho}}

\newcommand{\Ymr}{Y^{m,\rho}}
\newcommand{\Ymh}{\hat{Y}^{m}}
\newcommand{\tZmr}{\tilde{Z}^{m,\rho}}
\newcommand{\Zmr}{Z^{m,\rho}}
\newcommand{\Emr}{E^{m,\rho}}
\newcommand{\Nmr}{N^{m,\rho}}
\newcommand{\tmkm}{\tau^m_{k-1}}
\newcommand{\tmk}{\tau^m_k}
\newcommand{\tmim}{\tau^m_{i-1}}
\newcommand{\tmi}{\tau^m_i}
\newcommand{\tmjm}{\tau^m_{j-1}}
\newcommand{\tmj}{\tau^m_j}
\newcommand{\tml}{\tau^m_l}
\newcommand{\tmlm}{\tau^m_{l-1}}
\newcommand{\tepmr}{\tilde{\ep}^{m,\rho}}
\newcommand{\hepm}{\hat{\ep}^{m}}
\newcommand{\epm}{\ep^m}

\newcommand{\dmJ}{\delta^{m}(J)}
\newcommand{\Dm}{D_m}
\newcommand{\bcdot}{\boldsymbol{\cdot}}
\newcommand{\hyp}{\mathchar`-}
\newcommand{\tw}{\tilde{w}}
\newcommand{\rM}{\mathrm{M}}
\newcommand{\rIM}{\mathrm{IM}}
\newcommand{\rCN}{\mathrm{CN}}
\newcommand{\rFE}{\mathrm{FE}}
\newcommand{\Omegam}{\Omega_0^{(m,d^m)}}
\numberwithin{equation}{section}

\makeatletter
\def\address#1#2{\begingroup
\noindent\parbox[t]{8.8cm}{%
\small{\scshape\ignorespaces#1}\par\vskip1ex
\noindent\small{\itshape E-mail address}%
\/: #2\par\vskip4ex}\hfill%
\endgroup}%
\makeatother
\title
{An interpolation of discrete rough differential
equations
and its applications to
analysis of error distributions}
\author{Shigeki Aida and Nobuaki Naganuma}
\date{}
\begin{document}

\maketitle

\begin{abstract}
 	We consider the solution $Y_t$ $(0\le t\le 1)$ and several approximate solutions
	$\hat{Y}^m_t$ of 
	a rough differential equation driven by 
	a fractional Brownian motion $B_t$ 
	with the Hurst parameter $1/3<H\leq 1/2$
	associated with a dyadic partition of $[0,1]$.
	We are interested in analysis of asymptotic error distribution of
	$\hat{Y}^m_t-Y_t$ as $m\to\infty$.
	In the preceding results, it was proved that the weak limit of 
	$\{(2^m)^{2H-1/2}(\hat{Y}^m_t-Y_t)\}_{0\le t\le 1}$ coincides with
	the weak limit of $\{(2^m)^{2H-1/2}J_tI^m_t\}_{0\le t\le 1}$,
	where $J_t$ is the Jacobian process of $Y_t$ and $I^m_t$ is a certain
	weighted sum process of Wiener chaos of order $2$ defined by $B_t$.
	However, it is non-trivial to reduce a problem about $\hat{Y}^m_t-Y_t$
	to one about $J_t$ and $I^m_t$.
	In this paper, we introduce an interpolation process between
	$Y_t$ and $\hat{Y}^m_t$, and give several estimates
	of the interpolation process itself and its associated processes.
	The analysis provides a framework to deal with the reduction problem
	and provides a stronger result that the difference $R^m_t=\hat{Y}^m_t-Y_t-J_tI^m_t$
	is really small compared to the main term $J_tI^m_t$.
	More precisely, we show that
	$(2^m)^{2H-1/2+\varepsilon}\sup_{0\le t\le 1}|R^m_t|\to 0$ almost surely and in $L^p$
	(for all $p>1$) for certain explicit positive number $\varepsilon>0$.
	As a consequence, we obtain an estimate of the convergence rate
	of $\sup_{0\le t\le 1}|\Ymh_t-Y_t|\to 0$ in $L^p$ also.
\end{abstract}


\textbf{MSC2020 subject classifications}: 60F05; 60H35; 60G15.

\tableofcontents

\section{Introduction}

In this paper, we study asymptotic error distributions 
for several approximation schemes
of rough differential equations(=RDEs).
Typical driving processes of RDEs are long-range correlated Gaussian processes
and we cannot use several important tools in the study of
stochastic differential equations driven by standard Brownian motions.
For example, martingale central limit theorems cannot be applied to
the study of asymptotic error distributions.
However, the fourth moment theorem can be applicable for the study of
long-range correlated Gaussian processes and several limit theorems of
weighted sum processes of Wiener chaos have been established
(\cite{NourdinNualartTudor2010, liu-tindel2020, nourdin-peccati}
and references therein).
Furthermore, these limit theorems are important in the study of
asymptotic error distributions of RDEs
(\cite{aida-naganuma2020, hu-liu-nualart2016,hu-liu-nualart2021, 
liu-tindel, naganuma2015, ueda}).
However, it is not trivial to reduce
the problem of asymptotic error distributions of solutions of RDEs
to that of weighted sum processes of Wiener chaos.
We study this problem
by introducing certain interpolation processes 
between the solution and the 
approximate solutions of RDEs.

More precisely,
we explain our main results 
and the relation with previously known results.
We consider a solution $Y_t$ of a multidimensional RDE driven by fractional Brownian
motion(=fBm) $B_t$ with the Hurst parameter $\frac{1}{3}<H\le\frac{1}{2}$,
\begin{align*}
 Y_t=\xi+\int_0^t\sigma(Y_s)d\RPB_s+\int_0^tb(Y_s)ds,\quad\quad 0\le t\le 1,
\end{align*}
where $\RPB_t$ is a naturally lifted
geometric rough path of $B_t$.
The precise meanings of rough paths and RDEs will be given in Section~\ref{statement}.
Let $\hat{Y}^m_t$ be an
approximate solution associated with the dyadic partition
$\Dm=\{\tmk\}_{k=0}^{2^m}$, where $\tmk=k2^{-m}$.
Actually there are many approximation schemes, \textit{e.g.},
the implementable Milstein, Crank-Nicolson, Milstein
and first-order Euler schemes of RDEs.
The first-order Euler scheme was introduced by 
Hu-Liu-Nualart~\cite{hu-liu-nualart2016} and
further studied by Liu-Tindel~\cite{liu-tindel}.
Among them, we explain the result in Liu and Tindel~\cite{liu-tindel}
which is closely related to our main results.
For the first-order Euler approximate solution $\Ymh_t$, they proved that
$\{(2^m)^{2H-\frac{1}{2}}(\Ymh_t-Y_t)\}_{0\le t\le 1}$ 
weakly converges to the weak limit of
$\{(2^m)^{2H-\frac{1}{2}}J_tI^m_t\}_{0\le t\le 1}$ as $m\to\infty$
in $D([0,1])$ with respect to the Skorokhod topology.
Here $J_t(=\partial_{\xi}Y_t(\xi))$ is the Jacobian (derivative) process of $Y_t$ and
$I^m_t$ is a certain weighted sum process of Wiener chaos of order 2
defined by fBm $B_t$.
Note that the weak convergence of $\{(2^m)^{2H-\frac{1}{2}}I^m_t\}$
can be proved by using the fourth moment theorem.
Their limit theorem of the error $\Ymh_t-Y_t$ is the first result 
for solutions of multidimensional RDEs with
the Hurst parameter $\frac{1}{3}<H<\frac{1}{2}$.
We are interested in the difference 
$R^m_t=\Ymh_t-Y_t-J_tI^m_t$.
The convergence results of $\{(2^m)^{2H-\frac{1}{2}}(\Ymh_t-Y_t)\}$
and $\{(2^m)^{2H-\frac{1}{2}}I^m_t\}$ suggests that
$R^m_t$ might be a small term in a certain sense as $m\to\infty$.
Conversely, if one can prove 
$\lim_{m\to\infty}E[(2^m)^{2H-\frac{1}{2}}\sup_{0\leq t\leq 1}|R^m_t|]=0$, then
the weak convergence of $\{(2^m)^{2H-\frac{1}{2}}J_tI^m_t\}$ 
immediately implies
the weak convergence of $\{(2^m)^{2H-\frac{1}{2}}(\Ymh_t-Y_t)\}$ to
the same limit distribution.

In this paper, in the case of fBm, for the four schemes mentioned above,
we prove that \\
$(2^m)^{2H-\frac{1}{2}+\vep}\sup_{0\leq t\leq 1}|R^m_t|$
converges to 0 almost surely and in $L^p$ for all $p\ge 1$.
Here $0<\vep<3H-1$ is an arbitrary constant.
This is one of our main theorems (Theorem~\ref{main theorem fbm}).
Our proof of this result does not rely on 
the weak convergence of $\{(2^m)^{2H-\frac{1}{2}}I^m_t\}$
but the uniform $L^p$ estimate of the H\"older norm
of $\{(2^m)^{2H-\frac{1}{2}}I^m_t\}$ independent of $m$.
Our result shows that the remainder term $R^m_t$
is really small compared to the term $J_tI^m_t$
and that it suffices to establish the limit theorem of weighted sum process of
Wiener chaos to obtain a limit theorem of the error of
$\Ymh_t-Y_t$ in certain cases.
In addition, we can give an estimate of the convergence rate of
$\sup_{0\le t\le 1}|\Ymh_t-Y_t|\to 0$ in $L^p$ sense (see Remark~\ref{rem849320482091}).
To the best of the authors' knowledge,
$L^p$ convergence rate does not appear in the literature
concerning fBm with the Hurst parameter $\frac{1}{3}<H<\frac{1}{2}$.

Our idea to obtain the estimate of $R^m_t$ is as follows.
The approximate solutions considered in this paper
are essentially defined at the discrete times $\Dm$.
We denote the solution and approximate solution at the discrete times
$\Dm$ by $\{Y_t\}_{t\in\Dm}$ and
$\{\Ymh_t\}_{t\in\Dm}$ respectively.
We note that all four schemes are given by similar recurrence relations.
More precisely, the recurrence relations of three schemes,
implementable Milstein, Crank-Nicolson and first-order Euler schemes,
can be obtained by adding extra two terms containing $d^m$
and $\hepm$ to the recurrence relation
of the Milstein scheme as we will see in \eqref{approximate solution}.
Based on this observation, we introduce 
an interpolation process $\{\Ymr_t\}_{t\in\Dm}$
which is parameterized by $\rho\in [0,1]$
and satisfies $Y^{m,0}_t=Y_t$ and $Y^{m,1}_t=\Ymh_t$ for all $t\in\Dm$.
Note that $\Ymr_t$ is different from 
the standard linear interpolation $(1-\rho)Y_t+\rho\Ymh_t$.
We define $\{\Ymr_t\}_{t\in\Dm}$ by \eqref{def_interplation}.
Let $\Zmr_t=\partial_{\rho}\Ymr_t$.
We can represent the process $\{\Zmr_t\}_{t\in\Dm}$ 
by a constant
variation method by using
a certain matrix valued process $\{\tJmr_t\}_{t\in\Dm}$
which approximates the derivative process $J_t$.
The important point is that all processes 
$\{(\Ymr_t, \Zmr_t, \tJmr_t, (\tJmr_t)^{-1})\}_{t\in\Dm}$
are solutions of certain discrete RDEs
and we can get good estimates of them.
We study the error process by the estimates and the expression
$\Ymh_t-Y_t=\int_0^1\Zmr_t d\rho$.
More precisely, we show that the main part of the right-hand side of this identity
is given by $J_tI^m_t$
and prove our main theorems.

We revisit Liu-Tindel's result \cite{liu-tindel}.
They also obtained an expression of 
$\Ymh_t-Y_t$ by using the process $\Phi^m_t$ which also approximates $J_t$.
See Lemma 6.4 in \cite{liu-tindel}.
Their proof for the convergence of $\{(2^m)^{2H-\frac{1}{2}}(\Ymh_t-Y_t)\}$ is based on the expression.
The process $\Phi^m_t$ is defined by using the standard linear interpolation process
$(1-\rho)Y_t+\rho\Ymh_t$ and
$\Phi^m_t$ is different from our $\tJmr_t$.
For the sake of conciseness of the paper, they did not get into the detailed study of
the integrability of $\Phi^m_t$ but they believed the integrability of 
it and its inverse.
Hence they could provide only the almost sure convergence rate of $\sup_{0\le t\le 1}|\Ymh_t-Y_t|\to 0$,
but not the $L^p$ convergence rate.
One may prove the integrabilities, but, we introduce different kind of
interpolation process $\Ymr_t$ and prove the integrability of $\tJmr$
to obtain our main results including the $L^p$ convergence rate.

We now explain how to implement our idea mentioned above.
In fact, Theorem~\ref{main theorem fbm} is deduced
from more general results (Theorem~\ref{main theorem} and Corollary~\ref{corollary for main theorem}).
As we already explained, the recurrence relations of the three schemes contain
extra terms containing $d^m$ and $\hepm$, which are not contained in the recurrence relation of the Milstein scheme.
Recall that the Milstein approximation solution
converges to the solution in pathwise sense in \cite{davie, friz-victoir}.
Hence we expect that if these extra terms are sufficiently small in a certain sense
then the approximate solutions
converge to the solution, not to mention the case of the four schemes.
In Theorem~\ref{main theorem}, 
we are concerned with such more general approximate solutions and general driving Gaussian processes
and provide estimates of the errors at discrete times $\Dm$.
More precisely, in such a setting,
we give the estimate of the remainder term $R^m_t$ $(t\in D_m)$
under Conditions~\ref{condition on R} 
and \ref{condition on dm}\,$\sim$\,\ref{condition on K}.
Condition~\ref{condition on R} is a natural condition on the covariance of the 
driving Gaussian process $B$ which ensures that $B$ can be lifted to 
a geometric rough path.
The other conditions 
are smallness conditions on $d^m$ and $\hepm$.
The main non-trivial condition among them
is Condition~\ref{minimal moment estimate} on $I^m$,
that is, the uniform estimate
of the $L^p$ norm of the H\"older norm of
$(2^m)^{2H-\frac{1}{2}}I^m$ independent of $m$.
In the case of the implementable Milstein, Milstein, and first-order
Euler schemes whose driving process is an fBm,
all conditions can be checked.
Hence, after establishing the continuous time version of Theorem~\ref{main theorem},
in Corollary~\ref{corollary for main theorem},
Theorem~\ref{main theorem fbm} for the three schemes follows from these results.
In the case of the Crank-Nicolson scheme, some of the conditions are not satisfied,
so Theorem~\ref{main theorem fbm} requires additional arguments to be established.
Here we mention how to show that
Conditions~\ref{condition on dm}\,$\sim$\,\ref{condition on K} are satisfied.
These conditions can be checked for the four schemes (as mentioned above, only partially, in the case of Crank-Nicolson scheme)
whose driving process is an fBm
by using the previously known results, \textit{e.g.}, in \cite{liu-tindel}.
We can also prove that these conditions hold
by a different idea based on the Malliavin calculus and estimates for multidimensional
Young integrals although we need more smoothness assumption on $\sigma$ and $b$ to prove
Condition~\ref{minimal moment estimate} than the previous study in \cite{liu-tindel}.
To make the paper reasonable size, 
we study these problems in a separate paper \cite{aida-naganuma2024LP}.

This paper is organized as follows.
In Section~\ref{statement}, 
we recall basic notions and estimates of rough path analysis
and the definition of the typical four schemes.
We next state our main theorems and make remarks on them.
After that we prove Theorem~\ref{main theorem fbm}
assuming Theorem~\ref{main theorem} and Corollary~\ref{corollary for main theorem}.
We close this section by introducing notion of small order nice discrete process
which includes the process of $d^m$ and $\hepm$ as examples.
The estimates of discrete Young integrals with respect to these processes
play an important role in this study.
In Section~\ref{interpolation}, 
we introduce processes $\{(\Ymr_t, \Zmr_t, \tJmr_t, (\tJmr_t)^{-1})\}$
and put the list of notations which we will use in this paper.
In Section~\ref{ymrk and Mmrk}, we give estimates for
$\{(\Ymr_t, \Zmr_t, \tJmr_t, (\tJmr_t)^{-1})\}$
by using Davie's argument in \cite{davie}.
We next give $L^p$ estimates for
$\tJmr_t$ and $(\tJmr_t)^{-1}$ by using the estimate of
Cass-Litterer-Lyons~\cite{cll}.
Thanks to this integrability, we can obtain good enough estimates 
of several quantities to prove our main theorems.
In Section~\ref{proof of main theorem},
we give a more precise estimate of $\{\Zmr_t\}$.
In the final part of this section, 
we give proofs of Theorem~\ref{main theorem} 
and Corollary~\ref{corollary for main theorem}.

\section{Main results, remarks, and preliminaries}\label{statement}
This section begins with a collection of the notation that will be used later.
Throughout this paper, $m$ denotes a positive integer.
Set $\Delta_m=2^{-m}$ and $\tmk=k 2^{-m}$ ($0\le k\le 2^m$)
and write $\Dm=\{\tmk\}_{k=0}^{2^m}$ for the dyadic partition of $[0,1]$. 
We identify the set of partition points and the partition.
The standard basis of $\RR^d$ is denoted by $\{e_{\alpha}\}_{\alpha=1}^d$
and $\floor{x}=\max\{n\in \mathbb{Z}~|~n\le x\}$ for $x\geq 0$.

Let us consider a process $F=\{F_t\}_{t\in I}$ for $I=[0,1]$ or $\Dm$.
We say that $F$ is a discrete process if $I=\Dm$, namely $F_t$ is evaluated at $t\in\Dm$.
We write $F_{s,t}=F_t-F_s$ for $s<t$ and,
for $0<\theta<1$, define the (discrete) $\theta$-H{\"o}lder norm by
\begin{align}\label{def Holder norm}
  \|F\|_\theta
   =\max_{s, t\in I, s<t}\frac{|F_{s,t}|}{|t-s|^\theta}.
\end{align}
For two-parameter functions $F=\{F_{s,t}\}_{s<t}$,
we define the $\theta$-H{\"o}lder norm in the same way.
In addition, the H\"older norm of $F$ on the interval $J\subset I$
is denoted by $\|F\|_{J,\theta}$.

When we are given a sequence of random variables $\{\eta_{\tmim,\tmi}\}_{i=1}^{2^m}$,
we define a discrete stochastic process $\{\eta_t\}_{t\in \Dm}$ 
and its increment process $\{\eta_{s,t}\}_{s\le t, s,t\in \Dm}$ by
\begin{align}
	\eta_t
	&=
		\sum_{i=1}^{2^mt}
			\eta_{\tmim,\tmi},\quad 
	&
	\eta_{s,t}
	&=
		\eta_t-\eta_s
	\label{eq89034802914832}
\end{align}
with the convention $\eta_0=0$.
In our study, such an $\{\eta_{\tmim,\tmi}\}$ 
arises as a small increment in the time interval
$[\tmim,\tmi]$.

The remainder of this section is structured as follows.
In Section~\ref{2.1}, we recall basic notion in rough path analysis and introduce
a condition (Condition~\ref{condition on R}) on the covariance of the driving 
Gaussian process $B$ under which $B$ can be lifted to a rough path.
We next introduce the small remainder term $\epm_{\tmkm,\tmk}$ of the solution.
In Section~\ref{2.2}, we explain four approximation schemes of 
RDE and introduce two important quantities $d^m_{\tmkm,\tmk}$ which belongs to
Wiener chaos of order 2 and $\hepm_{\tmkm,\tmk}$ 
which is defined as a small remainder term of approximate solution 
similarly to $\epm_{\tmkm,\tmk}$.
We next explain that the approximation equations can be written as 
common recurrence equations using $d^m_{\tmkm,\tmk}$ and $\hepm_{\tmkm,\tmk}$.
This observation is important for our study.
In Section~\ref{2.3}, 
taking the common recurrence equations into account,
we consider more general approximation equations.
We next introduce Conditions~\ref{condition on dm}\,$\sim$\,\ref{condition on K}
on $d^m$, $\hepm$ and iterated integrals of
$B$ and state our main theorems 
(Theorem~\ref{main theorem}, Corollary~\ref{corollary for main theorem}, and Theorem~\ref{main theorem fbm}).
In Section~\ref{sec proof of main thm fbm part 1}, we show Theorem~\ref{main theorem fbm}
in the case of the implementable Milstein, Crank-Nicolson, Milstein and first-order Euler schemes,
assuming Theorem~\ref{main theorem} and Corollary~\ref{corollary for main theorem}.
In Section~\ref{small order process}, we define a class of discrete
processes, small order nice discrete processes, which includes
$d^m, \epm, \hepm$.

\subsection{Rough paths and solutions to RDEs}\label{2.1}

Here we recall some basic notions of rough path analysis.
For details, see \cite{friz-victoir, friz-hairer, lq}.

Let $\frac{1}{3}<\theta\le \frac{1}{2}$.
Let $X=\{X_{s,t}\}_{0\leq s<t\leq 1}$ and $\XX=\{\XX_{s,t}\}_{0\leq s<t\leq 1}$
be two-parameter functions with values 
in $\RR^d$ and $\RR^d\otimes\RR^d$, respectively.
\begin{definition}
	\begin{enumerate}
		\item	We say that the pair $\RPX=(X,\XX)$ is a $\theta$-H\"older rough path
				if $\|X\|_{\theta}<\infty$, $\|\XX\|_{2\theta}<\infty$
				and $X_{s,t}=X_{s,u}+X_{u,t}$, $\XX_{s,t}=\XX_{s,u}+\XX_{u,t}+X_{s,u}\otimes X_{u,t}$
				for $0\leq s<u<t\leq 1$ (Chen's identity).
		\item	We say that a $\theta$-H{\"o}lder rough path $\RPX=(X,\XX)$
				is geometric if it satisfies the following:
				there exists a sequence of smooth paths $X^m$ 
				such that its natural lift $\RPX^m=(X^m,\XX^m)$,
				where $\XX^m_{s,t}=\int_s^t X^m_{s,u}\otimes dX^m_{0,u}$, approximates $\RPX=(X,\XX)$
				in the rough path metric, that is,
				\begin{align*}
					\lim_{m\to\infty}
						\{\|X-X^m\|_{\theta}+\|\XX-\XX^m\|_{2\theta}\}
					=
						0.
				\end{align*}
				We denote by $\mathscr{C}^\theta_g$ the set of all $\theta$-H{\"o}lder geometric rough paths.
	\end{enumerate}
\end{definition}
We denote by $X^\alpha_{s,t}$ the $e_\alpha$-component of $X_{s,t}$
and by $X^{\alpha,\beta}_{s,t}$ the $e_\alpha\otimes e_\beta$-component 
of $\XX_{s,t}$.
Namely we write $X_{s,t}=\sum_{\alpha=1}^d X^\alpha_{s,t}e_{\alpha}$
and $\XX_{s,t}=\sum_{1\le \alpha,\beta\le d}X^{\alpha,\beta}_{s,t}e_{\alpha}
\otimes e_{\beta}$.
Recall that we can construct the third level rough paths
from the first and second level rough paths.
The $e_{\alpha}\otimes e_{\beta}\otimes e_{\gamma}$-component
of the third level rough paths will be denoted by $X^{\alpha,\beta,\gamma}_{s,t}$.

Next we introduce the notion of controlled paths and integration of controlled paths.
\begin{definition}
	Let $X=\{X_t\}_{0\leq t\leq 1}$ be a $\theta$-H{\"o}lder function with values in $\RR^d$.
	A $\theta$-H{\"o}lder function $Z=\{Z_t\}_{0\leq t\leq 1}$ with values in $\RR^K$
	is said to be a path controlled by $X$
	if there exist a $\theta$-H{\"o}lder function $Z'=\{Z'_t\}_{0\leq t\leq 1}$ with valued in $\mathcal{L}(\RR^d, \RR^K)$
	and a $(2\theta)$-H{\"o}lder function $R=\{R_{s,t}\}_{0\leq s<t\leq 1}$
	satisfying
	$
		Z_t-Z_s
		=
			Z'_s (X_t-X_s)+R_{s,t}
	$
	$(0\leq s<t\leq 1)$.
	The set of all pairs $(Z,Z')$ is denoted by $\mathscr{D}_X^{2\theta}([0,1],\RR^K)$.
\end{definition}
Let $\RPX=(X,\XX)$ be a geometric $\theta$-H\"older rough path
and identify $X$ with a one-parameter function by $X_t=X_{0,t}$.
We can define an integration of a path $(Z,Z')$ controlled by $X$ against $\RPX=(X,\XX)$
as follows.
\begin{theorem}[{\cite[Theorem~4.10]{friz-hairer}}]
	Let $(Z,Z')\in \mathscr{D}_X^{2\theta}([0,1],\mathcal{L}(\RR^d, \RR^K))$.
	We can define an integration of $(Z,Z')$ along $\RPX=(X,\XX)$ by
	\begin{align*}
		\int_s^t
			Z_u
			d\RPX_u
		=
			\lim_{|\mathcal{P}|\to 0}
				\sum_{i=1}^M
					\{
						Z_{t_{i-1}} X_{t_{i-1},t_i}
						+Z'_{t_{i-1}} \XX_{t_{i-1},t_i}
					\}.
	\end{align*}
	Here $\mathcal{P}=\{t_i\}_{i=0}^M$ denotes a partition of the interval $[s,t]$
	and $|\mathcal{P}|=\max\{t_i-t_{i-1}|1\leq i\leq M\}$.
	We call the left-hand side a rough integral.
\end{theorem}
Let $P_{\alpha}$ be the projection operator on $\RR^d$ onto
the subspace spanned by $e_{\alpha}$.
Then $\int_0^tZ_ud\RPX_u=\sum_{\alpha=1}^d\int_0^tZ_uP_{\alpha}d\RPX_u$ holds.
We may write $\int_0^tZ_uP_{\alpha}d\RPX_u=\int_0^tZ_ue_{\alpha}dX^{\alpha}_u$.
Actually, the rough integral $\int_0^t \tilde{Z}_ud\bar{Z}_u$ can be defined for any paths
$\tilde{Z}_t, \bar{Z}_t$
controlled by $X$ (see \cite[Remark~4.12]{friz-hairer}).
Also note that $Z_te_{\alpha}$ and $X^{\alpha}_t$ are $\theta$-H\"older
paths controlled by $X$.
It is easy to check that
$\int_0^tZ_ue_{\alpha}dX^{\alpha}_u$ coincide with the rough integral in that sense.
Note that 
the process 
$
	\big\{
		\big(
			\int_0^t
				Z_u
				d\RPX_u,
			Z_t
		\big)
	\big\}_{0\leq t\leq 1}
$
is also a path controlled by $X$
and we can define iterated integrals in the sense of rough integrals.
Furthermore, we have the following formula: for any $f\in C^3_b(\RR^K,\RR^L)$,
\begin{align}
	\label{eq483290842}
	f(Z_t)-f(Z_s)
	=
		\int_0^t (Df)(Z_u)Z'_ud\RPX_u
		+\int_0^t (Df)(Z_u)d\Gamma_u
\end{align}
if $(Z,Z')\in \mathscr{D}_X^{2\theta}([0,1],\RR^K)$ satisfies $Z_t=Z_0+\int_s^t Z'_u d\RPX_u+\Gamma_t$ 
for some $(Z',Z'')\in \mathscr{D}_X^{2\theta}([0,1],\mathcal{L}(\RR^d, \RR^K))$
and smooth function $\Gamma$ with values in $\RR^K$.
For detail, see \cite[Theorem~7.7]{friz-hairer}.

Next we introduce the notion of solutions to RDEs.
Let $\xi\in\RR^n$,
$\sigma\in C^4_b(\RR^{n}, \mathcal{L}(\RR^d, \RR^{n}))$,
$b\in C^2_b(\RR^{n},\RR^{n})$
and consider an RDE driven by $X$ on $\RR^{n}$,
\begin{align}
	Y_t
	=
		\xi
		+\int_0^t\sigma(Y_s)d\RPX_s
		+\int_0^t b(Y_s)ds,
		\qquad 0\le t\le 1.
	\label{generalRDE}
\end{align}
Here the first integral should be understood as a rough integral.
We also write $Y_t(\xi,X)=Y_t$ if the solution $Y_t$ exists.
We have several notion of solution, which are equivalent.
To state them, we set
\begin{align}
  ((D\sigma)[\sigma])(y)[v\otimes w]
  =
    D\sigma(y)[\sigma(y)v]w,
	\qquad
	y\in\RR^n, v,w\in\RR^d.
  \label{eqSimplifiedNotationDef}
\end{align}
In this notation, we have
\begin{align}
	\label{eqSimplifiedNotation02}
	((D\sigma)[\sigma])(y)\XX_{s,t}
	=
		\sum_{\alpha,\beta=1}^d
		  (D\sigma)(y)[\sigma(y)e_{\alpha}]e_{\beta}
		  X^{\alpha,\beta}_{s,t}.
\end{align}

\begin{theorem}[{\cite[Theorem~8.3 and Proposition~8.10]{friz-hairer}}]
	The following are equivalent and both are valid.
	\begin{enumerate}
		\item	There exists a unique $(Y,Y')\in \mathscr{D}_X^{2\theta}([0,1],\RR^n)$
				satisfying \eqref{generalRDE} with $Y'=\sigma(Y)$.
		\item	There exists a unique process $Y\colon[0,1]\to\RR^n$ satisfying 
				\begin{align}
					|Y_t-Y_s-\sigma(Y_s)X_{s,t}-((D\sigma)[\sigma])(Y_s)\XX_{s,t}-b(Y_s)(t-s)|
					\leq
						C(t-s)^{3\theta}
					\label{definition of davie}
				\end{align}
				for $0\leq s<t\leq 1$.
				Here $C$ can be estimated by a polynomial function of
				$\|X\|_{[0,1],\theta}$ and $\|\XX\|_{[0,1],2\theta}$.
				This is called a solution in the sense of Davie~\cite{davie}.
	\end{enumerate}
\end{theorem}
Note that we can choose $C$ in \eqref{definition of davie} so that 
it can be estimated by a polynomial function of
$\|X\|_{[s,t],\theta}$ and $\|\XX\|_{[s,t],2\theta}$.
We will record this estimate in Lemma~\ref{lem8345901809321} later.
Although the estimate on $C$ in (\ref{definition of davie}) and 
the unique existence of solution hold under weaker assumption
that $\sigma\in C^3_b$ and $b\in C^1_b$
(see \cite{friz-hairer}), we need to assume the above condition on
$\sigma$ and $b$ in our study.

We now introduce a condition to construct a rough path associated to 
a Gaussian process under which we will work.
Let $\Omega=C_0([0,1],\RR^d)$ be the set of $\RR^d$-valued continuous functions
on $[0,1]$ starting at the origin,
$B$ be the canonical process on $\Omega$, that is, $B_t(\omega)=\omega(t)$ 
($\omega\in \Omega$),
and $\mu$ be a centered Gaussian probability measure on $\Omega$.
Throughout this paper, we put the next condition on $B$:
\begin{condition}\label{condition on R}
Let $\frac{1}{3}<H\le \frac{1}{2}$.
	Let $B^{\alpha}_t$ be the $\alpha$-th component of 
$B_t$ $(1\le \alpha\le d)$.
Then $B^1_t,\ldots, B^d_t$ are independent centered continuous Gaussian processes.
	Let $R^{\alpha}(s,t)=E[B^\alpha_sB^\alpha_t]$.
	Then $V_{(2H)^{-1}}(R^\alpha ; [s,t]^2 )\le C_\alpha|t-s|^{2H}$
	holds for all $1\le \alpha\le d$ and $0\le s<t\le 1$.
	Here $V_p(R^{\alpha} ; [s,t]^2)$ denotes the $p$-variation norm of 
$R^{\alpha}$ on $[s,t]^2$.
\end{condition}

Note that Condition~\ref{condition on R} holds for the fBm
	with the Hurst parameter
	$\frac{1}{3}<H\le \frac{1}{2}$.

\begin{remark}\label{remark on lift}
It is known that under Condition~\ref{condition on R},
$B$ can be naturally lifted to $\RPB=(B,\BB)\in \mathscr{C}^\theta_g$ for any $\frac{1}{3}<\theta<H$.
More precisely, we can prove the following property
(Remark 10.7 in \cite{friz-hairer}, Theorem 15.33 in \cite{friz-victoir}).
We consider a sequence of smooth rough path $\RPB^m(\omega)=(B^m(\omega),\BB^m(\omega))$ defined by
a piecewise linear approximation of $B(\omega)$ such that
$\lim_{m\to\infty}\max_{0\le t\le 1}|B^m_t(\omega)-B_t(\omega)|=0$ for all 
$\omega\in \Omega$.
Then $\RPB^m(\omega)=(B^m(\omega),\BB^m(\omega))\in \mathscr{C}^\theta_g$ converges in probability in the 
$\theta$-H\"older rough path
metric for any $\frac{1}{3}<\theta<H$.
This implies that there exists a subset
$\Omega_0$ with $\mu(\Omega_0)=1$ such that, if necessary choosing a subsequence,
the limit
$\RPB(\omega)=(B(\omega),\BB(\omega))$ belongs to $\mathscr{C}^\theta_g$
for any $\omega\in \Omega_0$ and any $\frac{1}{3}<\theta<H$.
Of course, this rough path depends on the selected versions, but, 
note that any versions are almost surely identical.
We consider solutions to RDEs driven by this rough path obtained by Gaussian process
satisfying Condition~\ref{condition on R}.
\end{remark}

Here we fix $\frac{1}{3}<H^-<H$. 
For later use, we introduce a random variable $C(B)$ by 
	\begin{align}
	  C(B)
	  &=
		\max
			\left\{
				\|B(\omega)\|_{H^-}, 
				\|\BB(\omega)\|_{2H^-}
			\right\},
			\qquad
			\omega\in \Omega_0,\label{CB}
	\end{align}
	and a subset $\Omega^{(m)}_0$ of $\Omega_0$ by
	\begin{align*}
		\Omega^{(m)}_0
		=
		  \bigg\{
			  \omega\in \Omega_0
			  ~\bigg|~
			  \sup_{|t-s|\le 2^{-m}}
				\left|\frac{B_{s,t}(\omega)}{(t-s)^{H^{-}}}\right|
			  \le
				\frac{1}{2},\quad
			  \sup_{|t-s|\le 2^{-m}}
				\left|\frac{\BB_{s,t}(\omega)}{(t-s)^{2H^{-}}}\right|
			  \le
				\frac{1}{2}
		  \bigg\}.
	\end{align*}
	Under Condition~\ref{condition on R},
	$C(B)\in \cap_{p\ge 1}L^p$ holds.
	We refer the readers for this to 
	\cite{friz-victoir2010, friz-victoir, friz-hairer}.
Therefore, under Condition~\ref{condition on R}, we see that
\begin{align}
 \mu((\Omega^{(m)}_0)^\complement)\le C_{p}
2^{-mp}\quad\quad \text{for any $p>1$}
\label{complement of Omega0}
\end{align}
which eventually implies that the complement set is negligible for our problem.
Below, we actually consider analogous subset 
$\Omega^{(m, d^m)}_0$ which will be introduced 
in Section~\ref{small order process}.
The proof of $(\ref{complement of Omega0})$ is as follows.
Let $\kappa>0$ be a positive number satisfying 
$H^-+\kappa<H$.
Let $C(B)_{H^-+\kappa}$ denote the number 
obtained by replacing $H^-$ by $H^-+\kappa$ in the definition 
(\ref{CB}).
Then we have
\begin{align*}
 		  \sup_{|t-s|\le 2^{-m}}
			\left|\frac{B_{s,t}(\omega)}{(t-s)^{H^{-}}}\right|+
\sup_{|t-s|\le 2^{-m}}
			\left|\frac{\BB_{s,t}(\omega)}{(t-s)^{2H^{-}}}\right|
&\le 2^{1-m\kappa}C(B)_{H^-+\kappa}.
\end{align*}
Hence we obtain
$\liminf_{m\to\infty}\Omega^{(m)}_0=\Omega_0$ and
\[
 \mu((\Omega^{(m)}_0)^\complement)\le \mu(C(B)_{H^-+\kappa}\ge 2^{m\kappa-2})
\le 2^{-p(m\kappa-2)}\|C(B)_{H^-+\kappa}\|_{L^p}^p,
\]
which is the desired result.

\begin{remark}[About the constants in the estimates]
When a positive constant $C$ can be written as a polynomial function of
the sup-norm of some functions $\sigma, b, c$ and their derivatives,
we may say $C$ depends on $\sigma, b, c$ polynomially.
Similarly, 
when a constant $C$ can be written as a polynomial
of some positive random variable $X$, 
the sup-norms of $\sigma, b, c$ and their derivatives,
we say that $C$ depends on $\sigma, b, c, X$ polynomially.
Of course the coefficients of the polynomial should not depend on $\omega$.
When $X=C(B)$, we may denote such a constant $C$ by $\tilde{C}(B)$.
\end{remark}

Throughout this paper, we assume $B$ satisfies Condition~\ref{condition on R}
and $\RPB=(B,\BB)$ is the canonically defined rough path as explained above.
Let $Y_t=Y_t(\xi,B)$ be the solution to RDE on $\RR^n$ driven by $B$:
\begin{align}
	Y_t(\xi,B)
	=
		\xi
		+\int_0^t\sigma(Y_s(\xi,B))d\RPB_s
		+\int_0^t b(Y_s(\xi,B))ds,
		\qquad 0\le t\le 1.\label{rde}
\end{align}
We may omit writing the starting point $\xi$ 
and the driving process $B$ in $Y_t(\xi,B)$.
Note that $J_t=\partial_{\xi}Y_t(\xi)\in \mathcal{L}(\RR^n)$ 
and its inverse $J^{-1}_t$ are the solutions to the following RDEs:
\begin{align}
  \label{Jacobian}
 J_t&=I+\int_0^t(D\sigma)(Y_u)[J_u]d\RPB_u+\int_0^t(Db)(Y_u)[J_u]du,\\
 \label{Jaconian inverse}
J^{-1}_t&=I-\int_0^tJ^{-1}_u(D\sigma)(Y_u)d\RPB_u-\int_0^tJ^{-1}_u(Db)(Y_u)du.
\end{align}

We conclude this section by presenting a lemma and making a remark.
For every $1\leq k\leq 2^m$, 
define $\ep^m_{\tmkm,t}(\xi)$ $(\tmkm\le t\le \tmk)$ by
\begin{align}
Y_t&=Y_{\tmkm}+\sigma(Y_{\tmkm})B_{\tmkm,t}
    +((D\sigma)[\sigma])(Y_{\tmkm})\BB_{\tmkm,t}
    +b(Y_{\tmkm})(t-\tmkm)
  +\ep^m_{\tmkm,t}(\xi).
	\label{definition of epm}
\end{align}
We may use the notation $\epm_{\tmkm,t}$ instead of $\epm_{\tmkm,t}(\xi)$
for simplicity.
As we explained in the inequality (\ref{definition of davie}), we have the following.
\begin{lemma}\label{lem8345901809321}
\begin{enumerate}
 \item There exists a constant $C>0$ such that
	\begin{align}
		\label{eq459301490123}
	  |\epm_{\tmkm,t}|
	\le C (t-\tmkm)^{3H^-}
	\qquad 
	\text{for\,\,\,all}
	\qquad 
	1\le k\le 2^m,\quad \omega\in \Omega_0.
	\end{align}
	Here $C$ depends on $\|B\|_{[\tmkm,t],H^-}$, $\|\BB\|_{[\tmkm,t],2H^-}$,
	$\sigma, b$ polynomially.
\item There exists a constant $C>0$ depending on $\sigma,b$ polynomially
and bounded Lipschitz continuous
functions $F_{\alpha,\beta,\gamma}$, 
$F^1_{\alpha}$, $F^{2}_{\alpha}$ from
$\RR^n$ to $\RR^n$ such that for all $1\le k\le 2^m$
and $\tmkm\le t\le \tmk$,
\begin{multline}
	\left|
		\epm_{\tmkm,t}
		-
		\sum_{\alpha,\beta,\gamma}
			F_{\alpha,\beta,\gamma}(Y_{\tmkm})
			B^{\alpha,\beta,\gamma}_{\tmkm,t}
		-
		\sum_{\alpha}
			F^1_{\alpha}(Y_{\tmkm})
			B^{0,\alpha}_{\tmkm,t}
		-
		\sum_{\alpha}
			F^2_{\alpha}(Y_{\tmkm})
			B^{\alpha,0}_{\tmkm,t}
	\right|\\
	\le C (t-\tmkm)^{4H^-},
	\qquad 
	\qquad 
	\omega\in \Omega^{(m)}_0,\label{eq8940289042}
\end{multline}
where 
\begin{align}
 B^{0,\alpha}_{\tmkm,t}=\int_{\tmkm}^t(s-\tmkm)dB^{\alpha}_s,\quad\quad
B^{\alpha,0}_{\tmkm,t}=\int_{\tmkm}^tB^{\alpha}_{\tmkm,s}ds.\label{B0alpha}
\end{align}
\end{enumerate}
\end{lemma}
\begin{proof}
	We need only to prove \eqref{eq8940289042}.
	First we give an expression of $\ep^m_{\tmkm,t}$.
	Note that the solution $Y_t$ to \eqref{rde} satisfies
	$(Y,\sigma(Y)) \in \mathscr{D}_X^{2\theta}([0,1],\RR^n)$ 
	and $(\sigma(Y),((D\sigma)[\sigma])(Y))\in \mathscr{D}_X^{2\theta}([0,1],\mathcal{L}(\RR^d, \RR^n))$.
	Hence we can use \eqref{eq483290842}.
Then by applying the formula to $f(Y_t)-(Y_{\tmkm})$ for $f\in C^3_b(\RR^n,\RR^L)$
successively, we can decompose
$\epsilon^m_{\tmkm,t}$ in the following way.
This calculation is possible because $\sigma\in C^4_b, b\in C^2_b$.
We need the following functions to state it:
\begin{gather*}
  \begin{aligned}
      F^0(y)&=(Db)(y)[b(y)], & 
      F^1_{\alpha}(y)&=(Db)(y)[\sigma(y)e_{\alpha}], &
      F^2_{\alpha}(y)&=(D\sigma(y)e_{\alpha})[b(y)],
  \end{aligned}\\
  \begin{aligned}
    F_{\alpha,\beta,\gamma}(y)
    &=
      D\Bigl\{(D\sigma(y)e_{\gamma})[\sigma(y)e_{\beta}]\Bigr\}
      [\sigma(y)e_{\alpha}], &
    G_{\alpha,\beta}(y)
    &=
      D\Bigl\{(D\sigma(y)e_{\beta})[\sigma(y)e_{\alpha}]\Bigr\}[b(y)].
  \end{aligned}
\end{gather*}
The decomposition formula is as follows,
\begin{align}
&\ep^m_{\tmkm,t}
=
\sum_{\alpha,\beta,\gamma}
\int_{\tmkm}^{t}\left\{\int_{\tmkm}^{s}\left(\int_{\tmkm}^{u}
F_{\alpha,\beta,\gamma}(Y_v)
dB^{\alpha}_v\right)dB^{\beta}_u\right\}
dB^{\gamma}_s\nonumber\\
&\quad
+\sum_{\alpha,\beta,\gamma}
\int_{\tmkm}^{t}\left\{\int_{\tmkm}^{s}\left(\int_{\tmkm}^{u}
G_{\alpha,\beta}(Y_v)
dv\right)dB^{\alpha}_u\right\}
dB^{\beta}_s
+\int_{\tmkm}^{t}\left(\int_{\tmkm}^sF^0(Y_u)du\right)ds\nonumber\\
&\quad +\sum_{\alpha}\int_{\tmkm}^{t}\left(
\int_{\tmkm}^sF^1_{\alpha}(Y_u)du\right)dB^{\alpha}_s
+\sum_{\alpha}\int_{\tmkm}^{t}\left(
\int_{\tmkm}^sF^2_{\alpha}(Y_u)dB^{\alpha}_u\right)ds\nonumber\\
&:=I_1+\cdots+I_5.
\label{explicit form of epm}
\end{align}

By using estimates of rough integrals,
we have the following estimates: for all $\omega\in\Omega^{(m)}_0$, it holds that
\begin{gather}
	\Bigl|I_1-\sum_{\alpha,\beta,\gamma}F_{\alpha,\beta,\gamma}(Y_{\tmkm})
	B^{\alpha,\beta,\gamma}_{\tmkm,t}\Bigr|
	\le C(t-\tmkm)^{4H^-},\\
	\left|I_2\right|\le C(t-\tmkm)^{1+2H^-}, 
	\qquad
	\qquad
	\quad
	\left|I_3\right|\le C(t-\tmkm)^2, \\
	\left|I_4-F^1_{\alpha}(Y_{\tmkm})B^{0,\alpha}_{\tmkm,t}\right|
+\left|I_5-F^2_{\alpha}(Y_{\tmkm})B^{\alpha,0}_{\tmkm,t}\right|\le C(t-\tmkm)^{1+2H^-},\label{I4 I5}
\end{gather}
where $C$ depends on $\sigma$ and $b$ polynomially.
	This completes the proof.
\end{proof}

\begin{remark}\label{rem8349189021412}
For every $s,t\in \Dm$ with $s\le t$, define $\epm_t$ and $\epm_{s,t}$
in the same way as \eqref{eq89034802914832} with 
$\eta_{\tmim,\tmi}=\epm_{\tmim,\tmi}$.
Note that the identity
$
  \epm_{s,t}
  =
    Y_t-Y_s-\sigma(Y_s)B_{s,t}
    -((D\sigma)[\sigma])(Y_s)\BB_{s,t}
    -b(Y_s)(t-s)
$
does not hold
for general $s,t\in \Dm$ with $s\le t$.
\end{remark}

\subsection{Four approximation schemes}\label{2.2}
In this section, we introduce typical four approximation schemes.
That is, we introduce
the implementable Milstein approximate solution $Y^{\rIM,m}_t$, 
the Milstein approximate solution $Y^{\rM,m}_t$,
the first-order Euler approximate solution $Y^{\rFE,m}$,
and the Crank-Nicolson approximate solution $Y^{\rCN,m}_t$
associated to the dyadic partition $\Dm$.
The first three schemes are explicit scheme and defined inductively as follows:
$Y^{\rIM,m}_0=Y^{\rM,m}_0=Y^{\rFE,m}_0=\xi$ and 
\begin{align*}
  Y^{\rIM,m}_t
  &=
    Y^{\rIM,m}_{\tmkm}
    +\sigma(Y^{\rIM,m}_{\tmkm})B_{\tmkm,t}
    +((D\sigma)[\sigma])(Y^{\rIM,m}_{\tmkm})
      \left[
        \frac{1}{2}
        B_{\tmkm,t}
        \otimes
        B_{\tmkm,t}
      \right]\\
  &\qquad\qquad\qquad\qquad\qquad\qquad\qquad
    +b(Y^{\rIM,m}_{\tmkm})(t-\tmkm),\\
  Y^{\rM,m}_t
  &=
    Y^{\rM,m}_{\tmkm}
    +\sigma(Y^{\rM,m}_{\tmkm})B_{\tmkm,t}
    +((D\sigma)[\sigma])(Y^{\rM,m}_{\tmkm})\BB_{\tmkm,t}
    +b(Y^{\rM,m}_{\tmkm})(t-\tmkm),\\
  Y^{\rFE,m}_{t}
  &=
    Y^{\rFE,m}_{\tmkm}
    +\sigma(Y^{\rFE,m}_{\tmkm})B_{\tmkm,t}
    +((D\sigma)[\sigma])(Y^{\rFE,m}_{\tmkm})
      \left[
        \frac{1}{2}
        \sum_{\alpha=1}^d
          e_{\alpha}
          \otimes
          e_{\alpha}
          E[(B_{\tmkm,t}^\alpha)^2]
      \right]\\
  &\qquad\qquad\qquad\qquad\qquad\qquad\qquad
  +b(Y^{\rFE,m}_{\tmkm})(t-\tmkm),
\end{align*}
for every $\tmkm< t\le\tmk$ and $1\le k\le 2^m$.
In the above, we omit writing the initial value $\xi$ for the solution.
With the notation \eqref{eqSimplifiedNotationDef},
we have
\begin{align}
  ((D\sigma)[\sigma])(y)
    \left[
      \frac{1}{2}
      B_{s,t}
      \otimes
      B_{s,t}
    \right]
  &=
    \sum_{\alpha,\beta=1}^d
      \frac{1}{2}(D\sigma)(y)[\sigma(y)e_{\alpha}]e_{\beta}
      B^{\alpha}_{s,t}
      B^{\beta}_{s,t},\\
  \label{eqSimplifiedNotation03}
  ((D\sigma)[\sigma])(y)
  \left[
    \frac{1}{2}
    \sum_{\alpha=1}^d
      e_{\alpha}
      \otimes
      e_{\alpha}
      E[(B_{s,t}^\alpha)^2]
 \right]
 &=
 \sum_{\alpha=1}^d
 \frac{1}{2}
  (D\sigma)(y)
  [\sigma(y)e_{\alpha}]e_{\alpha}E[(B_{s,t}^\alpha)^2].
\end{align}

Next we introduce the Crank-Nicolson scheme.
Since the Crank-Nicolson scheme is an implicit scheme
and an equation stated later with respect to $Y^{\rCN,m}_t$ must be solvable.
For that purpose, we already introduced the set $\Omega^{(m)}_0$.
Since $D\sigma$ and $Db$ are bounded function,
the mapping
\[
  v
  \mapsto
    \eta
    +\frac{1}{2}\left(\sigma(\eta)+\sigma(v)\right)B_{\tmkm,t}
    +\frac{1}{2}\left(b(\eta)+b(v)\right)(t-\tmkm),
    \qquad
    \tmkm\leq t\leq \tmk,
\]
is a contraction mapping for any $\eta\in\RR^n$ and $\omega\in \Omega^{(m)}_0$
for large $m$.
Therefore, for $\omega\in \Omega^{(m)}_0$ for large $m$, 
the Crank-Nicolson scheme $Y^{\rCN,m}_t$ is uniquely defined
as the following inductive equation:
$Y^{\rCN,m}_0=\xi$ and
\begin{align}
	Y^{\rCN,m}_t
	&=
		Y^{\rCN,m}_{\tmkm}
		+\frac{1}{2}\left(\sigma(Y^{\rCN,m}_{\tmkm})+\sigma(Y^{\rCN,m}_t)\right)B_{\tmkm,t} \notag \\
	  &\qquad\qquad\qquad\qquad\qquad
	  +\frac{1}{2}\left(b(Y^{\rCN,m}_{\tmkm})+b(Y^{\rCN,m}_t)\right)(t-\tmkm)
\label{CN equation}
\end{align}
for every $\tmkm< t\le\tmk$ and $1\le k\le 2^m$.
For the completeness of definition, 
we set $Y^{\rCN,m}_t\equiv \xi$
for $\omega\in \Omega_0\setminus \Omega^{(m)}_0$.

In what follows, we discuss how to address the four schemes collectively.
This is one of the key ingredients of this paper.
We use the common notation $\{\hat{Y}^m_t\}_{t\in [0,1]}$ to 
denote these four approximate solutions.
The four approximate solutions $\{\hat{Y}^m_t\}_{t\in [0,1]}$
also satisfy similar but a little bit different equations to
(\ref{definition of epm}).
Indeed, by choosing a function 
$c\in C^3_b(\RR^n, L(\RR^d\otimes \RR^d,\RR^n))$ 
and random variables 
$d^m=\{d^m_{\tmkm,t}\}_{1\leq k\leq 2^m,\tmkm<t\le \tmk}\subset \RR^d\otimes \RR^d$ 
and $\hepm(\xi)=\{\hepm_{\tmkm,t}(\xi)\}_{1\leq k\leq 2^m,\tmkm<t\le \tmk}\subset \RR^n$
defined on $\Omega_0$,
these approximate equations can be written as the 
following common form on $\Omega_0$: $\Ymh_0=\xi$ and
\begin{align}
  \Ymh_t
  &=
    \Ymh_{\tmkm}+\sigma(\Ymh_{\tmkm})B_{\tmkm,t}
    +((D\sigma)[\sigma])(\Ymh_{\tmk})\BB_{\tmkm,t}
    +b(\Ymh_{\tmkm})(t-\tmkm)\nonumber\\
  &\qquad\qquad\qquad\qquad
    +c(\Ymh_{\tmkm})d^m_{\tmkm,t}
    +\hat{\ep}^m_{\tmkm,t}(\xi), \quad \tmkm<t\le \tmk.\label{approximate solution}
\end{align}

We explain more precisely what $c, d^m, \hepm(\xi)$ are for all cases.
In all cases, $c$ is given by
\begin{align*}
	c(y)[v\otimes w]
	=
		((D\sigma)[\sigma])(y)[v\otimes w]
	=
		D\sigma(y)[\sigma(y)v]w,
		\qquad
		y\in\RR^n, v,w\in\RR^d.
\end{align*}
and $d^m_{\tmkm,t}$ arises from the difference between the second level rough paths and
their approximations in each scheme.
Furthermore, $\hepm_{\tmkm,t}(\xi)$ denotes a smaller term in each scheme.
We may use the notation $\hepm_{\tmkm,t}$ for $\hepm_{\tmkm,t}(\xi)$ 
if there is no confusion.
For $Y^{\rIM,m}$, $Y^{\rM,m}$ and $Y^{\rFE,m}$,
the pairs of $d^m$ and $\hepm$ are given by
\begin{align*}
	d^{\rIM,m}_{\tmkm,t}
	&=
		\frac{1}{2}
		B_{\tmkm,t}
		\otimes
		B_{\tmkm,t}
		-
		\BB_{\tmkm,t},
	&
	\hat{\epsilon}^{\rIM,m}_{\tmkm,t}
	&=
		0,\\
	d^{\rM,m}_{\tmkm,t}
	&=
		0,
	&
	\hat{\epsilon}^{\rM,m}_{\tmkm,t}
	&=
		0,\\
	d^{\rFE,m}_{\tmkm,t}
	&=
	  \frac{1}{2}
	  \sum_{\alpha=1}^d
		e_{\alpha}
		\otimes
		e_{\alpha}
		E[(B^\alpha_{\tmkm,t})^2]
	  -
	  \BB_{\tmkm,t},
  	&
  	\hat{\epsilon}^{\rFE,m}_{\tmkm,t}
  	&=
		0.
\end{align*}

The Crank-Nicolson scheme leads to a slightly complicated situation.
For the Crank-Nicolson scheme $Y^{\rCN,m}$,
we set
$
	d^{\rCN,m}_{\tmkm,t}
	=
		\frac{1}{2}
		B_{\tmkm,t}
		\otimes
		B_{\tmkm,t}
		-
		\BB_{\tmkm,t}
$,
that is, the same one as the case of implementable Milstein scheme.
Once $d^{\rCN,m}_{\tmkm,t}$ is defined,
$\hat{\ep}^{\rCN,m}_{\tmkm,t}$ is automatically determined by 
the identity \eqref{approximate solution}.
For $\omega\in\Omega_0\setminus \Omega_0^{(m)}$, 
from $Y^{\rCN,m}_t=\xi$ and \eqref{approximate solution},
we easily see
\begin{align}
  \hat{\ep}^{\rCN,m}_{\tmkm,t}
  = 
  	-\sigma(\xi)B_{\tmkm,t}
    -((D\sigma)[\sigma])(\xi)\BB_{\tmkm,t}
    -c(\xi)d^{\rCN,m}_{\tmkm,t}
    -b(\xi)(t-\tmkm).
	\label{definition of hepm for CN}
\end{align}
For $\omega\in\Omega_0^{(m)}$, we set
\begin{align}
  \hat{\ep}^{\rCN,m}_{\tmkm,t}
  &=
    \frac{1}{2}
    \bigg(
        \int_0^1
          \bigg(
           (D\sigma)(Y^{\rCN,m}_{\tmkm}+\theta Y^{\rCN,m}_{\tmkm,t})[Y^{\rCN,m}_{\tmkm,t}]\nonumber\\
  &\qquad\qquad\qquad
          -(D\sigma)(Y^{\rCN,m}_{\tmkm})[\sigma(Y^{\rCN,m}_{\tmkm})B_{\tmkm,t}]\bigg)d\theta
    \bigg)
    B_{\tmkm,t}\nonumber\\
  &\qquad\qquad
    +
    \frac{1}{2}
    \left(
      \int_0^1(Db)
        (Y^{\rCN,m}_{\tmkm}+\theta Y^{\rCN,m}_{\tmkm,t})[Y^{\rCN,m}_{\tmkm,t}]
        d\theta
    \right)
	(t-\tmkm).
 \label{explicit form of hepmCN}
\end{align}
Then we see that the recurrence relation \eqref{approximate solution} holds
and that $\hat{\ep}^{\rCN,m}_{\tmkm,t}$ admits good estimates
as follows.
\begin{lemma}\label{lem89430222}
	Let $\omega\in \Omega_0^{(m)}$.
	\begin{enumerate}
		\item	The Crank-Nicolson approximate solution satisfies
				\eqref{approximate solution}
				with $d^m_{\tmkm,t}=d^{\rCN,m}_{\tmkm,t}$
				and $\hepm_{\tmkm,t}=\hat{\ep}^{\rCN,m}_{\tmkm,t}$.
		\item	There exists a positive constant $C$ such that
				\begin{align*}
					|\hat{\ep}^{\rCN,m}_{\tmkm,t}|
					\le
						C|t-\tmkm|^{3H^-}
					\qquad 
					\text{for all}
					\qquad
					1\le k\le 2^m.
				\end{align*}
				Here, $C$ depends on $\sigma$ and $b$ polynomially.			
		\item	There exist bounded Lipschitz continuous functions 
				$\varphi_{\alpha,\beta,\gamma} :\RR^n\to\RR^n$
				and $\psi_{\alpha} :\RR^n\to\RR^n$
				such that
				\begin{multline*}
					\Big|
						\hat{\ep}^{\rCN,m}_{\tmkm,\tmk}
						-
							\sum_{1\le \alpha,\beta,\gamma\le d}
								\varphi_{\alpha,\beta,\gamma}(Y^{\rCN,m}_{\tmkm})
								B^{\alpha,\beta,\gamma}_{\tmkm,\tmk}
						-
							\sum_{1\le \alpha\le d}
								\psi_{\alpha}(Y^{\rCN,m}_{\tmkm})
								B^{\alpha}_{\tmkm,\tmk}\Delta_m
					\Big|\\
				\le C \Delta_m^{4H^-}
				\qquad
				\text{for all}
				\qquad
				1\le k\le 2^m.
				\end{multline*}
				Here, $C$ depends on $\sigma$ and $b$ polynomially.
	\end{enumerate}
\end{lemma}

\begin{proof}
	We show (1).
From \eqref{CN equation}, we have
\begin{multline}
	\label{eq4382908401}
	Y^{\rCN,m}_t-Y^{\rCN,m}_{\tmkm}
	=
		\sigma(Y^{\rCN,m}_{\tmkm})B_{\tmkm,t}
		+b(Y^{\rCN,m}_{\tmkm})(t-\tmkm)\\
		+
		\frac{1}{2}
		\left(
			\sigma(Y^{\rCN,m}_t)-\sigma(Y^{\rCN,m}_{\tmkm})
		\right)
		B_{\tmkm,t}
		+
		\frac{1}{2}
		\left(
			b(Y^{\rCN,m}_t)-b(Y^{\rCN,m}_{\tmkm})
		\right)
		(t-\tmkm).
\end{multline}
Hence
applying the Taylor formula and writing $Y^{\rCN,m}_{\tmkm,t}=Y^{\rCN,m}_t-Y^{\rCN,m}_{\tmkm}$,
we have
\begin{multline*}
	Y^{\rCN,m}_t
	-Y^{\rCN,m}_{\tmkm}
	-\sigma(Y^{\rCN,m}_{\tmkm})B_{\tmkm,t}
	-b(Y^{\rCN,m}_{\tmkm})(t-\tmkm)\\
	\begin{aligned}
	&=
		\frac{1}{2}
		\left(
			\int_0^1
			(D\sigma)(Y^{\rCN,m}_{\tmkm}+\theta Y^{\rCN,m}_{\tmkm,t})[Y^{\rCN,m}_{\tmkm,t}]d\theta
		\right)
		B_{\tmkm,t}    \\
		&\phantom{=}\qquad \qquad \qquad 
		+
		\frac{1}{2}
		\left(
			\int_0^1
			(Db)(Y^{\rCN,m}_{\tmkm}+\theta Y^{\rCN,m}_{\tmkm,t})[Y^{\rCN,m}_{\tmkm,t}]d\theta
		\right)
		(t-\tmkm)  \\
	&=
		((D\sigma)[\sigma])(Y^{\rCN,m}_{\tmkm})
		\left[
			\frac{1}{2}
			B_{\tmkm,t}
			\otimes
			B_{\tmkm,t}
		\right]
		+
		\hat{\ep}^{\rCN,m}_{\tmkm,t}\\
	&=
		((D\sigma)[\sigma])(Y^{\rCN,m}_{\tmkm})
		\BB_{\tmkm,t}
		+
		c(Y^{\rCN,m}_{\tmkm})d^{\rCN,m}_{\tmkm,t}
		+
		\hat{\ep}^{\rCN,m}_{\tmkm,t}.
	\end{aligned}
\end{multline*}

	We show (2).
	From \eqref{CN equation}, we have
	\begin{align*}
		\max_k
		\sup_{\tmkm\leq t\leq\tmk}
			|Y^{\rCN,m}_t-Y^{\rCN,m}_{\tmkm}|
		\le
			C
			|t-\tmkm|^{H^-}.
	\end{align*}
	This estimate and \eqref{eq4382908401} imply
	\begin{align*}
		\max_k
		\sup_{\tmkm\leq t\leq\tmk}
			|Y^{\rCN,m}_t-Y^{\rCN,m}_{\tmkm}-\sigma(Y^{\rCN,m}_{\tmkm})B_{\tmkm,t}|
		\le 
			C
			|t-\tmkm|^{2H^-}.
	\end{align*}
	Hence, by substituting
	\begin{multline*}
		(D\sigma)(Y^{\rCN,m}_{\tmkm}+\theta Y^{\rCN,m}_{\tmkm,t})[Y^{\rCN,m}_{\tmkm,t}]\\
		\begin{aligned}
			&=
				(D\sigma)(Y^{\rCN,m}_{\tmkm}+\theta Y^{\rCN,m}_{\tmkm,t})[\sigma(Y^{\rCN,m}_{\tmkm})B_{\tmkm,t}]
				+
				O(|t-\tmkm|^{2H^-})\\
			&=
				(D\sigma)(Y^{\rCN,m}_{\tmkm})[\sigma(Y^{\rCN,m}_{\tmkm})B_{\tmkm,t}]
				+
				O(|t-\tmkm|^{2H^-})
		\end{aligned}
	\end{multline*}
into \eqref{explicit form of hepmCN}, we can estimate 
the first term in \eqref{explicit form of hepmCN}.
Because the second term can be estimated in the same way,
we arrive at (2).

We show (3).
By a similar calculation to the above, we have
\begin{multline*}
\Big|
\hat{\ep}^{\rCN,m}_{\tmkm,\tmk}
-\frac{1}{2}(D\sigma)(Y^{\rCN,m}_{\tmkm})
\left[
\frac{1}{2}(D\sigma)[\sigma](Y^{\rCN,m}_{\tmkm})
[(B_{\tmkm,\tmk})^{\otimes 2}]
+b(Y^{\rCN,m}_{\tmkm})\Delta_m\right]
B_{\tmkm,\tmk}\\
-\frac{1}{4}(D^2\sigma)(Y^{\rCN,m}_{\tmkm})
\left[
\sigma(Y^{\rCN,m}_{\tmkm})B_{\tmkm,\tmk},\sigma(Y^{\rCN,m}_{\tmkm})B_{\tmkm,\tmk}
\right]B_{\tmkm,\tmk}\\
-\frac{1}{2}(Db)(Y^{\rCN,m}_{\tmkm})\left[\sigma(Y^{\rCN,m}_{\tmkm})B_{\tmkm,\tmk}\right]\Delta_m
\Big|
\le
C\Delta_m^{4H^-}.
\end{multline*}
Note that the above constants depend 
on $\sigma,b$ polynomially because $\omega\in \Omega_0^{(m)}$.
The proof completed.
\end{proof}

\begin{remark}
	Let $d^m=d^{\rIM,m}$, $d^{\rM,m}$, $d^{\rFE,m}$, $d^{\rCN,m}$
	and $\hepm=\hat{\ep}^{\rIM,m}$, $\hat{\ep}^{\rM,m}$, $\hat{\ep}^{\rFE,m}$, $\hat{\ep}^{\rCN,m}$.
For every $s,t\in \Dm$ with $s\le t$ , define
$d^m_t$, $d^m_{s,t}$, $\hepm_t$ and $\hepm_{s,t}$
in the same way as \eqref{eq89034802914832} with 
$\eta_{\tmim,\tmi}=d^m_{\tmim,\tmi},\hepm_{\tmim,\tmi}$.
\end{remark}

\subsection{Statement of main results}\label{2.3}
Now we are in a position to state our main results
(Theorem~\ref{main theorem}, Corollary~\ref{corollary for main theorem} and Theorem~\ref{main theorem fbm}).
In Section~\ref{2.2}, we recalled four approximation schemes and we
wrote the solutions as $\Ymh_t$.
They are continuous processes but the values at the discrete times
$\{\Ymh_t\}_{t\in \Dm}$ well approximate $\{\Ymh_t\}_{t\in [0,1]}$.
Also it is natural to consider approximate schemes defined at discrete times
$\Dm$ only for implementation.
As stated in Introduction, 
in Theorem~\ref{main theorem}, we consider the recurrence relations
of $\{\Ymh_t\}_{t\in \Dm}$
can be obtained by adding extra two terms containing $d^m$
and $\hepm$ to the recurrence relation of the Milstein scheme.
Since the Milstein scheme converges, we can expect that $\Ymh_t$ also converges to $Y_t$
if $d^m$ and $\hepm$ are small in a certain sense.
Based on this idea, we introduce smallness conditions
as Conditions~\ref{condition on dm}\,$\sim$\,\ref{condition on K}
and address approximate solutions and estimates of the errors at discrete times $\Dm$.
This is stated as Theorem~\ref{main theorem},
which is a result in a general setting not limited to the four schemes and fBms.
Corollary~\ref{corollary for main theorem} is a continuous version of Theorem~\ref{main theorem}.
In Theorem~\ref{main theorem fbm}, we give estimates of errors for the four schemes and fBms.
Note that we can check Conditions~\ref{condition on dm}\,$\sim$\,\ref{condition on K}
to use Corollary~\ref{corollary for main theorem}
for the four schemes except Crank-Nicolson scheme in the case of
fBm with the Hurst parameter $\frac{1}{3}<H\leq\frac{1}{2}$.
Although the Crank-Nicolson scheme can also be reduced to 
a setting satisfying the conditions,
it requires additional considerations.

Here we reset the notation to state Theorem~\ref{main theorem}.
For $\omega\in\Omega_0$, we define $\{\Ymh_t\}_{t\in \Dm}$ by
the following recurrence relation:
$\Ymh_0=\xi$ and
\begin{align}
	\Ymh_{\tmk}
  &=
  \Ymh_{\tmkm}+\sigma(\Ymh_{\tmkm})B_{\tmkm,\tmk}
    +((D\sigma)[\sigma])(\Ymh_{\tmkm})\BB_{\tmkm,\tmk}
    +b(\Ymh_{\tmkm})\Delta_m\nonumber\\
  &\qquad\qquad\qquad\qquad
    +c(\Ymh_{\tmkm})d^m_{\tmkm,\tmk}
    +\hat{\ep}^m_{\tmkm,\tmk}, \quad 1\le k\le 2^m.
  \label{equation of ymh}
\end{align}
Here $c\in C^3_b(\RR^n, L(\RR^d\otimes \RR^d,\RR^n))$ is a function
and 
$d^m=\{d^m_{\tmkm,\tmk}\}_{1\leq k\leq 2^m}\subset \RR^d\otimes \RR^d$ 
and $\hepm=\{\hepm_{\tmkm,\tmk}\}_{1\leq k\leq 2^m}\subset \RR^n$
are random variables defined on $\Omega_0$.
We now state our smallness conditions on $d^m$ and $\hepm$.
\begin{condition}\label{condition on dm}
There exist two pairs of positive numbers $(\vep_0,2H^-)$ and $(\vep_1,\lambda_1)$
with $\vep_i>0$ $(i=0,1)$ and $\lambda_1+H^->1$
and non-negative random variables $G_0=G_0(\vep_0,2H^-)$ and
$G_1=G_1(\vep_1,\lambda_1)$
which belong to $\cap_{p\ge 1}L^p(\Omega_0)$
such that 
    \begin{align*}
      |d^m_{s,t}|
      \le
		\min
			\left\{
				\Delta_m^{\vep_0}G_0|t-s|^{2H^-},
				\Delta_m^{\vep_1}
				G_1
				|t-s|^{\lambda_1}
			\right\}
        \quad
        \text{for all $s,t\in \Dm$ with $s<t$.}
   \end{align*}
\end{condition}
Although the reader might be interested in the reason 
why two exponents $2H^-$ and $\lambda_1$ are introduced, 
we defer the explanation to Remark~\ref{rem4320914321424} and 
proceed to state the conditions.
We next explain a condition on $\hepm$.
In this condition, although (1-a) follows from (2), 
we state (1-a) independently because it is used in Section~\ref{ymrk and Mmrk}.
Below, $B^{\alpha,\beta,\gamma}_{s,t}$ $(0\le s\le t\le 1)$ denotes
the $e_{\alpha}\otimes e_{\beta}\otimes e_{\gamma}$-component
of the third level rough paths which are constructed from $(B,\BB)$.
\begin{condition}\label{condition for epsilon}
	\begin{enumerate}
		\item	\begin{enumerate}
					\item[\upshape{(a)}]
						There exists a positive constant $C$ such that
						\begin{align}
							\label{eq4839012814214}
							|\hepm_{\tmkm,\tmk}|
							\le
								C \Delta_m^{3H^-}
							\qquad 
							\text{for all}
							\qquad
							1\le k\le 2^m,\quad \omega\in \Omega_0^{(m)}.
						\end{align}
						Here, $C$ depends on $\sigma$, $b$ and $c$ polynomially.
					\item[\upshape{(b)}]
						There exists a positive constant $C$ such that
						\begin{align}
						\label{eq4583902014}
							|\hepm_{\tmkm,\tmk}|
							\le
								C \Delta_m^{3H^-}
							\qquad 
							\text{for all}
							\qquad
							1\le k\le 2^m,\quad
							\omega\in \Omega_0\setminus \Omega_0^{(m)}
						\end{align}
						Here, $C$ depends on $\sigma$, $b$, $c$ and $C(B)$ polynomially.
		\end{enumerate}
		\item	There exist bounded Lipschitz continuous functions 
				$\varphi_{\alpha,\beta,\gamma} :\RR^n\to\RR^n$
				and $\psi_{\alpha} :\RR^n\to\RR^n$
				such that
				\begin{multline*}
					\Big|
						\hepm_{\tmkm,\tmk}
						-
							\sum_{1\le \alpha,\beta,\gamma\le d}
								\varphi_{\alpha,\beta,\gamma}(\hat{Y}_{\tmkm})
								B^{\alpha,\beta,\gamma}_{\tmkm,\tmk}
						-
							\sum_{1\le \alpha\le d}
								\psi_{\alpha}(\Ymh_{\tmkm})
								B^{\alpha}_{\tmkm,\tmk}\Delta_m
					\Big|\\
				\le C \Delta_m^{4H^-}
				\qquad
				\text{for all}
				\qquad
				1\le k\le 2^m,
				\quad \omega\in \Omega_0^{(m)}.
				\end{multline*}
				Here, $C$ depends on $\sigma, b$ and $c$ polynomially.
	\end{enumerate}
\end{condition}

Here we state the main non-trivial condition assumed in our main results.
For $c\in C^3_b(\RR^n, L(\RR^d\otimes \RR^d,\RR^n))$, which is used in \eqref{equation of ymh},
set 
\begin{align}\label{defIm}
  I^m_t
  =
    I^m_t(c,d^m)
  &=
    \sum_{i=1}^{\floor{2^mt}}
      J_{\tmim}^{-1}
      c(Y_{\tmim})
      d^m_{\tmim,\tmi}.
\end{align}
Let $I^m|_{\Dm}$ denote the discrete process defined as the restriction of $I^m$ on $\Dm$.
\begin{condition}\label{minimal moment estimate}
	Let $I^m|_{\Dm}$ be as above.
	For all $p\geq 1$, we have
	\begin{align*}
	  \sup_m
		\|
		  \|(2^m)^{2H-\frac{1}{2}}I^m|_{\Dm}\|_{H^-}
		\|_{L^p}<\infty.
	\end{align*}
\end{condition}

We explain the final condition.
Let $d^{m,\alpha,\beta}_{\tmim,\tmi}
=(d^{m}_{\tmim,\tmi},e_{\alpha}\otimes e_{\beta})$.
We set
\begin{align*}
  \tilde{\mathcal{K}}^3_m
  &=
    \Big\{
      \big\{ d^{m,\alpha,\beta}_{\tmim,\tmi}B^{\gamma}_{\tmim,\tmi} \big\}_{i=1}^{2^m},
\,
      \big\{ B^{\alpha,\beta,\gamma}_{\tmim,\tmi} \big\}_{i=1}^{2^m}, 
\,
\big\{B^{0,\alpha}_{\tmim,\tmi}\big\},
\,
\big\{B^{\alpha,0}_{\tmim,\tmi}\big\}
      ~\Big|~1\le \alpha,\beta,\gamma\le d
    \Big\}.
\end{align*}
and
\begin{align}\label{eqK}
	\mathcal{K}^3_m
	&=
		\left\{
			\{K^m_t\}_{t\in \Dm} \,\middle|\,
			K^m_t=\sum_{i=1}^{\floor{2^mt}}K^m_{\tmim,\tmi}
			\,\, \text{for some}\,\, 
	\big\{ K^m_{\tmim,\tmi} \big\}_{i=1}^{2^m}\in \tilde{\mathcal{K}}^3_m
		\right\}.
\end{align}
Here we set $K^m_0=0$ with convention. 
Note that $B^{0,\alpha}_{\tmim,\tmi}$, $B^{\alpha,0}_{\tmim,\tmi}$ 
are defined in (\ref{B0alpha}).

\begin{condition}\label{condition on K} 
	There exist a pair of positive numbers $(\vep_2,\la_2)$ with $\lambda_2+H^->1$
	and a non-negative random variable $G_2=G_2(\vep_2,\lambda_2)\in \cap_{p\ge 1}L^p(\Omega_0)$
	such that for all discrete processes $\{K^m_t\}_{t\in \Dm}\in \mathcal{K}^3_m$,
	\begin{align*}
		\big|(2^m)^{2H-\frac{1}{2}}
		K^{m}_{s,t}\big|
		&\le
			\Delta_m^{\vep_2}
			G_2
			|t-s|^{\lambda_2}
		\quad
		\text{for all $s,t\in \Dm$.}
	\end{align*}
\end{condition}

In the above condition, we consider $B^{\alpha,\beta,\gamma}_{\tmim,\tmi}$ only
in a subset of 
Wiener chaos of order 3 which can be obtained by iterated integrals of $B$.
However, noting the relation,
\begin{align}
	\label{product of iterated integral}
	\left\{
		\begin{aligned}
			B^{\alpha,\beta}_{s,t}B^\gamma_{s,t}
			&=
				B^{\alpha,\beta,\gamma}_{s,t}
				+B^{\gamma,\alpha,\beta}_{s,t}
				+B^{\alpha,\gamma,\beta}_{s,t},\\
			B^{\alpha}_{s,t}B^{\beta}_{s,t}B^{\gamma}_{s,t}
			&=
			\frac{1}{2}
			\left(
				B^{\alpha,\beta,\gamma}_{s,t}+
				B^{\beta,\alpha,\gamma}_{s,t}+
				B^{\beta,\gamma,\alpha}_{s,t}+
				B^{\alpha,\gamma,\beta}_{s,t}+
				B^{\gamma,\alpha,\beta}_{s,t}+
				B^{\gamma,\beta,\alpha}_{s,t}
			\right),
		\end{aligned}
	\right.
\end{align}
which follows from the geometric property of $(B,\BB)$,
we obtain similar estimates for sum processes defined by
the above increments.

We now state our first main result.
Note that we always assume Condition~\ref{condition on R}
on $B$.

\begin{theorem}\label{main theorem}
	Let $Y_t$ be the solution to RDE \eqref{rde}.
Let $c\in C^3_b(\RR^n, L(\RR^d\otimes \RR^d,\RR^n))$.
Let $d^m=\{d^m_{\tmkm,\tmk}\}_{k=1}^{2^m}
\subset \RR^d\otimes \RR^d$ and 
$
\hepm=\{\hepm_{\tmkm,t}\}_{k=1}^{2^m}
\subset \RR^n
$
be random variables defined on $\Omega_0$.
Consider the approximate solution $\Ymh_t$ $(t\in \Dm)$ 
defined by \eqref{equation of ymh}.
Let $I^m$ be the weighted sum process defined by \eqref{defIm}.
Set
\begin{align}
	\label{remainder term discrete}
	R^m_t
	&=
	\Ymh_t-Y_t-J_tI^m_t,\quad t\in \Dm.
\end{align}
Let $\frac{1}{2}(H+\frac{1}{4})<H^-<H$.
Assume that Conditions~\upshape{\ref{condition on dm}}\,$\sim$\,\upshape{\ref{condition on K}} hold.
	Then for 
 $0<\vep<\min \{3H^--1,4H^- - 2H - \frac{1}{2},\vep_1,\vep_2\}$,
	we have
	$2^{m(2H-\frac{1}{2}+\vep)}
        \max_{t\in \Dm}|R^m_t|\to 0$ in $L^p$ 
        for all $p\ge 1$
	and almost surely.
\end{theorem}

The next is a remark on how to use Condition~\ref{condition on dm}.
\begin{remark}\label{rem4320914321424}
	In our proof, we will use 
	the H\"older estimate of $d^m$ given by the pair $(\vep_0,2H^-)$
	to estimate an approximation of the Jacobian and its inverse
	(we write them as $\tJmr$ and $(\tJmr)^{-1}$ later) 
	by using Cass-Litterer-Lyons' estimate.
	On the other hand, the H\"older estimate given by the pair $(\vep_1,\la_1)$
	determines the convergence rate of
	the remainder term $R^m_t$ in our main theorems.
	More precisely, $\vep_1$ is one of upper bounds of the convergence rate
	and we obtain a good convergence rate if we can choose large $\vep_1$.

	A trivial choice of $(\vep_1,\la_1)$ is $(\vep_0,2H^-)$.
	In general, there is a trade-off between
	the H\"older exponent and the value of the H\"older norm.
	Hence for $\la_1<2H^-$ we may be able to take $\vep_1>\vep_0$.
	This is a good situation for our application.
	In fact we can implement this situation in our application.
	Therefore we may be able to take large $\vep_1$ for small $\la_1$.
	We refer the readers for this 
	to Remark~\ref{remark on vep and lambda}.
\end{remark} 

In the above theorem, $d^m_{s,t}$ and $\hepm_{s,t}$ are defined only at the discrete
times $(s,t)=(\tmkm,\tmk)$ $(1\le k\le 2^m)$.
However, they are defined at 
$\{\{(s, t)\}_{s=\tmkm, t\in [\tmkm,\tmk]}\}_{k=1}^{2^m}$
in some cases as in the four schemes we explained.
As a corollary of this theorem, we have the following result in such a situation.

\begin{corollary}\label{corollary for main theorem}
We consider the same situation as in Theorem~\upshape{\ref{main theorem}}.
Further we assume $d^m_{s,t}$ and $\hepm_{s,t}$ are defined at
$\{\{(s, t)\}_{s=\tmkm, t\in [\tmkm,\tmk]}\}_{k=1}^{2^m}$
and assume that there exists a positive random variable
$\hat{X}\in \cap_{p\ge 1}L^p(\Omega_0)$ such that
\begin{align}
	\label{eq48902834}
 |d^m_{\tmkm,t}|&\le \hat{X}|t-\tmkm|^{2H^-},
 &
|\hepm_{\tmkm,t}|&\le \hat{X}|t-\tmkm|^{3H^-}
\end{align}
for all $\tmkm<t\leq \tmk$ and $1\le k\le 2^m$.
We define $\Ymh_t$ $(0\leq t\leq 1)$ as an extension of $\Ymh_t$ $(t\in \Dm)$
via \eqref{equation of ymh}, with $\tmk$ replaced by $t(\in [\tmkm,\tmk])$.
Set
\begin{align}
	\label{remainder term conti}
	R^m_t=\Ymh_t-Y_t-J_tI^m_t,
	\quad 0\leq t\leq 1.
\end{align}
Then for the same constant $\vep$ as in Theorem~$\ref{main theorem}$,
we have
	$2^{m(2H-\frac{1}{2}+\vep)}
        \sup_{0\le t\le 1}|R^m_t|\to 0$ in $L^p$ 
        for all $p\ge 1$
	and almost surely.
\end{corollary}

We will prove the above results in Section~\ref{proof of main theorem}.
We make a remark on the estimate of $\vep$ in the above theorem.

\begin{remark}\label{remark on vep}
	We fix $H^-$ and lift
	$B$ to an $H^-$-H\"older rough path.
	It is necessary to give the meaning of
	the solutions $Y_t$ and $J_t$ of the differential equations.
	That is, they depends on the choice of $H^-$.
	However, note that each $\Ymh, Y_t, I^m_t$ are all almost surely defined
	for any choice of $\frac{1}{3}<H^-<H$ in our problem
	because any versions of $(B,\BB)$ are identical almost all
	$\omega$ for any $H^-$ as noted in Remark~\ref{remark on lift}.
	Therefore, the optimal constant of the estimate of $\vep$ 
	in Theorem~\ref{main theorem} should be independent of the choice of $H^-$.
\end{remark}

We now return to the four schemes stated in Section~\ref{2.2}.
We assume that $B$ is an fBm.
The following is the second main theorem.
\begin{theorem}\label{main theorem fbm}
	Let $B$ be an fBm with the Hurst parameter $\frac{1}{3}<H\leq \frac{1}{2}$. 
	Let $Y_t$ be the solution to RDE \eqref{rde}.
	Consider the implementable Milstein, Crank-Nicolson, Milstein or first-order Euler scheme
	and let $\Ymh_t$ and $I^m_t$ be their counterparts.
	Let $R^m_t$ $(0\leq t\leq 1)$ be defined by \eqref{remainder term conti}.
	Then for $0<\vep<3H-1$, we have
	$2^{m(2H-\frac{1}{2}+\vep)}\sup_t|R^m_t|\to 0$ 
 in $L^p$ for all $p\ge 1$
	and almost surely.
\end{theorem}

We will show Theorem~\ref{main theorem fbm} for 
the four schemes in Section~\ref{sec proof of main thm fbm part 1} 
with the help of Corollary~\ref{corollary for main theorem}.
For the implementable Milstein, Milstein, and first-order Euler schemes, 
we can check the conditions assumed in Corollary~\ref{corollary for main theorem}.
The Crank-Nicolson scheme 
satisfies Condition~\ref{condition for epsilon} only partially.
Namely, while Lemma~\ref{lem89430222} 
implies that Condition~\ref{condition for epsilon} (1-a) and (2) holds,
expression \eqref{definition of hepm for CN} 
yields that Condition~\ref{condition for epsilon} (1-b) does not hold.
Hence we cannot use Corollary~\ref{corollary for main theorem} directly.
However, it is easy to reduce the problem of Crank-Nicolson scheme to
the case which can be treated in Corollary~\ref{corollary for main theorem}.

We conclude this section with remarks on Theorem~\ref{main theorem fbm}.
\begin{remark}\label{rem849320482091}
	When we consider the Milstein scheme,
	we have $d^m_{\tmkm,\tmk}=d^{\rM,m}_{\tmkm,\tmk}=0$.
	From Theorem~\ref{main theorem fbm}, for any $\kappa>0$, we have 
	$
		(2^m)^{5H-\frac{3}{2}-\kappa}\sup_{t}|\Ymh_t-Y_t|\to 0
	$
	as $m\to\infty$
	in $L^p$ for all $p\ge 1$ and almost surely.
	For the other schemes, we have $(2^m)^{2H-1/2-\kappa}\sup_t|\Ymh_t-Y_t|\to 0$ 
	in the same sense.
	We will explain related weaker results in
	Theorem~\ref{rough estimate of tzmr} and
	Remark~\ref{estimate for J and J^-1}.
\end{remark}

\begin{remark}
	We mention related study with the above results.
	Ueda~\cite{ueda} studied the estimate of the remainder term in
	one-dimensional case. By ``one-dimensional'',
	we mean that the solution $Y_t$ and the driving fBm $B_t$ is one-dimensional.
	In this case, $H$ can be arbitrary positive number less than $1$.
	His study also is based on analysis of interpolation processes between
	the solutions and approximate solutions.
\end{remark}

\begin{remark}\label{rem2944299884}
	We make remarks on
	weak convergence of $(2^m)^{2H-\frac{1}{2}}I^m_t$ in the case of fBm.
	Let $B$ be an fBm.
	Let $d^m=d^{\rIM,m}=d^{\rCN,m}$.
	In this case, $d^{m,\alpha,\beta}_{\tmkm,\tmk}=(d^{m}_{\tmkm,\tmk},e_{\alpha}\otimes e_{\beta})$
	is given by
	\begin{align*}
	  d^{m,\alpha,\beta}_{\tmkm,\tmk}=\frac{1}{2}B^{\alpha}_{\tmkm,\tmk}
	B^{\beta}_{\tmkm,\tmk}-B^{\alpha,\beta}_{\tmkm,\tmk}.
	\end{align*}
	Note that
	$d^{m,\alpha,\beta}_{\tmkm,\tmk}=-d^{m,\beta,\alpha}_{\tmkm,\tmk}$
	holds because the rough path is geometric.
	Furthermore, we see that 
	$\{(2^m)^{2H-\frac{1}{2}}J_tI^m_t\}_{0\le t\le 1}$ weakly converges to 
	\begin{align}
	  \left\{C\sum_{1\le \alpha, \beta\le d} 
	  J_t\int_0^tJ_s^{-1}(D\sigma)(Y_s)[\sigma(Y_s)e_{\alpha}]e_{\beta}
	  dW^{\alpha,\beta}_s\right\}_{0\le t\le 1} \label{limit process}
	\end{align}
	in $D([0,1],\RR^n)$ with respect to the Skorokhod $J_1$-topology.
	Here
	\begin{enumerate}
	  \item $\{W^{\alpha,\beta}_t\}$ $(1\le \alpha<\beta\le d)$
			is a $\frac{1}{2}d(d-1)$-dimensional standard Brownian motion which is
			independent of the fBm $(B_t)$
			and $W^{\beta,\alpha}_t=-W^{\alpha,\beta}_t$ $(\beta>\alpha)$, 
			$W^{\alpha,\alpha}_t=0$ $(1\le \alpha\le d)$.
	  \item Let $\alpha\ne\beta$. The constant $C$ is given by
			\begin{align*}
			  C=
			  &
				\Bigg\{
					E[(B^{\alpha,\beta}_{0,1})^2]
					+2\sum_{k=1}^{\infty}E[B^{\alpha,\beta}_{0,1}B^{\alpha,\beta}_{k,k+1}]
				  -\frac{1}{4}(E[(B^{\alpha}_{0,1})^2])^2
				  -\frac{1}{2}
				  \sum_{k=1}^{\infty}E[B^{\alpha}_{0,1}B^{\alpha}_{k,k+1}]^2\Bigg\}^{\frac{1}{2}}.
			\end{align*}
  \end{enumerate}
  We proved this convergence in \cite{aida-naganuma2024LP}
  under the assumption $\sigma, b\in C^{\infty}_b$.
  Note that $I^m_t\equiv 0$ in the case where $d^m=d^{\rM,m}$.
  Also a similar convergence is proved in the case
where $d^m=d^{\rFE,m}$ by Liu-Tindel~\cite{liu-tindel} too.
See also \cite{aida-naganuma2024LP}.
\end{remark}

\begin{remark}[Weak convergence via Remark~\ref{rem2944299884} 
and Theorem~\ref{main theorem fbm}]
 Combining Remark~\ref{rem2944299884} 
and Theorem~\ref{main theorem fbm},
we can prove $\{(2^m)^{2H-\frac{1}{2}}(\Ymh_t-Y_t)\}$ weakly converges to
the weak limit of $\{(2^m)^{2H-\frac{1}{2}}J_tI^m_t\}$
in $D([0,1],\RR^n)$ in the Skorokhod topology.
This follows from the following more general result.
Let $\{Z^m_t\}_{0\le t\le 1}$, 
$\{\tilde{Z}^m_t\}_{0\le t\le 1}$ and 
$\{R^m_t\}_{0\le t\le 1}$ be $\RR^n$-valued c\`adl\`ag processes
such that $Z^m_t=\tilde{Z}^m_t+R^m_t$ holds almost surely.
Suppose that $\tilde{Z}^m$ converges weakly in
$D([0,1],\RR^n)$ and $\lim_{m\to\infty}E[\sup_t|R^m_t|]=0$.
Then $Z^m$ also converges weakly to the same limit of $\tilde{Z}^m$.
The reason is as follows.
$D([0,1],\RR^n)$ is a Polish space with respect to 
a metric $\rho$ on $D([0,1],\RR^n)$ which satisfies
$\rho(x,y)\le \sup_t|x_t-y_t|$.
To prove the convergence and the coincidence of the limit,
it suffices to show that $\lim_{m\to\infty}E[\varphi(Z^m)-\varphi(\tilde{Z}^m)]=0$
for any bounded Lipschitz continuous function $\varphi$ on $D([0,1],\RR^n)$.
Clearly, this can be proved by using
\[
	|\varphi(Z^m)-\varphi(\tilde{Z}^m)|\le \|\varphi\|_{{\rm Lip}}
\rho(Z^m,\tilde{Z}^m)\le
\|\varphi\|_{{\rm Lip}}\sup_t|R^m_t|,
\]
and the assumption on $R^m$.
\end{remark}

\subsection{Proof of Theorem~\ref{main theorem fbm}}
\label{sec proof of main thm fbm part 1}
In this section, we show Theorem~\ref{main theorem fbm}.
First, in the case of the four schemes, 
the implementable Milstein, Crank-Nicolson, Milstein and 
first-order Euler schemes, we show that 
Conditions~\upshape{\ref{condition on dm}}, \upshape{\ref{condition on K}}
and \upshape{\ref{minimal moment estimate}} hold, in this order,
and then give a proof of Theorem~\ref{main theorem fbm}.

\begin{lemma}\label{Holder estimate of d}
  Assume that $B$ is a $d$-dimensional fBm with $\frac{1}{3}<H\leq\frac{1}{2}$.
  Let $d^m$ be $d^{\rIM,m}$, $d^{\rCN,m}$, $d^{\rM,m}$ or $d^{\rFE,m}$.
  Then Condition~\upshape{\ref{condition on dm}} is satisfied
  for the pairs $(\vep_1,\lambda_1)$ and $(\vep_0,2H^-)$,
where
$0<\vep_1<3H^--1$, $\lambda_1=1+2H-3H^-$ and
$0<\vep_0<2(H-H^-)$.
\end{lemma}

\begin{proof}
  Since 
  \begin{align*}
    d^{\rFE,m}_{\tmim,\tmi}
    &=
      -
      \sum_{1\le \alpha\neq \beta\le d}
        B^{\alpha,\beta}_{\tmim,\tmi}e_\alpha\otimes e_\beta
      -
      \sum_{\alpha=1}^d
        \frac{1}{2}
        \left\{(B^{\alpha}_{\tmim,\tmi})^2-\Delta_m^{2H}\right\}
        e_{\alpha}\otimes e_{\alpha},
  \end{align*}
  all components of $d^m_{\tmim,\tmi}$,
  $d^{m,\alpha,\beta}_{\tmim,\tmi}=(d^m_{\tmim,\tmi},e_\alpha\otimes e_\beta)$,
  are written by a linear combination of
  \begin{align}
	\label{eq43902498109}
    B^{\alpha}_{\tmim,\tmi}B^{\beta}_{\tmim,\tmi},
   \quad B^{\alpha,\beta}_{\tmim,\tmi},\quad
   (B^{\alpha}_{\tmim,\tmi})^2-\Delta_m^{2H},
   \quad \alpha\ne \beta.
   \end{align}
   Hence we may assume $d^m_{\tmim,\tmi}$ 
   to be one of the above without loss of generality.
	These quantities are considered in several papers; 
	for example \cite{aida-naganuma2024LP}, \cite{liu-tindel}, \cite{naganuma2016}, 
	and \cite{nourdin-peccati}.
	In what follows, we assume $\frac{1}{3}<H<\frac{1}{2}$.
	For the case $H=\frac{1}{2}$, we can easily modify the discussion.
   
   For $k<l$, we have
   \begin{gather*}
    \left|
      E\left[B^\alpha_{\tmkm,\tmk}B^\beta_{\tmkm,\tmk}B^\alpha_{\tmlm,\tml}B^\beta_{\tmlm,\tml}\right]
    \right|
    \le
      C
      \left(
        \frac{|k-l|^{2H-2}}{2^{2mH}}
      \right)^2,\\
    \left|
      E\left[B^{\alpha,\beta}_{\tmkm,\tmk}B^{\alpha,\beta}_{\tmlm,\tml}\right]
    \right|
    \le
      C
      \left(
        \frac{|k-l|^{2H-2}}{2^{2mH}}
      \right)^2,\\
    \left|
      E
        \left[
          (B^{\alpha}_{\tmkm,\tmk})^2-\Delta_m^{2H})
          (B^{\alpha}_{\tmlm,\tml})^2-\Delta_m^{2H})
        \right]
    \right|
    \leq 
      C
      \left(
        \frac{|k-l|^{2H-2}}{2^{2mH}}
      \right)^2.
  \end{gather*}
  For $k=l$, the terms above can be estimate by $C(2^{-2mH})^2$.
We refer the readers for these estimates to Lemma~3.4
  in \cite{liu-tindel}.
Also we can find these estimates in Lemma~7.2 (1) in \cite{aida-naganuma2024LP}.
  These estimates imply
  \begin{align*}
    E[|d^m_{s,t}|^{2}]&\le C
   \left(\frac{1}{2^m}\right)^{4H-1}(t-s)\qquad
   \text{for $s,t\in \Dm$ with $s<t$.}
   \end{align*}
	Note that all constants $C$ above are independent of $m$ and $H$.
   By using the hypercontractivity of the Ornstein-Uhlenbeck
   semigroup, we get
   \begin{align}
    E[|d^m_{s,t}|^{p}]&\le C_p
   \left(\frac{1}{2^m}\right)^{(2H-\frac{1}{2})p}(t-s)^{\frac{p}{2}}
   \qquad
   \text{for $s,t\in \Dm$ with $s<t$.}
   \label{moment estimate for dm}
   \end{align}
   This estimate implies the next assertion.
   For $0<\kappa<\frac{1}{2}$, set
   \begin{align*}
    G_{m,\kappa}=(2^m)^{2H-\frac{1}{2}}\max_{s,t\in \Dm, s\ne t}\frac{|d^m_{s,t}|}
   {|t-s|^{\frac{1}{2}-\kappa}}.
   \end{align*}
   Then
   \begin{align}
	&\qquad\qquad
	\sup_m\|G_{m,\kappa}\|_{L^p}<\infty & &\text{for all $p\ge 1$},
	\label{integrability of G}\\
	&
	|d^m_{s,t}|\le
	\Delta_m^{2H-\frac{1}{2}}
	|t-s|^{\frac{1}{2}-\kappa}G_{m,\kappa} & & \text{for all $s,t\in \Dm$ with $s<t$}
	\label{estimate of dm 2}.
	\end{align}
   This can be checked as follows.
   Since we see \eqref{estimate of dm 2} from the definition of $G_{m,\kappa}$,
	we show integrability \eqref{integrability of G}.
   Let $\{\tilde{d}^m_t\}_{t\in [0,1]}$ be the piecewise linear extension of
   $\{d^m_t\}_{t\in \Dm}$.
   By (\ref{moment estimate for dm}), we have
   \begin{align*}
    E[|\tilde{d}^m_{s,t}|^p]&\le 3^{p-1}C_p
   \left(\frac{1}{2^m}\right)^{(2H-\frac{1}{2})p}|t-s|^{\frac{p}{2}}.
   \end{align*}
   By the Garsia-Rodemich-Rumsey inequality, we have for any $p,\theta>0$
   \begin{align*}
    \left(\sup_{s,t, s\ne t}\frac{|\tilde{d}^m_{s,t}|}{|t-s|^{\theta}}\right)^p
   \le 2\int_0^1\int_0^t\frac{|\tilde{d}^m_{s,t}|^p}{|t-s|^{2+p\theta}}
   dsdt.
   \end{align*}
   Combining these two inequalities and setting $\theta=\frac{1}{2}-\kappa$, we get
   \begin{align*}
    E[G_{m,\kappa}^p]&\le 
   2\cdot(2^m)^{(2H-\frac{1}{2})p}
\int_0^1\int_0^t\frac{E[|\tilde{d}^m_{s,t}|^p]}{|t-s|^{2+p\theta}}dsdt
   \le
   2\cdot 3^{p-1} C_p\int_0^1\int_0^t|t-s|^{\kappa p-2}dsdt.
   \end{align*}
   If $p>\kappa^{-1}$, then the right-hand side is bounded and
   we get
   \begin{align*}
   E[G_{m,\kappa}^p]\le
   2\cdot 3^{p-1}C_p
   \left(\kappa p(\kappa p-1)\right)^{-1},
   \end{align*}
   which proves \eqref{integrability of G}.
	
   By using \eqref{integrability of G} and \eqref{estimate of dm 2},
   we show the assertion.
Let us choose $0<\vep<2H-\frac{1}{2}$ and 
$0<2\kappa<\vep$.
Using $\Delta_m\le t-s$, we get
\begin{align*}
\text{(RHS of \eqref{estimate of dm 2})}
&=\Delta_m^{\vep-\kappa}
\Delta_m^{2H-\frac{1}{2}-\vep+\kappa}|t-s|^{\frac{1}{2}-\kappa}G_{m,\kappa}\\
&\le \Delta_m^{\vep-\kappa}|t-s|^{2H-\vep}G_{m,\kappa}\\
&=\Delta_m^{\vep-2\kappa}|t-s|^{2H-\vep}\Delta_m^{\kappa}G_{m,\kappa}
\end{align*}
Let $G_1=\sum_{m=1}^{\infty}\Delta_m^{\kappa}G_{m,\kappa}$.
This infinite series converges for $\mu$ almost all $\omega$.
Because for all $p\ge 1$,
\begin{align*}
	\|G_1\|_{L^p}
		\le
			\sum_{m=1}^{\infty}
				\Delta_m^{\kappa}
				\sup_m\|G_{m,\kappa}\|_{L^p}
		<
			\infty.
	\end{align*}
Combining the trivial estimate $\Delta_m^{\kappa}G_{m,\kappa}\le G_1$,
we get
   \begin{align*}
|d^m_{s,t}|
&\le \Delta_m^{\vep-2\kappa}|t-s|^{2H-\vep}G_1.
   \end{align*}
To check the validity of the statements for
the pairs $(\vep_1,\lambda_1)$ and $(\vep_0,2H^-)$,
it suffices to set $\vep=3H^--1(<2H-\frac{1}{2})$ 
and $\vep=2(H-H^-)(<2H-\frac{1}{2})$
respectively and choose $\kappa$ to be sufficiently small.
This completes the proof.
\end{proof}

\begin{remark}\label{remark on vep and lambda}
	We make a remark on the numbers appeared in Lemma~\ref{Holder estimate of d}.
	Recall that $\lambda_1=1+2H-3H^-$ and 
	that $3H^--1$ and $2(H-H^-)$ are the upper bounds of $\vep_1$ and $\vep_0$, respectively.
	We see that both inequalities $\lambda_1<2H^-$ and $3H^--1>2(H-H^-)$
	are equivalent to $5H^--2H>1$.
	The inequality $5H^--2H>1$ holds true if $H^-$ is sufficiently close to $H$
	because $H>\frac{1}{3}$.
	Hence we see that the good situation stated in Remark~\ref{rem4320914321424}
	is fulfilled.
\end{remark}

\begin{lemma}\label{K satisfies Condition}
  Assume that $B$ is a $d$-dimensional fBm with $\frac{1}{3}<H\leq\frac{1}{2}$.
  Let $d^m$ be $d^{\rIM,m}$, $d^{\rCN,m}$, $d^{\rM,m}$ or $d^{\rFE,m}$.
  Then Condition~\upshape{\ref{condition on K}} is satisfied
	for $\vep_2<3H^--1+(\frac{1}{2}-H)$ and $\la_2=1+2H-3H^-$.
\end{lemma}

\begin{proof}
	In what follows, we assume $\frac{1}{3}<H<\frac{1}{2}$.
	In the case where $H=\frac{1}{2}$, we can easily modify the discussion.
	Let $(K^m_t)\in \mathcal{K}^3_m$. 
First, we give estimates for variance of $K^m_{s,t}$.
We have for $s,t\in \Dm$ with $s<t$,
\begin{numcases}
 {E[|K^m_{s,t}|^2]\le}
 C \Delta_m^{6H-1}|t-s|  &\quad
\text{if $K^m_{\tmim,\tmi}=d^{m,\alpha,\beta}_{\tmim,\tmi}B^{\gamma}_{\tmim,\tmi}
\,\,\text{or}\,\, B^{\alpha,\beta,\gamma}_{\tmim,\tmi},$}\label{K estimate 1}\\
 C \Delta_m^{2H+1}|t-s|
& 
\quad
\text{if $K^m_{\tmim,\tmi}=B^{0,\alpha}_{\tmim,\tmi}$\,\,\text{or}\,\,
$K^m_{\tmim,\tmi}=B^{\alpha,0}_{\tmim,\tmi}$.}\label{K estimate 2}
\end{numcases}

Note that 
if the schemes are implementable Milstein or Crank-Nicolson scheme, then 
it is enough to consider the case $K^m=B^{\alpha,\beta,\gamma}$ only
for the proof of (\ref{K estimate 1}) because of the identities
(\ref{product of iterated integral}).
  Therefore, in those cases, from \cite[Lemma~4.3]{liu-tindel2020}, we see
(\ref{K estimate 1}) holds.
  In \cite{aida-naganuma2024LP},
  the same estimates are obtained in a little bit different way.
If the scheme is the first-order Euler scheme, then by the same reasoning as above,
it is sufficient to estimate $E[(\Delta_m^{2H} B^{\gamma}_{s,t})^2]$.
For this, we have
\begin{align*}
 E[(\Delta_m^{2H} B^{\gamma}_{s,t})^2]&\le C\Delta_m^{4H}|t-s|^{2H}\\
&=C\Delta_m^{4H}\cdot|t-s|^{2H-1}|t-s|\\
&\le C\Delta_m^{4H}\Delta_m^{2H-1}|t-s|=C\Delta_m^{6H-1}|t-s|.
\end{align*}
Actually we use Condition~\ref{condition on R} only to obtain this estimate.

Now we consider (\ref{K estimate 2}).
Let $K^m_{\tmim,\tmi}=B^{\alpha,0}_{\tmim,\tmi}=\int_{\tmim}^{\tmi}B^{\alpha}_{\tmim,u}du$.
	By using
	$
		|E[B^{\alpha}_{\tmim,u}B^{\alpha}_{\tmjm,v}]|
		\leq
			|E[B^{\alpha}_{\tmim,\tmi}B^{\alpha}_{\tmjm,\tmj}]|
	$
	for $\tmim\leq u\leq \tmi\leq \tmjm\leq v\leq \tmj$,
	we have
	\begin{align*}
		|
			E
				[
					K^m_{\tmim,\tmi}
					K^m_{\tmjm,\tmj}
				]
		|
		&\leq
			\int_{\tmim}^{\tmi}
				du
			\int_{\tmjm}^{\tmj}
				dv
				|E[B^{\alpha}_{\tmim,\tmi}B^{\alpha}_{\tmjm,\tmj}]|\\
		&\leq
			2^{-2m}
			2^{-2Hm}
			|E[B^{\alpha}_{0,1}B^{\alpha}_{j-i-1,j-i}]|.
	\end{align*}
	Noting $E[B^{\alpha}_{0,1}B^{\alpha}_{k-1,k}]\sim -H(1-2H)k^{2H-2}$ 
as $k\to\infty$,
	we have for $k2^{-m}=s<t=l2^{-m}$,
\begin{align*}
 E[(K^m_{s,t})^2]\le C\sum_{i,j=k+1}^l
2^{-2m}2^{-2Hm}|j-i|^{2H-2}\le
C(2^{-m})^{2H+1}|t-s|.
\end{align*}
As for $B^{\alpha,0}_{\tmim,\tmi}$, we have
$B^{0,\alpha}_{\tmim,\tmi}=B^{\alpha,0}_{\tmim,\tmi}-\Delta_mB^{\alpha}_{\tmim,\tmi}$.
Hence, we need to estimate $E[(\Delta_m B^{\alpha}_{s,t})^2]$.
Since $\Delta_m\le \Delta_m^{2H}$, this term is smaller than $E[(\Delta_m B^{\alpha}_{s,t})^2]$
and we get desired estimate.

Because $6H-1\le 2H+1$,
consequently, for all cases, we have
$E[|K^m_{s,t}|^2]\le C\Delta_m^{6H-1}|t-s|$.
  Combining the hypercontractivity of the Ornstein-Uhlenbeck semigroup
  and the estimates above, for all $p\ge 2$, we obtain
      \begin{align*}
        E[|K^{m}_{s,t}|^{p}]
        \le C_p\left(2^{-m}\right)^{(3H-\frac{1}{2})p}(t-s)^{\frac{p}{2}}
        \quad \text{for all $s,t\in \Dm$}.
      \end{align*}
      From the same argument as in (\ref{estimate of dm 2}), 
      for any $\frac{1}{2}>\kappa>0$ and $m$, 
      there exists a positive random variable $G'_{m,\kappa}$ satisfying
      $\sup_m\|G'_{m,\kappa}\|_{L^p}<\infty$ for all
      $p\ge 1$ such that
      \begin{align*}
        |K^{m}_{s,t}|&\le
		\Delta_m^{3H-\frac{1}{2}}|t-s|^{\frac{1}{2}-\kappa}G'_{m,\kappa}\quad
      \text{for all $s,t\in \Dm$},
      \end{align*}      
which implies
\begin{align}
        |(2^{m})^{2H-\frac{1}{2}}K^{m}_{s,t}|&\le
		\Delta_m^{\frac{1}{2}-H} \Delta_m^{2H-\frac{1}{2}}
|t-s|^{\frac{1}{2}-\kappa}G'_{m,\kappa}\quad
      \text{for all $s,t\in \Dm$}.\label{estimate of 2^-mK^m}
\end{align}
Note that $\Delta_m^{2H-\frac{1}{2}}|t-s|^{\frac{1}{2}-\kappa}$ appears
in the proof of Lemma~\ref{Holder estimate of d} (see \eqref{estimate of dm 2}).

Let us choose $0<\vep<2H-\frac{1}{2}$ and 
$0<2\kappa<\vep$.
Then again using $\Delta_m\le t-s$ and similarly to the estimate of $d^m_{s,t}$,
we get
\begin{align}
    |(2^{m})^{2H-\frac{1}{2}}K^{m}_{s,t}|
	&\le
		\Delta_m^{\frac{1}{2}-H}
		\Delta_m^{\vep-2\kappa}
		|t-s|^{2H-\vep}
		\Delta_m^{\kappa}
		G'_{m,\kappa}
\end{align}
and set $G_2=\sum_{m=1}^{\infty}\Delta_m^{\kappa}G'_{m,\kappa}$
which converges $\mu$-a.s. $\omega$ and
$\|G_2\|_{L^p}<\infty$ for all $p\ge 1$.
Again by using the trivial estimate $\Delta_m^{\kappa}G'_{m,\kappa}\le G_2$,
we get
\begin{align*}
	|(2^{m})^{2H-\frac{1}{2}}K^{m}_{s,t}|
	&\le
		\Delta_m^{\frac{1}{2}-H}
		\Delta_m^{\vep-2\kappa}
		|t-s|^{2H-\vep}
		G_2.
\end{align*}
Putting $\vep=3H^--1(<2H-\frac{1}{2})$, 
we completes the proof.
\end{proof}

\begin{lemma}\label{lem438928401}
	Assume that $B$ is a $d$-dimensional fBm with $\frac{1}{3}<H\leq\frac{1}{2}$. 
	Let $d^m$ be $d^{\rIM,m}$, $d^{\rCN,m}$, $d^{\rM,m}$ or $d^{\rFE,m}$.
	Then Condition~\upshape{\ref{minimal moment estimate}} holds.
\end{lemma}
\begin{proof}
	Recall $c=(D\sigma)[\sigma]\in C^3_b$.
	We show the case $\frac{1}{3}<H<\frac{1}{2}$.
	We use the result by Liu-Tindel \cite{liu-tindel}.
	They considered similar problems (Proposition~4.7 and Corollary~4.9 in \cite{liu-tindel}).
	We can use their result to show the assertion as follows.
	Note that $f_t=J^{-1}_tc(Y_t)\in \mathcal{L}(\RR^d\otimes\RR^d,\RR^n)$ and
	$g_t\in \mathcal{L}(\RR^d,\mathcal{L}(\RR^d\otimes\RR^d,\RR^n))$ 
	defined by $g_tv=(-J^{-1}_tD\sigma(Y_t)v)c(Y_t)+J^{-1}_tDc(Y_t)[\sigma(Y_t)v]$ for $v\in\RR^d$ 
	satisfy \cite[(4.12)]{liu-tindel} because $Y$ and $J^{-1}$
	are solutions to \eqref{rde} and \eqref{Jaconian inverse} respectively
	and they belong to $L^p$ for all $p\ge 1$.
	The integrability of $J_t^{-1}$ is due to \cite{cll}
	(see also Remark~\ref{estimate for J and J^-1}).
	Hence from Corollary~4.9 in \cite{liu-tindel}, we get
	$
		\|(2^m)^{2H-\frac{1}{2}}I^m_{s,t}\|_{L^p}
		\leq
			C(t-s)^{\frac{1}{2}}
	$
	for some constant $C$.
	This and the Garsia-Rodemich-Rumsey inequality imply the assertion.
	While the above proof is based on the result by Liu-Tindel \cite{liu-tindel},
	we can provide another proof of the assertion
	under the assumption that $\sigma, b\in C^{\infty}_b$ 
(see \cite{aida-naganuma2024LP}).
   
	Finally, we consider the case where $H=\frac{1}{2}$.
Actually, it is not difficult to check this case by
using the It\^o calculus.
	For the reader's convenience, we include the proof.
	Recall that $I^m_t$ in Condition~\upshape{\ref{minimal moment estimate}} is
defined by
$
 I^m_t=\sum_{i=1}^{2^mt}F_{\tmim}d^m_{\tmim,\tmi}
$
$(t\in \Dm)$,
where $F_{t}=J_t^{-1}c(Y_t)$.
We give an estimate of $E[|I^m_{s,t}|^{2p}]$ by applying martingale theory.
Since all components of $d^{m}_{\tmim,\tmi}$, 
$d^{m,\alpha,\beta}_{\tmim,\tmi}=(d^m_{\tmim,\tmi},e_\alpha\otimes e_\beta)$,
are written by a linear combination of \eqref{eq43902498109},
the desired estimates follow from those of
\begin{align}
	\label{eq439028014214}
	&
	\sum_{i=1}^{2^mt}
		F^{\alpha,\beta}_{\tmim}
		B^{\alpha}_{\tmim,\tmi}B^{\beta}_{\tmim,\tmi},
	&
	&
	\sum_{i=1}^{2^mt}
		F^{\alpha,\beta}_{\tmim}
		B^{\alpha,\beta}_{\tmim,\tmi},
	&
	&
	\sum_{i=1}^{2^mt}
		F^{\alpha,\alpha}_{\tmim}
		\{
			(B^{\alpha}_{\tmim,\tmi})^2-\Delta_m
		\},
\end{align}
where $F^{\alpha,\beta}_t=F_t(e_\alpha\otimes e_\beta)$ and $\alpha\neq\beta$.
For $t\in[0,1]$, let
\[
	\tilde{I}^m_t
	=
		\sum_{i=1}^{2^{m}}
			F^{\alpha,\beta}_{\tmim}
			B^{\alpha,\beta}_{\tmim\wedge t,\tmi\wedge t}.
\]
Clearly $I^m_t=\tilde{I}^m_t$ $(t\in \Dm)$ holds.
Note that
\[
 B^{\alpha}_{s,t}B^{\beta}_{s,t}=B^{\alpha,\beta}_{s,t}+B^{\beta,\alpha}_{s,t}
\quad (\alpha\ne\beta),
\quad (B^{\alpha}_{s,t})^2-(t-s)=\int_s^tB^{\alpha}_{s,u}dB^{\alpha}_u,
\]
where the integral in the second identity is the It\^o integral.
Therefore, for all cases in (\ref{eq439028014214}), 
it suffices to give the moment estimate of
\begin{align*}
 \tilde{I}^{m,\alpha,\beta}_t&=
\int_0^tF^{m,\alpha,\beta}_udB^{\beta}_u,\qquad 1\le \alpha,\beta\le d,
\end{align*}
where the integral is an It\^o integral and
$
F^{m,\alpha,\beta}_u=\sum_{i=1}^{2^m}F^{\alpha,\beta}_{\tmim}
B^{\alpha}_{\tmim,u}1_{[\tmim,\tmi)}(u).
$
Let $p>1$. We have
\begin{align}
	\nonumber
E\left[|\tilde{I}^{m,\alpha,\beta}_{s,t}|^{2p}\right]
&\le
CE\left[\left(\int_s^t|F^{m,\alpha,\beta}_u|^{2}du\right)^p\right]\\
&\le C(t-s)^{p-1}E\left[\int_s^t|F^{m,\alpha,\beta}_u|^{2p}du\right]
\le C'\left(\frac{t-s}{2^m}\right)^p,\label{moment estimate H=1/2}
\end{align}
where we have used the Burkholder-Davis-Gundy and the H\"older inequalities,
and
the estimate 
\[
 E[|F^{m,\alpha,\beta}_u|^{2p}]\le CE[|F^{\alpha,\beta}_{\tmim}|^{2p}]E[(B^{\alpha}_{\tmim,u})^{2p}]
\le C2^{-pm}\sup_tE[|F_t|^{2p}],
\quad \tmim\le u<\tmi.
\]
By the estimate (\ref{moment estimate H=1/2}) and a similar argument to 
the estimate (\ref{integrability of G}) of $d^m_{s,t}$,
we see that the assertion holds.

We conclude this proof with mentioning that,
under the assumption of this lemma,
Condition~\upshape{\ref{minimal moment estimate}} holds for all $H^-<\frac{1}{2}$
and that we can choose $H^-$ close to $\frac{1}{2}$.
\end{proof}

We now prove Theorem~\ref{main theorem fbm}.

\begin{proof}[Proof of Theorem~\ref{main theorem fbm}]
First, we prove the case of the implementable Milstein, Milstein
	and first-order Euler schemes.
Note that in these cases, $\hepm\equiv 0$ holds for the approximate solution
$\Ymh_t$.
Hence Condition~\upshape{\ref{condition for epsilon}} is clearly satisfied.
	From Lemmas~\ref{Holder estimate of d}, 
\ref{K satisfies Condition}, and \ref{lem438928401},
	we see that Conditions~\upshape{\ref{condition on dm}}, 
\upshape{\ref{condition on K}}
	and \upshape{\ref{minimal moment estimate}} hold.
	From the definition, \eqref{eq48902834} also holds.
	Hence the conditions assumed in Corollary~\ref{corollary for main theorem} 
are satisfied.
 By Corollary~\ref{corollary for main theorem}, 
for any $\vep<\min\{3H^--1,4H^--2H-\frac{1}{2}\}$, we have
$(2^m)^{2H-\frac{1}{2}+\vep}\sup_{0\leq t\leq 1}|R^m_t|\to 0$ in $L^p$ $(p\ge 1)$
and almost surely.
Since $H^-$ can be any positive number less than $H$
and $3H-1\leq 2H-\frac{1}{2}$, the proof is completed.

We consider the case of the Crank-Nicolson approximate solution
$Y^{\rCN,m}_t$.
We cannot directly apply Corollary~\ref{corollary for main theorem} to the 
Crank-Nicolson scheme 
	since it satisfies only Condition~\ref{condition for epsilon} (1-a) and (2).
	However we can reduce it to Corollary~\ref{corollary for main theorem}.
To this end,
we introduce an auxiliary approximate
solution $\Ymh_t$ defined via \eqref{equation of ymh}, with $\tmk$ replaced by $t(\in [\tmkm,\tmk])$
and
\begin{align*}
	d^m_{\tmkm,t}
	&=
		d^{\rCN,m}_{\tmkm,t},
	&
	\hepm_{\tmkm,t}
	&=
	\begin{cases}
		\hat{\epsilon}^{\rCN,m}_{\tmkm,t}, & \omega\in\Omega_0^{(m)},\\
		0, & \omega\in (\Omega_0^{(m)})^{\complement}.
	\end{cases}
\end{align*}
Lemmas~\ref{lem89430222} and the definition above imply
that $\hepm_{\tmkm,t}$ satisfies Condition~\upshape{\ref{condition for epsilon}}.
From Lemmas~\ref{Holder estimate of d}, \ref{K satisfies Condition}, and \ref{lem438928401},
we see that Conditions~\upshape{\ref{condition on dm}}, 
\upshape{\ref{condition on K}}
and \upshape{\ref{minimal moment estimate}} hold.
We see that $d^m_{\tmkm,t}$ and $\hepm_{\tmkm,t}$ satisfy \eqref{eq48902834}.
Hence we can apply
Corollary~\ref{corollary for main theorem} to $\Ymh_t$ defined above.
By using $\sup_{0\le t\le 1}|J_t|\in L^p$ $(p\ge 1)$
which is due to \cite{cll}, 
as a consequence of Corollary~\ref{corollary for main theorem}, we see that
$\sup_m\sup_{0\le t\le 1}|\Ymh_t|\in L^p$ $(p\ge 1)$.
Note that we will give selfcontained proof of the integrability of
$J_t$ and $J_t^{-1}$ in Remark~\ref{estimate for J and J^-1} (2)
and the integrability
of $\Ymh_t$ holds under weaker assumption as in Lemma~\ref{Ist general}
since (\ref{eq48902834}) holds.
Let $R^m_t=\Ymh_t-Y_t-J_tI^m_t$ and
$R^{\rCN,m}_t=Y^{\rCN,m}_t-Y_t-J_tI^m_t$.
Then using $\Ymh_t=Y^{\rCN,m}_t$ $(\omega\in \Omega_0^{(m)})$
and $Y^{\rCN,m}_t\equiv \xi$ $(\omega\in (\Omega_0^{(m)})^{\complement})$,
we have
\begin{align*}
	R^{\rCN,m}_t
	&=
		R^m_t+Y^{\rCN,m}_t-\Ymh_t\\
	&=
		R^m_t
		+
		\{Y^{\rCN,m}_t-\Ymh_t\}
		1_{\Omega_0^{(m)}}
		+
		\{Y^{\rCN,m}_t-\Ymh_t\}
		1_{(\Omega_0^{(m)})^\complement}\\
	&=
		R^m_t+\{\xi-\Ymh_t\} 1_{(\Omega_0^{(m)})^\complement}.
\end{align*}
By Corollary~\ref{corollary for main theorem}, we have
$(2^m)^{2H-\frac{1}{2}+\vep}\sup_{0\le t\le 1}|R^m_t|\to 0$ for all
$p\ge 1$ and almost surely.
By the integrability of $\sup_m\sup_{0\le t\le 1}|\Ymh_t|$ and
the estimate (\ref{complement of Omega0}),
we have
$(2^m)^{2H-\frac{1}{2}+\vep}\sup_{0\le t\le 1}|(\xi-\Ymh_t)
1_{(\Omega_0^{(m)})^\complement}|\to 0$ in $L^p$ 
and almost surely.
This completes the proof.
\end{proof}

\subsection{Small order nice discrete process}\label{small order process}
We introduce a class of discrete stochastic processes,
which includes $d^m_t$ satisfying Condition~\ref{condition on dm}.
Before doing so, we need to define a subset of $\Omega_0^{(m)}$.
For a positive number $\lambda_1$ satisfying $\lambda_1+H^->1$,
we introduce the following set:
\begin{align*}
\Omegam
&=
	\{\,
		\omega\in \Omega_0^{(m)}
		~|~
		\|d^m(\omega)\|_{2H^-}\le 1,\quad \|d^m(\omega)\|_{\lambda_1}\le 1
	\}.
\end{align*}
Similarly to the estimate of the complement of $\Omega^{(m)}_0$,
if Condition~\ref{condition on dm} holds with the same exponent 
$\lambda_1$ in the definition of $\Omegam$,
we can prove that for any $p\ge 1$,
there exists $C_p>0$ such that
\begin{align}\label{complement of Omegam}
\mu\left((\Omegam)^\complement\right)\le C_{p}2^{-mp}
\end{align}
which implies the complement of $\Omegam$ is also negligible set for our 
problem.

\begin{definition}
	\phantomsection
	\label{nice discrete process}
 \begin{enumerate}
  \item Let $\eta=\{(\eta^m_t)_{t\in \Dm}; m\ge m_0\}$ 
be a sequence of Banach space
valued random variables such that
$\eta^m_0=0$ and $\{\eta^m_t\}_{t\in \Dm}$ is
defined on $\Omegam$ for each $m$,
where $m\ge m_0$ and $m_0$ is a non-random constant and 
depends on the sequence.
Let $\{a_m\}$ be a positive sequence which converges to $0$.
Let $\la$ be a positive number such that $\la+H^->1$.
We say that $\eta=(\eta^m)$ is a $\{a_m\}$-order nice discrete process 
with the H\"older exponent $\la$ if
there exists a positive random variable $X\in \cap_{p\ge 1}L^p(\Omega_0)$ 
which is independent of $m$ such that
\begin{align}
 \|\eta^m_t-\eta^m_s\|\le a_m X(\omega)|t-s|^{\la}
\qquad \text{for all $m\ge m_0$,\,
$t,s\in \Dm$,\,
$\omega\in \Omegam$}.
\label{small order discrete process inequality}
\end{align}

\item Let $\{v^m_{\theta}\}_{\theta\in \Theta}$
be a family of Banach space valued
random variables defined on $\Omegam$,
where $m\ge m_0$.
Let $\{a_m\}$ be a positive sequence which converges to $0$.
If there exists a non-negative random variable 
$X\in \cap_{p\ge 1}L^p(\Omega_0)$ 
which does not depend on $m$ such that
\begin{align*}
 \sup_{\theta\in \Theta}\|v^m_{\theta}\|\le a_m X(\omega) \quad 
\text{for all $m$ and $\omega\in \Omegam$},
\end{align*}
then we write
\begin{align*}
 \sup_{\theta\in \Theta}\|v^m_{\theta}\|=O(a_m).
\end{align*}
 \end{enumerate}
\end{definition}

\begin{remark}\label{rem89034111343}
	Here we give examples of small order nice discrete processes.
	\begin{enumerate}
		\item	Let $\epm_{\tmkm,\tmk}$ be given by 
\eqref{definition of epm}. Assume that 
Conditions~\ref{condition on dm}, \ref{condition for epsilon} (1) and 
\ref{condition on K} are satisfied.
Let $\vep_1,\la_1, \vep_2, \la_2$ be the numbers appeared 
in Condition~\ref{condition on dm} and \ref{condition on K}.
Set $a_m=\max\{\Delta_m^{3H^{-}-1},\Delta_m^{\vep_1},\Delta_m^{\vep_2}\}$
and $\lambda=\min\{2H^-,\lambda_1,\lambda_2\}$.
				Let $\omega\in\Omega_0$.
				Then there exists a non-negative random variable
		$X\in \cap_{p\ge 1}L^p(\Omega_0)$ which is independent of $m$ 
				such that
				\begin{align}
				|d^m_{s,t}|
				+|\epm_{s,t}|
				+|\hepm_{s,t}|
				+|(2^m)^{2H-\frac{1}{2}}K^m_{s,t}|
				\le
					a_m
					X
					|t-s|^\lambda
				\quad
				\text{for all}
				\quad
					s, t\in \Dm.
				\label{condition on dm ep eph}
				\end{align}
				In particular, $d^m_t$, $\ep^m_t$, $\hat{\ep}^m_t$
				and $(2^m)^{2H-\frac{1}{2}}K^m_t$
				are $\{a_m\}$-order 
nice discrete processes with the H{\"o}lder exponent $\lambda$.
We need to check $\epm$ and $\hepm$ satisfy the inequality.
		For $s=\tml$ and $t=\tmk$, Lemma~\ref{lem8345901809321}
and Condition~\ref{condition for epsilon} (1)
				imply
				\begin{align*}
				|\epm_{s,t}|
				+
				|\hepm_{s,t}|
				=
					\sum_{i=l+1}^k
						\left\{
							|\epm_{\tmim,\tmi}|
							+
							|\hepm_{\tmim,\tmi}|
						\right\}
				\le C(k-l)\Delta_m^{3H^-}
				\le C\Delta_m^{3H^--1}|t-s|^{2H^-},
				\end{align*}
where the constant $C$ depends
$\sigma$, $b$, $c$ and $C(B)$ polynomially.
If we consider the pair $(\vep_0, 2H^-)$, we can prove
that there exist $\tilde{X}\in \cap_{p\ge 1}L^p(\Omega_0)$ and 
$\tilde{a}_m=\max\{\Delta_m^{\vep_0},\Delta_m^{3H^--1}\}$ such that
\begin{align*}
 |d^m_{s,t}|+|\epm_{s,t}|+|\hepm_{s,t}|\le \tilde{a}_m
\tilde{X}|t-s|^{2H^-}.
\end{align*}
We use the estimate (\ref{condition on dm ep eph}) in Sections~\ref{4.2} and \ref{4.4}.
	\item	In the above definition of $\{a_m\}$-order nice discrete processes,
				we assume the strong assumption on $X$ such that
				$X\in \cap_{p\ge 1}L^p(\Omega_0)$.
			Under Conditions~\ref{condition on R} and~\ref{condition on dm},
			we have many examples which satisfy this strong conditions.
   \end{enumerate}
\end{remark}

\begin{remark}\label{discrete young integral}
	Suppose a Banach space valued discrete process
	$F=\{(F^m_t)_{t\in \Dm};$ $m\ge m_0\}$ defined on $\Omegam$ satisfy
	the H\"older continuity
	\begin{align*}
	& \|F^m_t-F^m_s\|\le X_F(\omega)|t-s|^{H^{-}}\quad \,\,
	\text{for all $m\ge m_0$,\, $s,t\in \Dm$,\, $\omega\in \Omegam$}, \\
	& \sup_m\|F^m_0(\omega)\|\le Y_F(\omega)\quad\quad \text{for}~~
	\omega\in \Omegam.
	\end{align*}
	Here $X_F, Y_F\in \cap_{p\ge 1}L^p(\Omega_0)$ are random variables
	independent of $m$.
	If  $\eta=(\eta^m)$ is a real valued $\{a_m\}$-order nice
	discrete process with the H\"older exponent $\la$, 
	then
	\begin{align*}
		\tilde{\eta}^m_{\tmk}=\sum_{i=1}^kF^m_{\tmim}\eta^m_{\tmim,\tmi}
	\end{align*}
	is also a $\{a_m\}$-order nice discrete process with the H\"older exponent
        $\la$ by the estimate of 
	the (discrete) Young integral (see \cite{friz-victoir}):
	\begin{align*}
	\|\tilde{\eta}^m\|_{\la}&\le C\left(\|F_0^m\|+\|F^m\|_{H^-}\right)
	\|\eta^m\|_{\la},
	\end{align*}
	where $C$ is a constant depending only on $H^-$ and $\la$.
    Note that we used $\la+H^->1$.

	This property is very nice for our purpose.
	However, in our application, since the estimate on $F^m$ is satisfied only on $\Omegam$,
	we cannot require \eqref{small order discrete process inequality}
	for all $\omega\in\Omega_0$
	to be nice discrete processes.
\end{remark}

\begin{remark}
  In what follows,
  we use the following elementary summation by parts formula several times:
  For sequences $\{f_i\}_{i=0}^n$, $\{g_i\}_{i=0}^n$, we have
  \begin{align}
    \sum_{i=1}^n
      f_{i-1}g_{i-1,i}
    &=
      f_ng_n-f_0g_0
      -
      \sum_{i=1}^nf_{i-1,i}g_i.
  \label{summation by parts formula}
  \end{align}
  We will use this formula when we give estimates of discrete Young integral.
\end{remark}

 \section{An interpolation of discrete rough differential equations}
\label{interpolation}
Let $Y_t$ and $\Ymh_t$ be a solution to \eqref{rde} and
an approximate solution given by \eqref{equation of ymh}, respectively.
In previous section, we observe that the discrete stochastic processes 
$\{Y_t\}_{t\in \Dm}$ and $\{\Ymh_t\}_{t\in \Dm}$
corresponding to the solution and our approximate solutions
respectively 
of the RDE 
satisfy the following common recurrence form:
$Y_0=\Ymh_0=\xi$ and, for $1\leq k\leq 2^m$,
\begin{align*}
	Y_{\tmk}
	&=
		Y_{\tmkm}
		+\sigma(Y_{\tmkm})B_{\tmkm,\tmk}
		+((D\sigma)[\sigma])(Y_{\tmkm})\BB_{\tmkm,\tmk}
		+b(Y_{\tmkm})\Delta_m
		+\epm_{\tmkm,\tmk},\\
	\Ymh_{\tmk}
	&=
		\Ymh_{\tmkm}
		+\sigma(\Ymh_{\tmkm})B_{\tmkm,\tmk}
		+((D\sigma)[\sigma])(\Ymh_{\tmkm})\BB_{\tmkm,\tmk}
		+b(\Ymh_{\tmkm})\Delta_m\\
	&\qquad\qquad\qquad
		+c(\Ymh_{\tmkm})d^m_{\tmkm,\tmk}
		+\hepm_{\tmkm,\tmk}.
\end{align*}
We now introduce an interpolation process between 
$\{Y_t\}_{t\in \Dm}$ and $\{\Ymh_t\}_{t\in \Dm}$
to study the difference $\Ymh_t-Y_t$.
Moreover, we introduce a matrix valued process $\tJmr_t$ which approximates
the derivative process $J_t$ when $m\to \infty$.
Note that, in this section, 
we do not use any specific forms of $d^m$ and $\hepm$
which were given in Section~\ref{statement}.
Taking a look at the recurrence equations, we see that
the different points between $\Ymh_t$ and $Y_t$ are the terms
$c(\Ymh_{\tmkm})d^m_{\tmkm,\tmk}$,
$\hat{\ep}^m_{\tmkm,\tmk}$ and $\ep^m_{\tmkm,\tmk}$.
In view of this, we define a sequence
$\{\Ymr_t\}_{t\in \Dm}$ by the following recurrence relation:
$\Ymr_0=\xi$ and, for $1\leq k\leq 2^m$,
\begin{align}
	\Ymr_{\tmk}
	&=
		\Ymr_{\tmkm}
		+\sigma(\Ymr_{\tmkm})B_{\tmkm,\tmk}
		+((D\sigma)[\sigma])(\Ymr_{\tmkm})\BB_{\tmkm,\tmk}
		+b(\Ymr_{\tmkm})\Delta_m \nonumber \\
	&\qquad\qquad\qquad
		+\rho c(\Ymr_{\tmkm})d^m_{\tmkm,\tmk}
		+\rho\hepm_{\tmkm,\tmk}
		+(1-\rho)\epm_{\tmkm,\tmk}.
	\label{def_interplation}
\end{align}
Note that $Y^{m,0}_t=Y_t$ and $Y^{m,1}_t=\Ymh_t$ ($t\in \Dm$).
In this paper, we call this recurrence relation a discrete RDE.
The function $[0,1]\ni\rho\mapsto \Ymr_t$ is smooth and
\begin{align*}
	\Ymh_t-Y_t=\int_0^1\partial_{\rho}\Ymr_t d\rho
\end{align*}
holds.
We give the estimate for
$\Ymh_t-Y_t$ by using the estimate of $\Zmr_t=\partial_{\rho}\Ymr_t$.
Then $\{\Zmr_t\}_{t\in \Dm}$ satisfies $\Zmr_0=0$ and, for $1\leq k\leq 2^m$,
\begin{align}
	\Zmr_{\tmk}
	&=
		\Zmr_{\tmkm}
		+(D\sigma)(\Ymr_{\tmkm})[\Zmr_{\tmkm}]B_{\tmkm,\tmk}
		+\left(D((D\sigma)[\sigma])\right)(\Ymr_{\tmkm})[\Zmr_{\tmkm}]\BB_{\tmkm,\tmk} \nonumber \\
	&\phantom{=}\qquad
		+(Db)(\Ymr_{\tmkm})[\Zmr_{\tmkm}]\Delta_m 
		+\rho (Dc)(\Ymr_{\tmkm})[\Zmr_{\tmkm}]d^m_{\tmkm,\tmk}\nonumber \\
	&\phantom{=}\qquad
		+c(\Ymr_{\tmkm})d^m_{\tmkm,\tmk}
		+\hepm_{\tmkm,\tmk}
		-\epm_{\tmkm,\tmk},
  \label{zmr identity}
\end{align}
where
\begin{align}
  \label{DDsigma2}
  (D((D\sigma)[\sigma]))(y)[\eta] v\otimes w
  =
    D^2\sigma (y)[\eta,\sigma(y)v]w+D\sigma(y)[D\sigma(y)[\eta]v]w
\end{align}
for $y,\eta\in\RR^n$ and $v,w\in\RR^d$ (see also \eqref{eqSimplifiedNotationDef}).

We introduce the $\mathcal{L}(\RR^n)$-valued, that is,
matrix valued process
$\{\tJmr_t\}_{t\in \Dm}$ 
to obtain the estimates of $\{\Zmr_t\}_{t\in \Dm}$.
Let $\{\tJmr_t\}_{t\in \Dm}$ be the solution to the following 
recurrence relation:
$\tJmr_0=I$ and, for $1\leq k\leq 2^m$,
\begin{align}
	\tJmr_{\tmk}
	&=
		\tJmr_{\tmkm}
		+[D\sigma](\Ymr_{\tmkm})[\tJmr_{\tmkm}]B_{\tmkm,\tmk}
		+\left(D((D\sigma)[\sigma])\right)(\Ymr_{\tmkm})[\tJmr_{\tmkm}]\BB_{\tmkm,\tmk}\nonumber\\
	&
		\qquad
		\qquad
		+(Db)(\Ymr_{\tmkm})[\tJmr_{\tmkm}]\Delta_m+\rho (Dc)(\Ymr_{\tmkm})[\tJmr_{\tmkm}]d^m_{\tmkm,\tmk}.
	\label{Mmr}
\end{align}
Clearly, we can represent $\{\Zmr_t\}_{t\in \Dm}$ by using $\{\tJmr_t\}_{t\in \Dm}$ and 
$\{(\tJmr_t)^{-1}\}_{t\in \Dm}$ if $\tJmr_t$ are invertible by a constant variation method.
Actually, such kind of representation holds in general case too.
To show this, and for later purpose, we consider 
discrete RDEs which are driven by time shift process
of $B_t$.

Let $u\in\Dm$ with $u\leq 1-\Delta_m$.
For $\tmk\leq 1-u$, we introduce time shift variables:
\begin{align*}
(\theta_uB)_{\tmkm,\tmk}
&=
	B_{u+\tmkm,u+\tmk},
&
(\theta_u\BB)_{\tmkm,\tmk}
&=\BB_{u+\tmkm,u+\tmk},\\
(\theta_ud^m)_{\tmkm,\tmk}
&=
	d^m_{u+\tmkm,u+\tmk},\\
 (\theta_u\epm)_{\tmkm,\tmk}
 &=\epm_{u+\tmkm,u+\tmk},
 &
 (\theta_u\hepm)_{\tmkm,\tmk}
 &=\hepm_{u+\tmkm,u+\tmk}.
\end{align*}
For general $x\in \RR^n$, we define 
a discrete process $\{\Ymr_t(x)\}_{t\in\Dm,0\leq t\leq 1-u}$
by $\Ymr_0(x)=x$ and, for $\tmk\leq 1-u$,
\begin{align*}
	\Ymr_{\tmk}(x)
	&=
		\Ymr_{\tmkm}(x)
		+\sigma(\Ymr_{\tmkm}(x))(\theta_uB)_{\tmkm,\tmk}\\
	&\qquad\qquad
		+((D\sigma)[\sigma])(\Ymr_{\tmkm}(x)) (\theta_u\BB)_{\tmkm,\tmk}
		+b(\Ymr_{\tmkm}(x))\Delta_m\\
	&\qquad\qquad
		+\rho c(\Ymr_{\tmkm}(x))(\theta_ud^m)_{\tmkm,\tmk}
		+\rho(\theta_u\hepm)_{\tmkm,\tmk}
		+(1-\rho)(\theta_u\epm)_{\tmkm,\tmk}.
\end{align*}
To make clear the dependence of the driving process,
we may denote the solution of the above equation
by $\Ymr_t(x,\theta_uB)$.
For simplicity, we write $\Ymr_t$ for $\Ymr_t(\xi,B)$.
Using these notation,
we have $\Ymr_t(\Ymr_u(\xi,B),\theta_uB)=\Ymr_{u+t}(\xi,B)$.
We consider the case where $x=\Ymr_u$ ($u\in\Dm$ with $u\leq 1-\Delta_m$) below.

We now explain explicit representation of $\tJmr_t$.
For given $x\in \RR^n$, let
\begin{align}\label{defEmr}
  \Emr(x,\theta_tB)
  &=
    I
    +(D\sigma)(x)B_{t,t+\Delta_m}
    +D((D\sigma)[\sigma])(x)\BB_{t,t+\Delta_m} \nonumber \\
  &\phantom{=}\qquad
    +(Db)(x)\Delta_m
    +\rho (Dc)(x)d^m_{t,t+\Delta_m}.
\end{align}
Then for $t\in D_m$ with $t>0$, we have
 \begin{align*}
  \tJmr_t&=
  \Emr(\Ymr_{t-\Delta_m},\theta_{t-\Delta_m}B)
  \Emr(\Ymr_{t-2\Delta_m},\theta_{t-2\Delta_m}B)
  \cdots
  \Emr(\xi,B).
\end{align*}
Since $\tJmr_t$ depends on $\xi$ and $B$, 
we may denote $\tJmr_t$ by $\tJmr_t(\xi,B)$.
Next we define $\tJmr_t(\Ymr_u,\theta_uB)$
similarly to $\Ymr_t(x,\theta_uB)$.
That is, $\tJmr_t(\Ymr_u,\theta_uB)$ is defined by substituting
$\Ymr_u(=\Ymr_u(\xi,B))$, 
$\theta_uB$, $\theta_u\BB$, $\theta_ud^m$ for
$\xi$, $B$, $\BB$, $d^m$ in the equation \eqref{Mmr} of $\tJmr_t(=\tJmr_t(\xi,B))$.
Using $\Ymr_t(\Ymr_u,\theta_uB)=\Ymr_{u+t}(\xi,B)$,
we see that
$\tJmr_t(\Ymr_u,\theta_uB)$
satisfies $\tJmr_0(\Ymr_u, \theta_uB)=I$ and, for $\tmk\leq 1-u$,
\begin{align*}
	\tJmr_{\tmk}(\Ymr_u,\theta_uB)
	&=
		\tJmr_{\tmkm}(\Ymr_u,\theta_uB)
		+[D\sigma](\Ymr_{u+\tmkm})[\tJmr_{\tmkm}(\Ymr_u,\theta_uB)]B_{u+\tmkm,u+\tmk} \\
	&\qquad\qquad
		+\left(D((D\sigma)[\sigma])\right)(\Ymr_{u+\tmkm})[\tJmr_{\tmkm}(\Ymr_u,\theta_uB)]\BB_{u+\tmkm,u+\tmk}\\
	&\qquad\qquad
		+(Db)(\Ymr_{u+\tmkm})[\tJmr_{\tmkm}(\Ymr_u,\theta_uB)]\Delta_m\\
	&\qquad\qquad
		+\rho (Dc)(\Ymr_{u+\tmkm})[\tJmr_{\tmkm}(\Ymr_u,\theta_uB)]d^m_{u+\tmkm,u+\tmk}.
\end{align*}
From this equation, we obtain
\begin{align}
	\tJmr_{\tmk}(\Ymr_u,\theta_uB)
	&=
		\Emr(\Ymr_{u+\tmkm},\theta_{u+\tmkm}B)
		\tJmr_{\tmkm}(\Ymr_u,\theta_uB),\label{EMmrk}
\end{align}
which implies
\begin{multline}
	\tJmr_t(\Ymr_u,\theta_uB)\\
	=
		\Emr(\Ymr_{u+t-\Delta_m},\theta_{u+t-\Delta_m}B)
		\Emr(\Ymr_{u+t-2\Delta_m},\theta_{u+t-2\Delta_m}B)
		\cdots
		\Emr(\Ymr_u,\theta_uB)\label{time shift of Mmrk}
\end{multline}
Also we have, for $s,t\in\Dm$ with $s+t\leq 1-u$,
\begin{align}
	\tJmr_{s+t}(\Ymr_u,\theta_uB)
	&=
		\tJmr_s(\Ymr_{u+t},\theta_{u+t}B)
		\tJmr_t(\Ymr_u,\theta_uB).\label{Mmrk'+k}
\end{align}
The proof of \eqref{Mmrk'+k} is as follows.
By \eqref{time shift of Mmrk}, we have
\begin{align*}
	\tJmr_{s+t}(\Ymr_u,\theta_uB)
	&=
		\Emr(\Ymr_{u+t+s-\Delta_m},\theta_{u+t+s-\Delta_m}B)
		\cdots
		\Emr(\Ymr_{u+t},\theta_{u+t}B)\\
		&\qquad
		\cdot
		\Emr(\Ymr_{u+t-\Delta_m},\theta_{u+t-\Delta_m}B)
		\cdots
		\Emr(\Ymr_{u},\theta_{u}B)\\
	&=
		\tJmr_s(\Ymr_{u+t},\theta_{u+t}B)
		\tJmr_t(\Ymr_{u},\theta_{u}B).
\end{align*}

We have the following lemma for the invertibility of $\tJmr_t$.

\begin{lemma}\label{expansion of Mmr}
	For $1\leq k\leq 2^m$, we have
\begin{align}
	\tJmr_{\tmk}
	&=
		\Emr(\Ymr_{\tmkm},\theta_{\tmkm}B)\tJmr_{\tmkm}\nonumber\\
	&=
		\Bigl(
			I
			+(D\sigma)(\Ymr_{\tmkm})B_{\tmkm,\tmk}
			+D((D\sigma)[\sigma])(\Ymr_{\tmkm})\BB_{\tmkm,\tmk}\nonumber\\
	&\qquad
	\qquad
	\qquad
	\qquad
	\qquad
			+\rho (Dc)(\Ymr_{\tmkm})d^m_{\tmkm,\tmk}
			+(Db)(\Ymr_{\tmkm})\Delta_m
		\Bigr)
		\tJmr_{\tmkm},
\label{expansion of Mmr2}
\end{align}
    and for large $m$, $\tJmr_t$ are invertible.
For example, for any $\omega\in \Omegam$, if $m$ satisfies
\begin{align}
 \Delta_m^{H^-}\|D\sigma\|+\Delta_m^{2H^-}
\|D\left((D\sigma)[\sigma]\right)\|+
\Delta_m^{2H^-}\|Dc\|+\Delta_m\|Db\|\le \frac{1}{2},
\label{assumption on m}
\end{align}
then $\Emr(\Ymr_{\tmkm},\theta_{\tmkm}B)$ is invertible and
it holds that
\begin{align}
 \left|\Emr(\Ymr_{\tmkm},\theta_{\tmkm}B)^{-1}-I+(D\sigma)(\Ymr_{\tmkm})B_{\tmkm,\tmk}\right|&\le C\Delta_m^{2H^-},
\qquad 1\le k\le 2^m,\label{estimate of Emr^-1}
\end{align}
where
$C$ depends on $\sigma,b,c$ polynomially.
   \end{lemma}
\begin{proof}
	Under the assumption, $\Emr(\Ymr_{\tmkm},\theta_{\tmkm}B)^{-1}$ 
	is given by the Neumann series of 
	$
		A^{m,\rho}_{\tmkm}
		=
			I
			-
			\Emr(\Ymr_{\tmkm},\theta_{\tmkm}B)
	$.
	The estimate of the residual terms implies \eqref{estimate of Emr^-1}.
\end{proof}

\begin{remark}\label{assumption348901824091}
 When we consider the inverse $(\tJmr_t)^{-1}$, we always assume
 that $\omega\in \Omegam$ and $m$ satisfies \eqref{assumption on m}.
\end{remark}

We have the following representation of $\Zmr_t$.

\begin{lemma}
For any $t\in\Dm$ with $t>0$, we have
\begin{align}
	\Zmr_t
	&=
		\sum_{i=1}^{2^mt}
			\tJmr_{t-\tmi}(\Ymr_{\tmi},\theta_{\tmi}B)
			\left(
				c(\Ymr_{\tmim})d^m_{\tmim,\tmi}
				+\hat{\ep}^m_{\tmim,\tmi}
				-\ep^m_{\tmim,\tmi}
			\right).
			\label{rep for zmrk}
\end{align}
If all $\Zmr_s(\xi,B)$~$(s\in\Dm,0\leq s\leq t)$ are invertible, 
\begin{align*}
	\Zmr_t
	&=
		\tJmr_t
		\sum_{i=1}^{2^mt}
			(\tJmr_{\tmi})^{-1}
			\left(
				c(\Ymr_{\tmim})d^m_{\tmim,\tmi}
				+\hat{\ep}^m_{\tmim,\tmi}
				-\ep^m_{\tmim,\tmi}
			\right).
\end{align*}
\end{lemma}

\begin{proof}
The second statement follows from \eqref{Mmrk'+k} and \eqref{rep for zmrk}.
We show \eqref{rep for zmrk}.
Write $k=2^mt$ and denote by $\zeta_k$ the quantity on the right-hand side of \eqref{rep for zmrk}.
For simplicity we write
\begin{align*}
  c_{i-1}d_{i-1,i}=
  c(\Ymr_{\tmim})d^m_{\tmim,\tmi},\quad
 \ep_{i-1,i}=
 \hat{\ep}^m_{\tmim,\tmi}-
 \ep^m_{\tmim,\tmi}.
 \end{align*}
 From \eqref{EMmrk},
we have
\begin{multline*}
\zeta_k-(\zeta_{k-1}+c_{k-1}d_{k-1,k}+\ep_{k-1,k})\\
\begin{aligned}
  &=
    \sum_{i=1}^{k-1}
      \left\{
			\tJmr_{\tau^m_{k-i}}(\Ymr_{\tmi},\theta_{\tmi}B)
			-
			\tJmr_{\tau^m_{k-i-1}}(\Ymr_{\tmi},\theta_{\tmi}B)
      \right\}
      (c_{i-1}d_{i-1,i}+\ep_{i-1,i}) \\
  &=
    \sum_{i=1}^{k-1}
      \{
       \Emr(\Ymr_{\tau^m_{k-1}},\theta_{\tau^m_{k-1}}B)-I
      \}
      \tJmr_{\tau^m_{k-i-1}}(\Ymr_{\tmi},\theta_{\tmi}B)
      (c_{i-1}d_{i-1,i}+\ep_{i-1,i})\\
  &=
  \{
	\Emr(\Ymr_{\tau^m_{k-1}},\theta_{\tau^m_{k-1}}B)-I
   \}
      \sum_{i=1}^{k-1}
	  \tJmr_{\tau^m_{k-i-1}}(\Ymr_{\tmi},\theta_{\tmi}B)
        (c_{i-1}d_{i-1,i}+\ep_{i-1,i})\\
  &=
  \{
	\Emr(\Ymr_{\tau^m_{k-1}},\theta_{\tau^m_{k-1}}B)-I
   \}
      \zeta_{k-1},
\end{aligned}
\end{multline*}
which implies
\begin{align*}
  \zeta_k
  =
  \Emr(\Ymr_{\tau^m_{k-1}},\theta_{\tau^m_{k-1}}B)\zeta_{k-1}
    +c_{k-1}d_{k-1,k}+\ep_{k-1,k}.
\end{align*}
Comparing the above with \eqref{zmr identity}, we complete the proof.
\end{proof}

\begin{remark}
\begin{enumerate}
\item We do not use the notation $J^{m,\rho}_t$ to denote
the solution of \eqref{Mmr}.
The reason is as follows.
It is natural to use $(\Ymr_t,J^{m,\rho}_t)$ to denote
the interpolation process between
$(Y_t, J_t)$ and its approximate solution,
that is, we expect that $(Y^{m,0}_t,J^{m,0}_t)$ and $(Y^{m,1}_t,J^{m,1}_t)$ 
coincide $(Y_t, J_t)$ and its approximate solution, respectively.
However, $\tJmr_t$ is not such an process.
In fact, $\tilde{J}^{m,0}_t$ is not equal to $J_t$.
Differently from this, in the case of the implementable Milstein and 
Milstein schemes, $(\Ymh_t,\tilde{J}^{m,1}_t)$ is identical to
the corresponding approximate solution of $(Y_t, J_t)$.
\item When we consider quantity associated with
$\{\Ymr_t\}$,
$\{a_m\}$-order nice discrete process
$\eta$ may depend on a parameter $\rho$ $(0\le \rho\le 1)$.
For $\eta^{\rho}=\{(\eta^{m,\rho}_t)_{t\in \Dm}; m=1,2,\ldots \}$,
if we can choose the random variable $X$ in
(\ref{small order discrete process inequality}) independently of
$\rho$, we say that
$\eta^{\rho}$ is a $\{a_m\}$-order nice discrete process independent of
$\rho$.
\end{enumerate}
\end{remark}

For later use, we introduce the following.

\begin{definition}\label{definition of tzmr}
 When $\tJmr_t$ is invertible, we define
$
 \tZmr_t=(\tJmr_t)^{-1}\Zmr_t
$
for $t\in \Dm$.
Explicitly,
\begin{align}\label{def tzmr}
	\tZmr_t
	=
		\sum_{i=1}^{2^mt}
			(\tJmr_{\tmi})^{-1}
			\left(
				c(\Ymr_{\tmim})d^m_{\tmim,\tmi}
				+\hepm_{\tmim,\tmi}
				-\epm_{\tmim,\tmi}
			\right).
\end{align}
\end{definition}

\begin{proposition}\label{expression of ymh-ym}
  We assume (\ref{assumption on m}) holds.
  For any $\omega\in \Omegam$, we obtain the following neat expression
  \begin{align*}
   \Ymh_t-Y_t&=\int_0^1 \tJmr_t \tZmr_t d\rho.
  \end{align*}    
\end{proposition}

Below, we prove that under appropriate assumptions:
as $m\to\infty$,
\begin{enumerate}
  \item $\tJmr_t\to J_t$, 
        $(\tJmr_t)^{-1}\to J_t^{-1}$,
        $\Ymr_t\to Y_t$ uniformly in $t\in\Dm$ for all
        $\omega\in \Omegam$.
  \item $
          (2^m)^{2H-\frac{1}{2}}
          \sum_{i=1}^{2^mt}
		  (\tJmr_{\tmi})^{-1}
            \big(
              \hat{\ep}^m_{\tmim,\tmi}-\ep^m_{\tmim,\tmi}
            \big)
        $
        converges to $0$ uniformly in $t\in\Dm$.
\end{enumerate}
Hence it is reasonable to conjecture the main theorem holds true by Proposition~\ref{expression of ymh-ym}.
We prove our main theorem by using estimates for $\tZmr$.

\begin{remark}[List of notations]

\noindent
 \begin{itemize}
	\item	$Y_t$: solution of RDE
	\item	$\Ymh_t$: discrete approximate solution of $Y_t$ 
	\item	$\Ymr_t$: an interpolated process between
			$Y_t(=Y^{m,0}_t)$ and $\Ymh_t(=Y^{m,1}_t)$
	\item	$J_t=\partial_{\xi}Y_t(\xi,B)$
	\item	$\tJmr_t$: $\mathcal{L}(\RR^n)$-valued process
			defined by $\Ymr_t$ which approximates $J_t$
	\item	$\tJm_t=\tilde{J}^{m,0}_t$
	\item	$\Zmr_t=\partial_{\rho}\Ymr_t$
	\item	$\tZmr_t=(\tJmr_t)^{-1}\Zmr_t$
			(see Definition~\ref{definition of tzmr})
	\item	$\Emr(\Ymr_s,\theta_sB)=\tJmr_t(\tJmr_s)^{-1}$ for $t-s=\Delta_m$
			(see \eqref{defEmr} and Lemma~\ref{expansion of Mmr})
 \end{itemize}
\end{remark}

\section{Estimates of $\Ymr_t$ and $\tJmr_t$}\label{ymrk and Mmrk}

In this section,
we give estimates for $\Ymr_t$, $\tJmr_t$
and $(\tJmr_t)^{-1}$ which do not depend on
$\rho$.
Recall that $\{\Ymr_t\}_{t\in \Dm}$ satisfies $\Ymr_0=\xi$ and
\eqref{def_interplation}.
This equation is defined by the data of random variables
$d^m=\{d^m_{\tmkm,\tmk}\}_{k=1}^{2^m}$,
$\hepm=\{\hepm_{\tmkm,\tmk}\}_{k=1}^{2^m}$ $(m=1,2,\ldots)$
and $c\in C^3_b(\RR^n,L(\RR^d\otimes\RR^d,\RR^n))$.
$d^m$ and $\hepm$ need not to be 
corresponding quantities
defined in Section~\ref{2.2}
and it is not necessary that $c=(D\sigma)[\sigma]$.
Note that we define $ d^m_{s,t},\hepm_{s,t}$ for general $s,t\in\Dm$ with $s<t$
by \eqref{eq89034802914832} with 
$\eta_{\tmim,\tmi}=d^m_{\tmim,\tmi},\hepm_{\tmim,\tmi}$.
We choose $0<\lambda_1<1$ so that $\lambda_1+H^->1$ arbitrarily
and fix it.
Note that $\|d^m\|_{\lambda_1}<\infty$
because $d^m_t$ is defined on the finite set $\Dm$.

In Section~\ref{4.1}, 
for $\omega\in \Omega_0$, by applying Davie's method \cite{davie},
we give an estimate for $\Ymr_t$
in terms of the three constants $C$ given in
\eqref{eq459301490123}, \eqref{eq4839012814214}, and \eqref{eq4583902014},
and $\|d^m\|_{\lambda_1}$.

In Section~\ref{4.2}, 
we give estimates for
$\max_{t\in \Dm}\big\{|\tJmr_t|,|(\tJmr_t)^{-1}|\big\}1_{\Omegam}$.
The coefficient of the discrete RDE for which $\tJmr$ satisfies 
is not bounded but linear growth.
Hence, we cannot apply the estimate in Section~\ref{4.1}.
To overcome the difficulty, we view
the $H^{-1}$-H\"older rough path
$(B_{s,t},\BB_{s,t})$ as a rough path of finite $(H^{-})^{-1}$-variation norm.
Note that we assume Condition~\ref{condition on R} on $B_t$ 
and so we can apply the result due to 
Cass-Litterer-Lyons~\cite{cll} (see Lemma~\ref{CLL} below) to
obtain the estimate of $\tJmr$ and $(\tJmr)^{-1}$ similarly to
$J_t$ and $(J_t)^{-1}$.
In Section~\ref{4.3}, we give estimates for
$J_t-\tJm_t$ and $J^{-1}_t-(\tJm_t)^{-1}$ on $\Omega_0^{(m)}$ by using 
the results in Section~\ref{4.2}.
In Section~\ref{4.4} , we give estimates for
$\tJmr_t-J_t$ and
$(\tJmr_t)^{-1}-J^{-1}_t$.

\subsection{Estimates of $\Ymr_t$ on $\Omega_0$}\label{4.1}
For $s,t\in \Dm$ with $s\le t$,
let
\begin{align}\label{Ist def}
  I_{s,t}
  =
    \Ymr_t
    -\Ymr_s
    -\sigma(\Ymr_s)B_{s,t}
    -((D\sigma)[\sigma])(\Ymr_s)\BB_{s,t}
    -\rho c(\Ymr_s)d^m_{s,t}
    -b(\Ymr_s)(t-s).
\end{align}
First, we prove the following.

 \begin{lemma}\label{Ist small}
	Assume that Condition~\upshape{\ref{condition for epsilon}} (1) holds
	and let $\omega\in \Omega_0$.
Let $\lambda_1$ be a positive number satisfying $\lambda_1+H^->1$.
Set $\lambda=\min\{\lambda_1,2H^-\}$.
 There exist $0<\delta\le 1$ and $C_1>0$ such that
 \begin{align}
  |I_{s,t}|\le C_1|t-s|^{\lambda+H^-},\quad \quad 
	\text{$s,t\in \Dm$ \,\, with \,\, $|t-s|\le \delta$}.
\label{estimate of Ist}
 \end{align}
Here $\delta^{-1}$ and $C_1$ depend only on $\sigma, b, c$, 
$C(B)$ and $\|d^m\|_{\lambda_1}$ polynomially.
\end{lemma}

 \begin{proof}
Below, $C$ is a constant depending only on 
$\sigma$, $b$, $c$, $C(B)$ and $\|d^m\|_{\lambda_1}$ polynomially.
By using $C$, we determine $\delta$ and $C_1$ so that
\eqref{estimate of Ist} holds.
For simplicity we write $\tmi=t_i$.
Let
 $s=t_k, t=t_{k+l}$.
By
$I_{t_k,t_{k+1}}=(1-\rho)\ep^m_{t_k,t_{k+1}}
   +\rho\hat{\ep}^m_{t_k,t_{k+1}}$
and the estimate of $\hat{\ep}^m$,
we see that
\eqref{estimate of Ist} holds for any $\delta$
and for the maximum of three constants $C$ stated in
\eqref{eq459301490123}, \eqref{eq4839012814214}, and \eqref{eq4583902014}.
Let $K\ge 1$.
Suppose the following estimate:
there exists $M>0$ such that
\begin{align*}
 |I_{s,t}|&\le M|t-s|^{\lambda+H^-}
\end{align*}
holds for
 $\{(s,t)=(t_k, t_{k+l})~|~0\le k\le 2^m-1,
 l\le K,
\,\, |t-s|\le \delta
 \}.$
Here $M$ should be larger than the number $C_1$ which
is determined by the case $K=1$.

We consider the case $K+1$.
We rewrite $s=t_k$ and $t=t_{k+K+1}$.
Choose maximum $u=t_l$ satisfying
 $|u-s|\le |t-s|/2$.
Then  $|t-t_{l+1}|\le |t-s|/2$
holds.
Note that $l-k\le K$ and $K+1-(l+1)\le K$.
Hence by the assumption, we have
 \begin{align}
  \max\{|I_{s,u}|, |I_{t_{l+1},t}|\}
  &  \le M\left|\frac{t-s}{2}\right|^{\lambda+H^-},\label{estimate of Isu}\\
  \max\{|\Ymr_u-\Ymr_s|,|\Ymr_{t}-\Ymr_{t_{l+1}}|\}
  &\le
  M\left|\frac{t-s}{2}\right|^{\lambda+H^-}
  +C|t-s|^{H^{-}}.\label{estimate of Ymrsu}
 \end{align}
 Next we estimate $(\delta I)_{s,u,t}=I_{s,t}-I_{s,u}-I_{u,t}$.
 Denote by $(\delta I)_{s,u,t}^\sigma$, $(\delta I)_{s,u,t}^b$ and $(\delta I)_{s,u,t}^c$
 the terms in $(\delta I)_{s,u,t}$ being concerned with $\sigma$, $b$ and $c$, respectively.
 Then
\begin{align*}
  (\delta I)_{s,u,t}^b
  &=
    -b(\Ymr_s)(t-s)
    +b(\Ymr_s)(u-s)
    +b(\Ymr_u)(t-u)\\
  &=
    \{b(\Ymr_u)-b(\Ymr_s)\}(t-u),\\
  (\delta I)_{s,u,t}^c
  &=
   \rho
   \{c(\Ymr_u)-c(\Ymr_s)\}
   d^m_{u,t}
\end{align*}
and
\begin{align*}
  (\delta I)_{s,u,t}^\sigma
  &=
    \{\sigma(\Ymr_u)-\sigma(\Ymr_s)\}B_{u,t}
    -((D\sigma)[\sigma])(\Ymr_s)[\BB_{s,t}-\BB_{s,u}-\BB_{u,t}]\\
  &\qquad
    -\big\{((D\sigma)[\sigma])(\Ymr_s)
    -((D\sigma)[\sigma])(\Ymr_u)\big\}\BB_{u,t}\\
  &=
    \big\{\sigma(\Ymr_u)-\sigma(\Ymr_s)-D\sigma(\Ymr_s)[\Ymr_u-\Ymr_s]\big\}B_{u,t}\\
  &\qquad
    +
    D\sigma(\Ymr_s)
      [
          I_{s,u}
          +((D\sigma)[\sigma])(\Ymr_s)\BB_{s,u}
          +\rho c(\Ymr_s)d^m_{s,u}
          +b(\Ymr_s)(u-s)
      ]
      B_{u,t}\\
  &\qquad
    -
    \big\{
        ((D\sigma)[\sigma])(\Ymr_s)
        -
        ((D\sigma)[\sigma])(\Ymr_u)
    \big\}
    \BB_{u,t}.
\end{align*}
Here we used Chen's identity and definition of $I_{s,u}$.
By (\ref{estimate of Isu}) and (\ref{estimate of Ymrsu}),
we obtain
  \begin{align*}
    |(\delta I)_{s,u,t}|
    &\le
    C\big\{1+M\delta^{H^-}+(M\delta^{H^-})^2\big\}|t-s|^{\lambda+H^-}.
  \end{align*}
Similarly, we obtain
$
   |(\delta I)_{t_l,t_{l+1},t}|
   \le
   C|t-s|^{3H^{-}}
$.
By
\begin{align*}
  I_{s,t}
  &=
    I_{s,u}+I_{t_l,t_{l+1}}+I_{t_{l+1},t}
    +(\delta I)_{t_{l},t_{l+1},t}
    +(\delta I)_{s,u,t},
\end{align*}
we have $|I_{s,t}|\leq f(C,M,\delta)|t-s|^{\lambda+H^-}$, where
\begin{align*}
  f(C,M,\delta)
  =
    2^{1-(\lambda+H^-)}M
    +C\big\{1+M\delta^{H^-}+(M\delta^{H^-})^2\big\}.
\end{align*}
Note that the function $f$ and $C$ do not depend on $K$.
Let $(M,\delta)$ be a pair such that
$f(C,M,\delta)\le M$ holds and $M$ is greater than or equal to 
the maximum of three constants $C$ stated in
\eqref{eq459301490123}, \eqref{eq4839012814214}, and \eqref{eq4583902014}.
Then (\ref{estimate of Ist}) holds for
$(C_1,\delta)=(M,\delta)$.
One choice is as follows.
\begin{align*}
 M=\frac{3C}{1-2^{1-(\lambda+H^-)}},\quad 
\delta=
\min\bigg\{\bigg(\frac{3C}{1-2^{1-(\lambda+H^-)}}\bigg)^{-\frac{1}{H^-}}, 1\bigg\},
\end{align*}
where 
$C$ is greater than or equal to 
the maximum of three constants $C$ stated in
\eqref{eq459301490123}, \eqref{eq4839012814214}, and \eqref{eq4583902014}.
This completes the proof.
 \end{proof}

  \begin{lemma}\label{Ist general}
	Assume that Condition~\upshape{\ref{condition for epsilon}} (1) holds
	and let $\omega\in \Omega_0$.
Let $\lambda_1$ be a positive number satisfying $\lambda_1+H^->1$.
Set $\lambda=\min\{\lambda_1,2H^-\}$.
Then
there exist a positive number $C_2$
which depends on $\sigma,b,c$, $C(B)$ and $\|d^m\|_{\lambda_1}$
polynomially
such that
   \begin{align*}
    |I_{s,t}|\le C_2|t-s|^{\lambda+H^-},
\qquad s,t\in \Dm.
   \end{align*}
      \end{lemma}
\begin{proof}
Below, $C$ denote constants
depending only on $\sigma, b, c$, $C(B)$ and $\|d^m\|_{\lambda_1}$ polynomially.
We have proved the case where $s,t$ with $t-s\le \delta$.
Suppose $t-s>\delta$.
In this case, from the definition of $I_{s,t}$
and $(\delta^{-1}|t-s|)^\lambda \geq 1$,
we have
\begin{align*}
	|I_{s,t}|
	\leq
		|\Ymr_{s,t}|
		+C|t-s|^{H-}
	\leq
		|\Ymr_{s,t}|
		+C\delta^{-1}|t-s|^{\lambda+H^-}.
\end{align*}
Here we wrote $\Ymr_{s,t}=\Ymr_t-\Ymr_s$.
In what follows, we will give an estimates of $|\Ymr_{s,t}|$.

First, we consider the case $2^{-m}\ge \delta$.
For $s=2^{-m}k<t=2^{-m}l$, we have
\begin{align*}
	|\Ymr_{s,t}|&=\left|\sum_{i=k+1}^l\Ymr_{\tmim,\tmi}\right|
\le C (l-k)\Delta_m^{H^-}
=C (2^m)^\lambda(l-k)^{1-(\lambda+H^-)}|t-s|^{\lambda+H^-}.
\end{align*}
Noting $(2^m)^\lambda\le \delta^{-\lambda}$, we obtain
$|\Ymr_{s,t}|\le C \delta^{-\lambda}|t-s|^{\lambda+H^-}$.

We next consider the case $2^{-m}<\delta$.
	Let $\tau^m_K=\max\{\tmk~|~\tmk\le \delta\}$.
	Then $2^{-1}\delta\le \tau^m_K$.
Let $s_i=s+i\tau^m_K$ $(0\leq i\leq N-1)$,
where $N$ is a positive integer such that
$0\le t-s_{N-1}<\tau^m_K$.
For notational simplicity, we set $s_N=t$.
Then we have
$
	N\le (\tau^m_K)^{-1}(t-s)+1\le 2(t-s)(\tau^m_K)^{-1}
\le 4\delta^{-1}(t-s).
$
By the estimate in Lemma~\ref{Ist small}, we have
\begin{align*}
	|\Ymr_{s_{i-1},s_i}|
	\leq
		C
		\big\{
	|t-s|^{\lambda+H^-}+|t-s|^{H^-}+
	|t-s|^{2H^-}+|t-s|^\lambda+|t-s|
	\big\}
	\leq
		C
		|t-s|^{H^-}.
\end{align*}
Hence
\begin{align*}
	|\Ymr_{s,t}|
	\le
		\sum_{i=1}^N|\Ymr_{s_i}-\Ymr_{s_{i-1}}|
	\le 
		\delta^{-1}|t-s|
		\cdot
		C
		|t-s|^{H^-}.
\end{align*}
	Since $1>\lambda$,
	we obtain $|\Ymr_t-\Ymr_s|\le C\delta^{-1}|t-s|^{\lambda+H^-}$.
Since $\delta^{-1}$ depends on 
$\sigma$, $b$, $c$, $C(B)$, $\|d^m\|_{\lambda_1}$
polynomially, 
we complete the proof.
	\end{proof}

 For $f\in C^2_b(\RR^n, \mathcal{L}(\RR^d,\RR^K))$,
 $g\in C^1_b(\RR^n, \RR^K)$,
  and $h\in C^1_b(\RR^n,\mathcal{L}(\RR^d\otimes\RR^d,\RR^K))$,
  and  $s,t\in \Dm$ with $s<t$,
we define an $\RR^K$-valued random variable by
  \begin{align*}
   \Xi(f,g,h)_{s,t}&=
   f(\Ymr_{s})B_{s,t}
   +(Df)[\sigma](\Ymr_{s})\BB_{s,t}
   +g(\Ymr_{s})(t-s)
   +h(\Ymr_{s})d^m_{s,t},
\end{align*}
where
$
  (Df)[\sigma](y)[v\otimes w]
  =
    Df(y)[\sigma(y)v]w
$
for $y\in\RR^n$, $v,w\in\RR^d$ (see also \eqref{eqSimplifiedNotationDef}).
For a sub-partition ${\cal P}=\{u_i\}_{i=0}^l\subset \Dm$
$(s=u_0,t=u_l)$,
let
\begin{align*}
 I(f,g,h ; {\cal P})_{s,t}=
 \sum_{i=0}^l
 \Xi(f,g,h)_{u_{i-1},u_i}.
\end{align*}

\begin{lemma}\label{estimate of discrete rough integral}
	Assume that Condition~\upshape{\ref{condition for epsilon}} (1) holds
	and let $\omega\in \Omega_0$.
Let $\lambda_1$ be a positive number satisfying $\lambda_1+H^->1$.
  Set $\lambda=\min\{\lambda_1,2H^-\}$.
  Then
  \begin{align*}
   |I(f,g,h ; {\cal P})_{s,t}-\Xi(f,g,h)_{s,t}|\le C|t-s|^{\lambda+H^-},
  \end{align*}
  where $C$ depends on $\sigma,b,c,C(B),\|d^m\|_{\lambda_1}$ polynomially.
\end{lemma}

\begin{proof}
  Let $I_{st}$ be the function defined in \eqref{Ist def}.
\begin{multline*}
 \delta \Xi(f,g,h)_{s,u,t}
 =
 \Xi(f,g,h)_{s,t}-\Xi(f,g,h)_{s,u}-\Xi(f,g,h)_{u,t}\\
 \begin{aligned}
  &=
    -\left\{f(\Ymr_u)-f(\Ymr_s)
    -(Df)(\Ymr_s)[\Ymr_u-\Ymr_s]\right\}B_{u,t}\\
  &\qquad 
    -(Df)(\Ymr_s)
      \left[
        I_{s,u}
        +((D\sigma)[\sigma])(\Ymr_s)\BB_{s,u}
        +\rho c(\Ymr_{s})d^m_{s,u}+b(\Ymr_s)(u-s)
      \right]
      B_{u,t}\\
  &\qquad
    +
    \left\{
      (Df)(\Ymr_s)[\sigma(\Ymr_s)]-(Df)(\Ymr_u)[\sigma(\Ymr_u)]
    \right\}\BB_{u,t}\nonumber\\
  &\qquad
    +\left\{g(\Ymr_s)-g(\Ymr_u)\right\}(t-u)
    +\left\{h(\Ymr_s)-h(\Ymr_u)\right\}d^m_{u,t}.
 \end{aligned}
 \end{multline*}
     Hence
     $
      |\delta I(f,g,h)_{s,u,t}|\le
      C|t-s|^{\lambda+H^-}
    $.
     By a standard argument (for example, use the sewing lemma (see \cite{friz-hairer})), we complete the proof
     of the lemma.
         \end{proof}

\subsection{Estimates of $\tJmr_t$ and
$(\tJmr_t)^{-1}$ on $\Omegam$}\label{4.2}
We next proceed to the estimate of 
$\tJmr_t(\omega)$ 
and their inverse.
From now on, we always assume that $\omega\in \Omegam$ and $m$ satisfies
\eqref{assumption on m}; see Remark~\ref{assumption348901824091}.
For $\omega\in \Omega_0^{(m,d^m)}$, both estimates
$\|d^m(\omega)\|_{2H^-}\le 1$ and $\|d^m(\omega)\|_{\lambda_1}\le 1$ hold.
However note that we use one or the other only of the two estimates
in the proofs of some statements
in this section.
Since $\tJmr$ is also a solution to a discrete RDE,
one may expect similar estimates for $\tJmr$ to $\Ymr$.
However, the coefficient of the RDE of $\tJmr$ is unbounded,
we cannot apply the same proof as the one of $\Ymr$
and we need to prove the boundedness of $\tJmr$ in advance.
We give an estimate of $\tJmr$ by combining the group property
of $\tJmr$ 
and a similar argument to the estimate of $\Ymr$.
The difference from $\Ymr$ is that we use the estimate $\|d^m(\omega)\|_{2H^-}\le 1$
and the variation norm of $(B,\BB)$ (see Definition~\ref{w})
to obtain the boundedness of $\tJmr$.
After obtaining the boundedness, we see estimates on
$\tJmr_t$ and their inverse
by using the estimate $\|d^m(\omega)\|_{\lambda_1}\le 1$
and the H\"older norm of $(B,\BB)$.

First, we observe the following.
For
  $s\le t$, $s,t,\tau\in \Dm$ with
$t+\tau\le 1$, let us define
  \begin{align*}
  I_{s,t}(\Ymr_\tau,\theta_\tau B)&=
\tJmr_{t}(\Ymr_\tau,\theta_\tau B)-\tJmr_{s}(\Ymr_\tau,\theta_\tau B)
\\
&\qquad\quad-
   (D\sigma)(\Ymr_{s}(\Ymr_\tau,\theta_\tau B))
[\tJmr_{s}(\Ymr_\tau,\theta_\tau B)]
(\theta_\tau B)_{s,t}\\
&\qquad\quad   -D\left((D\sigma)[\sigma]\right)(\Ymr_{s}(\Ymr_\tau,
\theta_\tau B))
[\tJmr_{s}(\Ymr_\tau,\theta_\tau B)]
(\theta_\tau\BB)_{s,t}\\
   &\qquad\quad
   -\rho (Dc)(\Ymr_{s}(\Ymr_\tau,\theta_\tau B))
[\tJmr_{s}(\Ymr_\tau,\theta_\tau B)]
(\theta_\tau d^m)_{s,t}\\
&\qquad\quad   -(Db)(\Ymr_{s}(\Ymr_\tau ,\theta_\tau B))
[\tJmr_{s}(\Ymr_\tau ,\theta_\tau B)](t-s).
  \end{align*}
We may write 
$
I_{s,t}(\xi,B)=I_{s,t}
$
for simplicity.
  Note that
  \begin{align}
   I_{0,t-u}(\Ymr_u,\theta_uB)&=
   \tJmr_{t-u}(\Ymr_u,\theta_uB)
   -I-(D\sigma)(\Ymr_u)[I]B_{u,t}
  -D((D\sigma)[\sigma])(\Ymr_u)[I]\BB_{u,t}\nonumber\\
  &\quad -
  \rho (Dc)(\Ymr_u)[I]d^m_{u,t}-
  (Db)(\Ymr_u)[I](t-u),\label{I0t-u}
  \end{align}
  where $I$ denotes the identity operator and
  we refer the notation
  $D((D\sigma)[\sigma])(\Ymr_u)[I]\BB_{u,t}$ to
  (\ref{DDsigma2}).
  By (\ref{I0t-u}),
  if $I_{0,t-u}(\Ymr_u,\theta_uB)$ and
  $t-u$ is sufficiently small, then we see
  $\tJmr_{t-u}(\Ymr_u,\theta_uB)$ is invertible.

 \begin{lemma}\label{decomposition of Ist}
      Let $s,t,\tau, \tau' \in \Dm$ with
$\tau'\le s\le t$ and $t+\tau\le 1$.
  Then
 \begin{align*}
 I_{s,t}(\Ymr_{\tau},\theta_\tau B)&
=I_{0,t-s}(\Ymr_{s+\tau},\theta_{s+\tau}B)\tJmr_{s}
(\Ymr_{\tau},\theta_{\tau}B)\nonumber\\
&=I_{s-\tau',t-\tau'}(\Ymr_{\tau'+\tau},\theta_{\tau'+\tau} B)
\tJmr_{\tau'}
(\Ymr_{\tau},\theta_\tau B).
   \end{align*}
 \end{lemma}

\begin{proof}
These follows from the definition and the following
 identity.
 Let $u\ge s$.
\begin{align*}
  \Ymr_u(\Ymr_{\tau},\theta_{\tau}B)
  &=
    \Ymr_{u+\tau}(\xi,B)
  =
    \Ymr_{u-s}(\Ymr_{s+\tau},\theta_{s+\tau} B),\\
  \tJmr_u(\Ymr_\tau,\theta_\tau B)
  &=
    \tJmr_{u-s}(\Ymr_{s+\tau},\theta_{s+\tau}B)
    \tJmr_s(\Ymr_\tau,\theta_\tau B),\\
  (\theta_\tau \Xi)_{u,t}
  &=
  (\theta_{s+\tau}\Xi)_{u-s,t-s}
  \qquad
   \text{for}
  \qquad \Xi=B,\BB,d^m.
\end{align*}
\end{proof}

\begin{definition}\label{w}
Let $p=(H^-)^{-1}$.
For $(1,B_{s,t},\BB_{s,t})_{0\le s\le t\le 1}$, we define
\begin{align*}
 w(s,t)&=\|B\|_{[s,t], p\hyp var}^{p}+
\|\BB\|_{[s,t], \frac{p}{2}\hyp var}^{\frac{p}{2}},\qquad
0\le s\le t\le 1,
\end{align*}
where $\|~\|_{[s,t], r\hyp var}$ denotes the $r$-variation norm.
Also we define
$\tilde{w}(s,t)=w(s,t)+|t-s|$.
\end{definition}

Note that the variables $s,t$ move in $[0,1]$ and
$B$ and $\BB$ are random variables defined on $\Omega_0$
and so are $w(s,t)$ and $\tilde{w}(s,t)$.

We give estimates for $\tJmr$ and $I_{s,t}(\Ymr_\tau,\theta_\tau B)$
by using $\tilde{w}$.
First we note that the following estimate.

\begin{lemma}\label{Ist small again}
	Assume that Condition~\upshape{\ref{condition for epsilon}} (1) holds
	and let $\omega\in \Omegam$.
There exist $0<\delta\le 1$ and $C_3>0$ such that for all $s,t\in \Dm$ with
$0\le s<t\le 1$ and $\tw(s,t)\le\delta$, the following estimate holds:
\begin{align*}
	&
\left|\Ymr_t-\Ymr_s-
\sigma(\Ymr_s)B_{s,t}-((D\sigma)[\sigma])(\Ymr_s)\BB_{s,t}-
\rho c(\Ymr_s)d^m_{s,t}-b(\Ymr_s)(t-s)
\right|\\
&\qquad
\le C_3\tw(s,t)^{3H^-},
\end{align*}
where $\delta$ 
and $C_3$ are constants depending only on $\sigma, b, c, H^-$.
\end{lemma}

\begin{proof}
	The proof of this lemma is similar to that of Lemma~\ref{Ist small}
	and is done by induction.
	The difference is that we do not use \eqref{eq459301490123} and \eqref{eq4583902014}
	and use \eqref{eq8940289042} and \eqref{eq4839012814214}.
	Here we give a sketch of the proof.
	Below, $\tmi=t_i$ and
	$C$ denotes a constant depending only on 
	$\sigma$, $b$, $c$, and $H^-$ polynomially.

	The first step of the induction is as follows.
	Note
	$
		I_{t_k,t_{k+1}}
		=
			(1-\rho)\ep^m_{t_k,t_{k+1}}
			+\rho\hat{\ep}^m_{t_k,t_{k+1}}
	$.
	The estimates \eqref{eq8940289042} and \eqref{eq4839012814214}
	imply
	$
		|\epm_{t_{k-1},t_k}|
		+
		|\hepm_{t_{k-1},t_k}|
		\le
			C \tilde{w}(t_{k-1},t_k)^{3H^-}
	$
	for all
	$1\le k\le 2^m$ and $\omega\in \Omega_0^{(m)}$.
	Hence $|I_{t_k,t_{k+1}}|\leq C\tilde{w}(t_{k-1},t_k)^{3H^-}$.
	The induction works well by noting
	\begin{align*}
	|B_{s,t}|&\le \tilde{w}(s,t)^{H^-}, &
	|\BB_{s,t}|&\le \tilde{w}(s,t)^{2H^-}, & 
	|d^m_{s,t}|&\le \tilde{w}(s,t)^{2H^-}
	\qquad \text{for all}\qquad s,t\in \Dm.
	\end{align*}
	The last estimate above follows from $\omega\in\Omegam$.
	For example, we need to change 
	the sentence ``maximum $u=t_l$ satisfying $|u-s|\le |t-s|/2$''
	to ``maximum $u=t_l$ satisfying $\tw(s,u)\le \tw(s,t)/2$''.
	For this $l$, we see $\tw(t_{l+1},t)\le \frac{1}{2}\tw(s,t)$.
	We omit the details.
\end{proof}

  \begin{lemma}\label{Ist small2}
	Assume that Condition~\upshape{\ref{condition for epsilon}} (1) holds
	and let $\omega\in \Omegam$.
   There exist $0<\delta\le 1$ and $C_4>0$
   such that for any
$t,\tau\in \Dm$ with
$\tw(\tau,\tau+t)\le \delta$ and
$t+\tau\le 1$, the following estimate holds.
   \begin{align}
|I_{0,t}(\Ymr_{\tau},\theta_\tau B)|\le
    C_4 \tw(\tau,\tau+t)^{3H^{-}},\label{Ist}
       \end{align}
where 
$\delta$ and $C_4$ are constants depending only on 
$\sigma, b, c, H^{-}$.
  \end{lemma}

    \begin{proof}
Below, we write $\tilde{w}_{\tau}(s,t)=\tilde{w}(s+\tau,t+\tau)$
and $C$ is a constant depending only on
     $\sigma,b,c, H^-$ which may change line by line.
The proof is similar to that of Lemma~\ref{Ist small}.
     We take $\delta$ smaller than $\delta$ in
     Lemma~\ref{Ist small again}.
    For simplicity we write $t_k=\tmk$.
  It suffices to consider the case where $\tau\le 1-2^{-m}$.
We consider the following claim depending on a positive integer
$K$.

\medskip

\noindent
(Claim $K$)~(\ref{Ist}) holds for all $\tau$ and $t_k$ satisfying
$\tau+t_k\le 1$, $\tw_\tau(0,t_k)\le \delta$ and $1\le k\le K$.

\medskip
Since
$
 I_{0,t_1}=
I_{0,t_1}(\Ymr_{\tau},\theta_{\tau}B)=0
$
holds for all $\tau$, 
(Claim 1) holds for $C_4=0$ and any $\delta$.
We assume (Claim $K$) holds and we will find the condition
on $C_4$ and $\delta$ independent of $K$ under which (Claim $K+1$) holds.
Assume the case $K$ holds for a positive constant $C_4$ and $\delta$.
Suppose $\tau+t_{K+1}\le 1$ and 
$\tw_\tau(0,t_{K+1})\le \delta$, where $K\ge 1$.
Define $0\le t_l<t_{K+1}$ as the maximum number 
such that $\tw_{\tau}(0,t_l)\le \tw_{\tau}(0,t_{K+1})/2$.
On the other hand, for $t_{l+1}$, 
we have $\tw_{\tau}(t_{l+1},t_{K+1})\le \tw_{\tau}(0,t_{K+1})/2$.
We will write $u=t_l$ and $t=t_{K+1}$.
By (Claim $K$), we have
\begin{align}
  \label{Itl}
  |I_{0,u}(\Ymr_\tau,\theta_\tau B)|
  &\le
    C_4(\tw_{\tau}(0,t)/2)^{3H^{-}},\\
  \label{Itl2}
  |I_{0,t-t_{l+1}}(\Ymr_{t_{l+1}+\tau},\theta_{t_{l+1}+\tau}B)|
  &\le 
    C_4(\tw_{\tau}(0,t)/2)^{3H^{-}}.
\end{align}
The estimate (\ref{Itl}) implies
\begin{gather}
  \label{JI}
  \begin{aligned}
    |\tJmr_{u}(\Ymr_\tau,\theta_\tau B)-I|
    &\le
      C_4(\tw_{\tau}(0,t)/2)^{3H^{-}}
      +C\tw_{\tau}(0,t)^{H^{-}}
      +C\tw_{\tau}(0,t)^{2H^{-}} \\
    &\le
      \{C_4(\delta/2)^{2H^{-}}+C\}
      \tw_{\tau}(0,t)^{H^{-}},
  \end{aligned}\\
  \label{JIDsigma}
  |
    \tJmr_{u}(\Ymr_\tau,\theta_\tau B)
    -I
    -(D\sigma)(\Ymr_{\tau})B_{\tau,u+\tau}
  |
  \le
      \{C_4(\delta/2)^{H^{-}}+C\}
      \tw_{\tau}(0,t)^{2H^{-}}.
\end{gather}
For simplicity, we write $I_{0,t}=I_{0,t}(\Ymr_\tau,\theta_\tau B)$
and set $(\delta I)_{0,u,t}=I_{0,t}-I_{0,u}-I_{u,t}$.
Hereafter we will estimate $(\delta I)_{0,u,t}$ and $I_{u,t}$.
By the results on them and the inductive assumption, we will obtain a bound of $I_{0,t}$

First we consider $(\delta I)_{0,u,t}$.
Denote by $(\delta I)_{0,u,t}^\sigma$, $(\delta I)_{0,u,t}^b$ and $(\delta I)_{0,u,t}^c$
the terms in $(\delta I)_{0,u,t}$ being concerned with $\sigma$, $b$ and $c$, respectively.
Then we have
\begin{align*}
  (\delta I)_{0,u,t}^b
  &=
    -(Db)(\Ymr_\tau)[I]t
    +(Db)(\Ymr_\tau)[I]u\\
  &\phantom{=}\qquad
    +(Db)(\Ymr_u(\Ymr_\tau,\theta_\tau B))[\tJmr_{u}(\Ymr_\tau,\theta_\tau B)](t-u)\\
  &=
    \big\{
        (Db)(\Ymr_{u+\tau})[\tJmr_{u}(\Ymr_\tau,\theta_\tau B)]
        -
        (Db)(\Ymr_\tau)[I]
    \big\}
    (t-u)\\
  (\delta I)_{0,u,t}^c
  &=
    \rho
    \big\{
        (Dc)(\Ymr_{u+\tau})[\tJmr_u(\Ymr_{\tau},\theta_{\tau}B)]
        -
        (Dc)(\Ymr_\tau)[I]
    \big\}
    d^m_{u+\tau,t+\tau}
\end{align*}
and
\begin{align*}
  (\delta I)_{0,u,t}^\sigma
  &=
    -(D\sigma)(\Ymr_{\tau})[I]B_{u+\tau,t+\tau}
    -
    D((D\sigma)[\sigma])(\Ymr_{\tau})[I]
      \left(
          \BB_{\tau,\tau+t}-\BB_{\tau,\tau+u}
      \right) \\
  &\qquad\qquad
    +(D\sigma)(\Ymr_{u+\tau})[\tJmr_{u}(\Ymr_\tau,\theta_\tau B)] B_{u+\tau,t+\tau} \\
  &\qquad\qquad
    +
    D((D\sigma)[\sigma])(\Ymr_{u+\tau})[\tJmr_{u}(\Ymr_\tau,\theta_\tau B)] \BB_{u+\tau,t+\tau}.
\end{align*}
Here by getting the first and third terms together, we have
\begin{multline*}
  (D\sigma)(\Ymr_{u+\tau})
  \Big[
      \tJmr_{u}(\Ymr_\tau,\theta_\tau B)-I-D\sigma(\Ymr_\tau)[I]B_{\tau,\tau+u}
  \Big] B_{u+\tau,t+\tau}\\
  +
  \Big\{
    (D\sigma)(\Ymr_{u+\tau})[I]
    -
    (D\sigma)(\Ymr_{\tau})[I]
    -
    D(D\sigma)(\Ymr_\tau)[\sigma(\Ymr_{\tau})B_{\tau,u+\tau}]
  \Big\}
  B_{u+\tau,t+\tau}\\
  +
  \uwave{
    (D\sigma)(\Ymr_{u+\tau})
    \left[
      D\sigma(\Ymr_\tau)[I]B_{\tau,\tau+u}
    \right]
    B_{u+\tau,t+\tau}  
  }
  +
  \uwave{D(D\sigma)(\Ymr_\tau)[\sigma(\Ymr_{\tau})B_{\tau,u+\tau}]B_{u+\tau,t+\tau}}.
\end{multline*}
Because of Chen's identity, the summation of the second and fourth terms gives
\begin{multline*}
  \Big\{
    D((D\sigma)[\sigma])(\Ymr_{u+\tau})[\tJmr_{u}(\Ymr_\tau,\theta_\tau B)]
    -
    D((D\sigma)[\sigma])(\Ymr_{\tau})[I]
  \Big\}
  \BB_{u+\tau,t+\tau}\\
  +
  \uwave{
    (
      -
      D((D\sigma)[\sigma])(\Ymr_{\tau})[I]
      \left\{
        B_{\tau,\tau+u}\otimes B_{\tau+u,\tau+t}
      \right\}
    )
  }.
\end{multline*}
Since the summation of terms with\, \uwave{\phantom{aaaa}}\, vanishes due to \eqref{DDsigma2},
we have
\begin{align*}
  (\delta I)_{0,u,t}^\sigma
  &=
    (D\sigma)(\Ymr_{u+\tau})
      \Big[\tJmr_{u}(\Ymr_\tau,\theta_\tau B)-I-(D\sigma)(\Ymr_{\tau})B_{\tau,u+\tau}\Big]
      B_{u+\tau,t+\tau}\\
  &\qquad
    +
      \Big\{
          (D\sigma)(\Ymr_{u+\tau})
          -
          (D\sigma)(\Ymr_{\tau})-D(D\sigma)(\Ymr_\tau)[\sigma(\Ymr_{\tau})B_{\tau,u+\tau}]
      \Big\}
      B_{u+\tau,t+\tau} \\
  &\qquad
    +
    \Big\{
      D((D\sigma)[\sigma])(\Ymr_{u+\tau})[\tJmr_{u}(\Ymr_\tau,\theta_\tau B)]
      -
      D((D\sigma)[\sigma])(\Ymr_{\tau})[I]
    \Big\}
    \BB_{u+\tau,t+\tau}.
\end{align*}
    Thus, combining Lemma~\ref{Ist small again},
    (\ref{JI}) and (\ref{JIDsigma}) ,
    we get
\begin{align*}
  |(\delta I)_{0,u,t}^\sigma|
  &\le
    C\tw_{\tau}(0,t)^{3H^{-}}
    +
    C\big\{ 1 + C_4(\delta/2)^{H^-} \big\}\tw_{\tau}(0,t)^{3H^-},\\
  |(\delta I)_{0,u,t}^b|
  &\le
    C\big\{ 1 + C_4(\delta/2)^{2H^-} \big\}\tw_{\tau}(0,t)^{1+H^-},\\
  |(\delta I)_{0,u,t}^c|
  &\le
    C\big\{ 1 + C_4(\delta/2)^{2H^-} \big\}\tw_{\tau}(0,t)^{3H^-}.
 \end{align*}
Hence,
\begin{align*}
  |(\delta I)_{0,u,t}|
  &\le
    C\big\{1+C_4\delta^{H^-}\big\}\tw_{\tau}(0,t)^{3H^-}.
\end{align*}

We estimate $I_{u,t}$.
We have
$
  I_{u,t}=I_{t_l,t}
  =
    (\delta I)_{t_l,t_{l+1},t}+
    I_{t_l,t_{l+1}}+I_{t_{l+1},t}
$.
It is clear that $I_{t_l,t_{l+1}}=0$.
First we consider $(\delta I)_{t_l,t_{l+1},t}$.
Using Lemma~\ref{decomposition of Ist} and (\ref{JI}), we get
\begin{align*}
  |(\delta I)_{t_l,t_{l+1},t}|
  &=
    \big|
      \big\{
        I_{0,t-t_l}(\Ymr_{t_l+\tau},\theta_{t_l+\tau}B)
        -I_{0,t_{l+1}-t_l}(\Ymr_{t_l+\tau}, \theta_{t_l+\tau}B)\\
  &\qquad\qquad\qquad
        -I_{t_{l+1}-t_l,t-t_l}(\Ymr_{t_l+\tau},\theta_{t_l+\tau}B)
      \big\}
      \cdot
      \tJmr_{t_l}(\Ymr_\tau,\theta_\tau B)
    \big|\\
 &\le
    C
    \big\{
        1 + C_4\delta^{H^-}
    \big\}
    \tw_{\tau+t_l}(0,t-t_l)^{3H^-}
    \big|\tJmr_{t_l}(\Ymr_\tau,\theta_\tau B)\big|,
\end{align*}
where we have used a similar estimate of
$(\delta I)_{0,t_{l+1}-t_l,t-t_l}$ to $(\delta I)_{0,u,t}$
and note $\tw_{\tau+t_l}(0,t-t_l)=\tw_{\tau}(t_l,t)\leq \tw_{\tau}(0,t)$.
Next we consider $I_{t_{l+1},t}$.
Lemma~\ref{decomposition of Ist} implies
\begin{align*}
  I_{t_{l+1},t}
  &=
    I_{0,t-t_{l+1}}(\Ymr_{t_{l+1}+\tau},\theta^m_{t_{l+1}+\tau}B)
    \tJmr_{t_{l+1}}(\Ymr_{\tau},\theta_\tau B)\\
  &=
    I_{0,t-t_{l+1}}(\Ymr_{t_{l+1}+\tau},\theta^m_{t_{l+1}+\tau}B)
    \Emr(\Ymr_{t_l+\tau},\theta_{t_l+\tau}B)
    \tJmr_{t_l}(\Ymr_\tau,\theta_\tau B).
\end{align*}
By \eqref{Itl2} and the definition of $\Emr$ (see \eqref{defEmr}), we obtain
\begin{align*}
  |I_{t_{l+1},t}|
  &\le
    C_4
    \left(
      \frac{1}{2}
      \tw_{\tau}(0,t)
    \right)^{3H^{-}}
    \big\{ 1 + C\tw_\tau(0,t)^{H^-} \big\}
    \big|\tJmr_{t_l}(\Ymr_\tau,\theta_\tau B)\big|.
\end{align*}
Hence noting 
$
  |\tJmr_{t_l}(\Ymr_\tau,\theta_\tau B)|
  \leq 
    1+C\{1+C_4\delta^{H^-}\}
$,
we have
\begin{align*}
  |I_{u,t}|
  &\leq
    \big\{
      C
      \big\{
        1+C_4\delta^{H^-}
      \big\}
      +
      C_4
      2^{-3H^-}
      \big\{ 1+C\delta^{H^-} \big\}
    \big\}
    \big\{
      1+C\{ 1 + C_4\delta^{H^-} \}
    \big\}
    \tw_\tau(0,t)^{3H^-}\\
  &\leq
    \big\{
      C_4
      2^{-3H^-}
      +
      C
      \big\{
        1 + C_4\delta^{H^-}
      \big\}
    \big\}
    \big\{
      1+C\{ 1 + C_4\delta^{H^-} \}
    \big\}
    \tw_\tau(0,t)^{3H^-}\\
  &\leq
    \big\{
      C_4
      2^{-3H^-}
      +
      C
      \big\{
        1+C_4\delta^{H^-}+(C_4\delta^{H^-})^2
      \big\}
    \big\}
    \tw_\tau(0,t)^{3H^-}.
\end{align*}

Consequently, noting $I_{0,t}=I_{0,u}+(\delta I)_{0,u,t}+I_{u,t}$, we obtain
\begin{align*}
  |I_{0,t}|
  &\le
    \big\{
      2C_4 2^{-3H^{-}}
      +
      C
      \big\{1+(C_4\delta^{H^-})+(C_4\delta^{H^-})^2\big\}
  \big\}\tw_{\tau}(0,t)^{3H^{-}}.
\end{align*}
Hence if $C_4$ and $\delta$ satisfies
$
  C_4 2^{1-3H^{-}}
  +
  C
  \big\{
      1+(C_4\delta^{H^-})+(C_4\delta^{H^-})^2
  \big\}
  \le
    C_4
$,
then (\ref{Ist}) holds in the case of $K+1$.
One choice of $C_4, \delta$ is 
\begin{align*}
  C_4
  &=
    \frac{3C}{1-2^{1-3H^-}},
  &
  \delta
  &=
    \min
      \bigg\{
        \bigg(\frac{3C}{1-2^{1-3H^-}}\bigg)^{-\frac{1}{H^-}}, 1
      \bigg\}.
\end{align*}

Under this choice, we see that (\ref{Ist}) holds for any
$t,\tau\in \Dm$ with $\tw(\tau,\tau+t)\le \delta$ and $t+\tau\le 1$.
This completes the proof.
    \end{proof}

In order to obtain $L^p$ estimate in Theorem~\ref{main theorem},
we need the estimate obtained by Cass-Litterer-Lyons~\cite{cll}.
To this end, we introduce the number $N_{\beta}(w)$ which is defined for
any control function $w$ and positive number $\beta$.
We already used the notation $w$ in Definition~\ref{w} and so
this is an abuse in a certain sense.
For a control function $w$ and a positive number $\beta$, let us define 
$N_{\beta}(w)$ and a nondecreasing sequence
$\{\sigma_i\}_{i=0}^{\infty}\subset [0,1]$
as follows.
\begin{enumerate}
  \item $\sigma_0=0$.
  \item Let $i\ge 0$ and write $A_i=\{s\in [0,1]~|~s\ge \sigma_{i},w(\sigma_{i},s)\ge \beta\}$.
        Set $\sigma_{i+1}=\inf A_i$ (resp. $1$) if $A_i\neq \emptyset$ (resp. $A_i=\emptyset$).
  \item  $N_{\beta}(w)=\sup\{i\ge 0~|~\sigma_i<1\}$.
\end{enumerate}

We have the following.

\begin{lemma}\label{number N}
Let $w, w'$ be any control functions and
$\beta,\beta'>0$.
\begin{enumerate}
  \item There exist positive 
constants $C_{\beta,\beta'}, C'_{\beta,\beta'}$ 
which are independent of $w$ such that
\[
 C_{\beta,\beta'}(N_{\beta'}(w)+1)\le N_{\beta}(w)+1\le 
C'_{\beta,\beta'}(N_{\beta'}(w)+1).
\]
\item If $w(s,t)\le w'(s,t)$ $(0\le s\le t\le 1)$ holds,
then $N_{\beta}(w)\le N_{\beta}(w')$.
\item Let $\tilde{w}(s,t)=w(s,t)+|t-s|$ $(0\le s\le t\le 1)$.
Then for any $\beta\ge 3$, we have
$N_{\beta}(\tilde{w})\le N_1(w)$.
 \end{enumerate}
\end{lemma}

\begin{proof}
	We show (1). We use $\sigma^\beta_i$ to denote the dependence of $\sigma_i$ on $\beta$.
	Assume $\beta'<\beta$. Then $\sigma^{\beta'}_i\leq \sigma^\beta_i$ for all $i\geq 0$,
	which implies $N_{\beta'}(w)\geq N_\beta(w)$.
	Conversely, by setting 
	$
		\Lambda_i
		=
			\{j:\sigma^\beta_i\leq \sigma^{\beta'}_j,\sigma^{\beta'}_{j+1}\leq \sigma^\beta_{i+1}\}
	$
	for $0\leq i\leq N_\beta(w)-1$,
	we have
	\begin{align*}
		\beta 
		=
			w(\sigma^\beta_i,\sigma^\beta_{i+1})
		\geq
			\sum_{j\in \Lambda_i}
				w(\sigma^{\beta'}_j,\sigma^{\beta'}_{j+1})
		=
			\sharp \Lambda_i
				\beta'.
	\end{align*}
	Since the number of $j$ such that $\sigma^\beta_i\in(\sigma^{\beta'}_j,\sigma^{\beta'}_{j+1})$
	for some $1\leq i\leq N_\beta(w)$ is bounded by $N_\beta(w)$ from above
	and the number of $j$ such that $(\sigma^{\beta'}_j,\sigma^{\beta'}_{j+1}]\subset(\sigma^\beta_{N_\beta(w)},1]$
	is bounded by $\beta/\beta'$,
	we have $\sum_{i=0}^{N_\beta(w)-1} \sharp \Lambda_i \geq  N_{\beta'}(w)-N_\beta(w)-\beta/\beta'$.
	Hence $\beta N_\beta(w)\geq \beta'(N_{\beta'}(w)-N_\beta(w)-\beta/\beta')$.
	Hence we see the assertion for $\beta'<\beta$. It can be generalized easily.
 	We can show (2) easily from the definition.
	We prove (3).
Let $\{\tilde{\sigma}_i\}_{i=0}^{N_{\beta}(\tilde{w})}$ 
and $\{\sigma_i\}_{i=0}^{N_1(w)}$ be corresponding increasing
sequences.
Then by the definition, we have
$
 w(\tilde{\sigma}_{i-1},\tilde{\sigma}_i)\ge 2
$ for $1\le i\le 
N_{\beta}(\tilde{w})$.
This implies $\sigma_i\le \tilde{\sigma}_i$~~
$(1\le i\le N_{\beta}(\tilde{w}))$ and so the proof is finished.
\end{proof}

In what follows, we write
\begin{align*}
  \tilde{N}(B)=2^{N_{\beta}(\tilde{w})+1}.
\end{align*}

\begin{lemma}\label{estimate for tJmr}
	Assume that Condition~\upshape{\ref{condition for epsilon}} (1) holds
	and let $\omega\in \Omegam$.
There exist a positive integer $m_0$ and a positive number $\beta$
which depend only on
$\sigma,b,c, H^-$ such that for all $m\ge m_0$ it holds that
$\tJmr_t$ are invertible for all $t\in \Dm$ and
\begin{align*}
 \max_{t\in \Dm}\big\{|\tJmr_t|,|(\tJmr_t)^{-1}|\big\}\le \tilde{N}(B).
\end{align*}
\end{lemma}

\begin{proof}
Let $\delta$ and $C_4$ be numbers given in Lemma~\ref{Ist small2}.
Let us take $m$ satisfying $2^{-m}\le \delta$.
Let $0<\vep\le\delta$.
By Lemma~\ref{Ist small2}, for $t,\tau$ satisfying
$\tw(\tau,\tau+t)\le\vep$ and $\tau+t\le 1$,
we have
\begin{align*}
 |\tJmr_{t}(\Ymr_{\tau},\theta_{\tau}B)-I|&\le
C_4\vep^{3H^-}+C(\vep^{H^-}+\vep^{2H^-}+\vep),
\end{align*}
where $C$ is a constant depending only on $\sigma, b,c$.
Hence, for sufficiently small $\vep$ which depends only on
$C_4, C$, that is, depends only on $\sigma,b,c$,
it holds that for any $t,\tau\in \Dm$ with $t+\tau\le 1$ and
$\tw(\tau,t+\tau)\le \vep$,
$\tJmr_t(\Ymr_{\tau},\theta_{\tau}B)$ are invertible 
and
\begin{align}
 \max\big\{|\tJmr_{t}(\Ymr_{\tau},\theta_{\tau}B)|,
|\tJmr_{t}(\Ymr_{\tau},\theta_{\tau}B)^{-1}|\big\}
&\le 2.\label{small tJmr}
\end{align}

By the definition of $w$, we see that there exists
a constant $C_{H^-}(\ge 1)$ such that for any $0\le s<u<t\le 1$
\begin{align*}
 w(s,t)\le C_{H^-}\left(w(s,u)+w(u,t)\right).
\end{align*}
For $\omega\in \Omega_0^{(m)}$,
$w\left(u,(u+2^{-m})\wedge 1\right)\le 2^{-m}$ 
holds for any $0\le u\le 1$.
Therefore, we get
\begin{align*}
 w\left(s,(u+2^{-m})\wedge 1\right)\le C_{H^-}\left(w(s,u)+2^{-m}\right),
\quad 0\le s\le u\le 1.
\end{align*}
By using this, we get 
\begin{align*}
 \tilde{w}\left(s,(u+2^{-m})\wedge 1\right)
\le C_{H^-}\left(\tilde{w}(s,u)+2^{1-m}\right),
\quad 0\le s\le u\le 1.
\end{align*}
Let us take a positive number $\beta$ and $m$ such that
\begin{align*}
 C_{H^-}\left(\beta+2^{1-m}\right)\le \vep.
\end{align*}
Note that $\beta$ and $m$ depends on $C_{H^-}$ and $\vep$.
Let $\{\tilde{\sigma}_i\}_{i=0}^{N_{\beta}(\tilde{w})}$ be the increasing
sequence defined by $\tilde{w}$ and $\beta$.
Let $\hat{\sigma}_i=\inf\{t\in \Dm~|~t\ge \tilde{\sigma}_i\}$
$(0\le i\le N_{\beta}(\tilde{w}))$.
Also set $\hat{\sigma}_{N_{\beta}(\tilde{w})+1}=1$.
Then we have for all $0\le i\le N_{\beta}(\tilde{w})$
\begin{align}
 \tilde{w}(\hat{\sigma}_{i},\hat{\sigma}_{i+1})
\le \tilde{w}\left(\tilde{\sigma}_{i},
(\tilde{\sigma}_{i+1}+2^{-m})\wedge 1\right)
\le C_{H^-}(\tilde{w}(\tilde{\sigma}_{i},\tilde{\sigma}_{i+1})+2^{1-m})
\le \vep.\label{sigmai and sigmai+1}
\end{align}
Take $t(\ne 0)\in \Dm$ and choose $j$ so that
$\hat{\sigma}_{j-1}<t\le \hat{\sigma}_j$~$(1\le j\le N_{\beta}(\tilde{w})+1)$.
We have
\begin{align}
 \tJmr_t(\xi,B)=
\tJmr_{t-\hat{\sigma}_{j-1}}(\Ymr_{\hat{\sigma}_{j-1}},
\theta_{\hat{\sigma}_{i-1}}B)\cdots
\tJmr_{\hat{\sigma}_2-\hat{\sigma}_1}(\Ymr_{\hat{\sigma}_1},
\theta_{\hat{\sigma}_1}B)
\tJmr_{\hat{\sigma}_1}(\xi,B).
\label{product of tJmr1}
\end{align}
By (\ref{small tJmr}), (\ref{sigmai and sigmai+1}) and
(\ref{product of tJmr1}),
We obtain
\begin{align*}
\max_{t\in \Dm}\big\{
|\tJmr_t(\xi,B)|, |\tJmr_t(\xi,B)^{-1}|
\big\}\le 2^{N_{\beta}(\tilde{w})+1},
\end{align*}
which completes the proof.
\end{proof}

 \begin{lemma}\label{R estimate for tJmr}
	Assume that Condition~\upshape{\ref{condition for epsilon}} (1) holds
	and let $\omega\in \Omegam$.
Set $\lambda=\min\{\lambda_1,2H^-\}$.
Let $m$ be a sufficiently large number as in
Lemma~\upshape{\ref{estimate for tJmr}}.
There exists a positive number $C_5$ which 
does not depend on $m$ and depends on
$\tilde{C}(B)$ and $\tilde{N}(B)$ polynomially
such that, for all $t,s\in \Dm$,
 \begin{multline}
\big|\tJmr_{t}-\tJmr_{s}
-(D\sigma)(\Ymr_{s})
[\tJmr_{s}]B_{s,t}
-D\left((D\sigma)[\sigma]\right)(\Ymr_{s})[\tJmr_{s}]\BB_{s,t}\\
-\rho (Dc)(\Ymr_{s})[\tJmr_{s}]d^m_{s,t}
-(Db)(\Ymr_{s})[\tJmr_{s}](t-s)\big|\le
C_5|t-s|^{\lambda+H^-}.\label{estimate of remainder term of tJmr}
 \end{multline}
\end{lemma}

\begin{proof}
We already proved that there exists $\tilde{N}(B)$ such that
$|\tJmr_t|\le \tilde{N}(B)$ for all sufficiently large $m$ and
$t\in \Dm$.
Noting this boundedness, we obtain desired result by the same proofs 
as in Lemmas~\ref{Ist small} and~\ref{Ist general}.
\end{proof}

$(\tJmr_t)^{-1}$ also satisfies a similar estimate.

\begin{lemma}\label{R estimate for tJmr^-1}
	For every $s,t\in \Dm$ with $s\le t$, set
	\begin{align*}
		\tilde{A}^{m,\rho}_{s,t}
		&=
			-
			\Big[
				(D\sigma)(\Ymr_{s})B_{s,t}\\
		&\qquad\qquad
			+
				\sum_{\alpha,\beta}
					\Bigl\{
						(D\sigma)(\Ymr_{s})
						[(D\sigma)(\Ymr_{s})e_{\beta}]e_{\alpha}-
						(D^2\sigma)(\Ymr_{s})
						[\cdot,\sigma(\Ymr_{s})e_{\alpha}]e_{\beta}
					\Bigr\}
					B^{\alpha,\beta}_{s,t}\\
		&\qquad\qquad
			+
			\rho (Dc)(\Ymr_{s})d^m_{s,t}
			+
			(Db)(\Ymr_{s})(t-s)
		\Big].
	\end{align*}
	Assume that Condition~\upshape{\ref{condition for epsilon}} (1) holds
	and let $\omega\in \Omegam$.
Set $\lambda=\min\{\lambda_1,2H^-\}$.
Let $m$ be a sufficiently large number as in
Lemma~\upshape{\ref{estimate for tJmr}}.
\begin{enumerate}
	\item We define $\tilde{\ep}^{m,\rho}_{\tmim,\tmi}$ by
	$
	\tepmr_{\tmim,\tmi}
	=
		(\tJmr_{\tmi})^{-1}
		-(\tJmr_{\tmim})^{-1}
		-(\tJmr_{\tmim})^{-1}\tilde{A}^{m,\rho}_{\tmim,\tmi}
	$.
Then it holds that
\begin{align}
& |\tepmr_{\tmim,\tmi}|\le
2\tilde{N}(B)
\Bigl(1+\|D\sigma\|+\|D\left((D\sigma)[\sigma]\right)\|+
\|Dc\|+\|Db\|\Bigr)^3\Delta_m^{\lambda+H^-}.
\end{align}
\item For all $s,t\in \Dm$ with $s\le t$, it holds that
there exists a constant $C_6$ which is defined by a polynomial function
of $\tilde{C}(B)$ and $\tilde{N}(B)$ such that
\begin{align}
\big|(\tJmr_{t})^{-1}-(\tJmr_{s})^{-1}
-(\tJmr_{s})^{-1}\tilde{A}^{m,\rho}_{s,t}\big|\le C_6|t-s|^{\lambda+H^-}.
\label{estimate of remainder term of tJmr^-1}
\end{align}
\end{enumerate}
\end{lemma}

\begin{proof}
	(1) Set
	$
		A^{m,\rho}_{\tmim,\tmi}
		=
		I-\Emr(\Ymr_{\tmim},\theta_{\tmim}B)
	$.
	By the equation \eqref{expansion of Mmr2}, we have
\begin{multline*}
	(\tJmr_{\tmi})^{-1}-(\tJmr_{\tmim})^{-1}
		=
			(\tJmr_{\tmim})^{-1}
			\left(\Emr(\Ymr_{\tmim},\theta_{\tmim}B)^{-1}-I\right)\\
		=
			(\tJmr_{\tmim})^{-1}
			\left((I-A^{m,\rho}_{\tmim,\tmi})^{-1}-I\right)
	=
		(\tJmr_{\tmim})^{-1}
		\bigg[
			A^{m,\rho}_{\tmim,\tmi}
			+
			\sum_{l=2}^{\infty}
				\left\{
					A^{m,\rho}_{\tmim,\tmi}
				\right\}^l
		\bigg].
\end{multline*}
By the geometric property 
$B^{\alpha,\beta}_{s,t}=B^{\alpha}_{s,t}B^{\beta}_{s,t}-B^{\beta,\alpha}_{s,t}$,
we have
	\begin{align*}
&  (D\sigma)(\Ymr_s)[(D\sigma)(\Ymr_s)B_{s,t}]
	B_{s,t}-
	(D\sigma)(\Ymr_s)[(D\sigma)(\Ymr_s)]\BB_{s,t}
	\nonumber\\
	&\quad=
	(D\sigma)(\Ymr_s)[(D\sigma)(\Ymr_s)e_{\alpha}]
	e_{\beta}B^{\alpha}_{s,t}B^{\beta}_{s,t}-
	(D\sigma)(\Ymr_s)[(D\sigma)(\Ymr_s)e_{\alpha}]e_{\beta}
	\BB^{\alpha,\beta}_{s,t}\nonumber\\
	&\quad=(D\sigma)(\Ymr_s)[(D\sigma)(\Ymr_s)e_{\alpha}]e_{\beta}
	\BB^{\beta,\alpha}_{s,t}.
	\end{align*}
Using this and by the assumption of (\ref{assumption on m})
and Lemma~\ref{estimate for tJmr},
we obtain the desired estimate.

\noindent
(2) We have proved that $(\tJmr_t)^{-1}$ satisfies a similar equation
to $\Ymr_t$ and the norm can be estimated as in 
Lemma~\ref{estimate for tJmr}.
Hence, we can complete the proof in the same way as in
Lemma~\ref{Ist general}.
\end{proof}

We now give an estimate of discrete rough integral
similarly to Lemma~\ref{estimate of discrete rough integral}.

\begin{lemma}\label{estimate of discrete rough integral 2}
 Let $\varphi$ be a
$C^2_b$ function on
$\RR^n\times \mathcal{L}(\RR^n)\times 
\mathcal{L}(\RR^n)$ with values in $\mathcal{L}(\RR^d,\RR^l)$
whose all derivatives and itself are 
at most polynomial order growth.
For $t\in\Dm$, set
\begin{multline*}
	I^{m,\rho}(\varphi)_t\\
	=
		\sum_{i=1}^{2^mt}
			\Big\{
				\varphi
				\left(
					\Ymr_{\tmim},
					\tJmr_{\tmim},
					(\tJmr_{\tmim})^{-1}
				\right)
				B_{\tmim,\tmi}
			+
			\varphi
				\left(
					\Ymr,
					\tJmr,
					(\tJmr)^{-1}
				\right)^{\bcdot}_{\tmim}
				\BB_{\tmim,\tmi}
			\Big\},
\end{multline*}
where
$
	\varphi
	(
		\Ymr,
		\tJmr,
		(\tJmr)^{-1}
	)^{\bcdot}_t
$
$(t\in\Dm)$
is the $\mathcal{L}(\RR^d\otimes\RR^d,\RR^l)$-valued process
such that
\begin{multline*}
	\varphi
	\left(
		\Ymr,
		\tJmr,
		(\tJmr)^{-1}
	\right)^{\bcdot}_t
	[v \otimes w]
  =
    (D_1\varphi)
	\left(
		\Ymr_t,
		\tJmr_t,
		(\tJmr_t)^{-1}
	\right)
	\left[\sigma\left(\Ymr_t\right)v\right]w \\
  \begin{aligned}
    &
    +
      (D_2\varphi)
	  \left(
		\Ymr_t,
		\tJmr_t,
		(\tJmr_t)^{-1}
		\right)
		\left[
			(D\sigma)\left(\Ymr_t\right)
			\left[\tJmr_t\cdot\right]v
		\right]
		w \\
    &
    -
      (D_3\varphi)
	  \left(
			\Ymr_t,
			\tJmr_t,
			(\tJmr_t)^{-1}
		\right)
		\left[
			(\tJmr_t)^{-1}
		    (D\sigma)\left(\Ymr_t\right)[\cdot]v
		\right]w
  \end{aligned}
\end{multline*}
for $v,w\in\RR^d$.
Here
$D_i$ denotes the derivative with respect to the $i$-th variable
of $\varphi$.

Assume that Condition~\ref{condition for epsilon} (1) holds
and let $\omega\in \Omegam$.
We have
$
 \|I^{m,\rho}\|_{H^-}\le
C_7,
$
where $C_7$ depends on $\sigma, b, c, \varphi, C(B), \tilde{N}(B)$ 
polynomially.
\end{lemma}

\begin{proof}
We already proved Lemma~\ref{R estimate for tJmr} 
and Lemma~\ref{R estimate for tJmr^-1}.
Hence the proof is similar to that of
Lemma~\ref{estimate of discrete rough integral}.
\end{proof}

So far, we have given deterministic estimates of our processes
based on
$\tilde{C}(B)$ and $\tilde{N}(B)$.
We now give $L^p$ estimate of our processes.
The following result is due to
\cite{cll}.
See \cite{friz-hairer} also.

\begin{lemma}\label{CLL}
Assume that the covariance $R$ satisfies 
Condition~\upshape{\ref{condition on R}}.
Let $w$ be the control function defined in Definition~\upshape{\ref{w}}.
Then for any $\beta>0$, there exist positive numbers
$c_1$ and $c_2$ depending only on $H$ and $\beta$ such that
\begin{align}
 \mu\left(N_{\beta}(w)\ge r\right)\le c_1 e^{-c_2 r^{4H}}.
\label{exp decay}
\end{align}
\end{lemma}

The following is an immediate consequence of
Lemma~\ref{number N} and Lemma~\ref{CLL}.
Note that $N_{\beta}(\tilde{w})$ is a random variable defined on $\Omega_0$.

\begin{corollary}\label{cor to CLL}
 Assume the same assumption in Lemma~\upshape{\ref{CLL}}.
A similar estimate to $(\ref{exp decay})$ holds for
$N_{\beta}(\tilde{w})$.
\end{corollary}

By these results, under additional assumption on the covariance
of $(B_t)$, we obtain $L^p$ estimate of several quantities.

\begin{lemma}\label{L^p estimate for tJmr}
	Assume that Condition~\upshape{\ref{condition for epsilon}} (1) holds.
Let $\tilde{N}(B), C_5, C_6$ and $C_7$ be the positive numbers
defined in Lemmas~\upshape{\ref{estimate for tJmr}},
\upshape{\ref{R estimate for tJmr}},
\upshape{\ref{R estimate for tJmr^-1}} and
\upshape{\ref{estimate of discrete rough integral 2}}.
Then we have
\begin{align*}
 \max\big\{\tilde{N}(B),C_5, C_6, C_7\big\}\in \cap_{p\ge 1}L^p(\Omega_0).
\end{align*}
In particular 
\begin{align*}
  \sup_m
    \Big\|
      \max_{0\leq \rho\leq 1, t\in \Dm}
        \big\{
          |\tJmr_t(\xi,B)|,
          |\tJmr_{t}(\xi,B)^{-1}|
        \big\}1_{\Omegam}
      \Big\|_{L^p}
  <\infty.
\end{align*}
\end{lemma}

Consequently we obtain the following estimate.
Note that $\tZmr_t$ is a discrete process defined by \eqref{def tzmr}.
Also recall that
the notion of $\{a_m\}$-order nice discrete process
was introduced and the definition
of $\sup_{t,\rho}|\Ymr_t-Y_t|=O(a_m)$
was given in
Definition~\ref{nice discrete process}.

 \begin{theorem}\label{rough estimate of tzmr}
 Assume that Conditions~\upshape{\ref{condition on dm}} 
and~\upshape{\ref{condition for epsilon}} (1) hold.
 Let $\vep_1$ be the constant given in Condition~\ref{condition on dm}.
Set $a_m=\max\{\Delta_m^{3H^--1},\Delta_m^{\vep_1}\}$.
Then we have the following.
\begin{enumerate}
  \item It holds that 
  $\{\tZmr\}_{m}$ is an $\{a_m\}$-order nice discrete process 
with the H\"older exponent $\lambda=\min\{\lambda_1,2H^-\}$ which
is independent of $\rho$.
  \item It holds that
        $
          \sup_{t,\rho}|\Ymr_t-Y_t|=O(a_m)
        $
        in the sense of Definition~\upshape{\ref{nice discrete process}}~\upshape{(2)}.
\item For any $p\ge 1$ and $\kappa>0$, we have
\begin{align*}
 \lim_{m\to\infty}\|(2^m)^{\min\{3H^--1,\vep_1\}-\kappa}
\max_{t\in D_m}|\Ymh_t-Y_t|\|_{L^p}=0.
\end{align*}
\end{enumerate}
\end{theorem}

\begin{proof}
(1)~Note that the processes $(\tJmr)^{-1}$ and $c(\Ymr)$ appeared in \eqref{def tzmr} 
admit the uniform H\"older estimates
and that $d^m$ and $\hepm-\epm$ are 
$\{a_m\}$-order nice discrete processes (see Remark~\ref{rem89034111343}).
Hence the assertion follows from Remark~\ref{discrete young integral}.
(2) follows from (1) and Proposition~\ref{expression of ymh-ym}.
We prove (3).
By (2), there exists $X\in \cap_{p\ge 1}L^p(\Omega)$ such that
$\max_t|\Ymh_t-Y_t|\le a_mX$ on $\Omegam$.
Also we have for any $R>0$, there exists $C_R>0$ such that
$\mu\Big((\Omegam)^\complement\Big)\le C_R2^{-mR}$.
Using these estimates and the Schwarz inequality, we have
\begin{align*}
&\|(2^m)^{\min\{3H^--1,\vep_1\}-\kappa}\max_t|\Ymh_t-Y_t|\|_{L^p}^p\\
&\le
 E\left[(2^m)^{-\kappa p}
X^p ; \Omegam\right]+
E\left[(2^m)^{(\min\{3H^--1,\vep_1\}-\kappa)p}
\max_t|\Ymh_t-Y_t|^p ;(\Omegam)^\complement\right]\\
&\le 2^{-mp\kappa}\|X\|_{L^p}^p+
(2^m)^{(\min\{3H^--1,\vep_1\}-\kappa)p-R/2}C_R^{\frac{1}{2}}E[\max_t|\Ymh_t-Y_t|^{2p}]^{\frac{1}{2}}.
\end{align*}
Combining this estimate and Lemma~\ref{Ist general}, 
we complete the proof.
\end{proof}

We remark some consequences of the above results
in the case of the Milstein approximate solution.

\begin{remark}\label{estimate for J and J^-1}
\begin{enumerate}
 \item Let us consider non-random case.
That is, we consider a $\theta$-H\"older geometric rough path $(X,\XX)$.
The Milstein approximation solution $\Ymh_t$ ($t\in D_m$) is defined by
the similar equation to that explained in Section~\ref{2.2}
replacing $(B,\BB)$ by $(X,\XX)$.
Let $C(X)=\max\{\|X\|_{\theta},\|\XX\|_{2\theta}\}$.
Also we define $\tilde{N}(X)$ similarly to $\tilde{N}(B)$.
Note that $d^m\equiv 0$ and $\hepm\equiv 0$ and
we have the estimate $|\epm_{\tmkm,\tmk}|\le C\Delta_m^{3\theta}$,
where $C$ depends on $\sigma, b, C(X)$ polynomially.
Let $\kappa$ be a small positive number and set
$\theta^-=\theta-\kappa$.
We can view $(X,\XX)$ as a $\theta^-$-H\"older rough path.
Then for sufficiently large $m$, we have
\begin{align*}
\sup_{|t-s|\le 2^{-m}}
				\left|\frac{X_{s,t}}{(t-s)^{\theta^{-}}}\right|
+			  \sup_{|t-s|\le 2^{-m}}
			\left|\frac{\XX_{s,t}}{(t-s)^{2\theta^{-}}}\right|
			  \le
2^{-m\kappa+1}C(X)\le \frac{1}{2}.
\end{align*}
We can define an interpolated process
$\Ymr_t$ and $\tJmr_t$ similarly.
By the same argument as in this section,
we obtain,
\begin{align}
 \max_{t\in D_m}|\Ymh_{t}-Y_t|\le C\Delta_m^{3\theta^--1}, \label{Milstein estimate}
\end{align}
where $C$ depends on $\sigma, b$ and $C(X), \tilde{N}(X)$ polynomially.
Similar estimate was obtained by Davie~\cite{davie} and Friz-Victoir~\cite{friz-victoir}.
As for implementable versions, one 
can find some information in \cite{liu-tindel}.
We think our estimate makes clear how $C$ depends on $(X,\XX)$
more explicitly in (\ref{Milstein estimate}).
In Theorem~\ref{rough estimate of tzmr}, we deal with an RDE driven by
random rough path $(B,\BB)$ for which 
$\tilde{N}(B), C(B)\in \cap_{p\ge 1}L^p(\Omega_0)$ holds.
Hence, we can obtain $L^p$ convergence in (3).

\item 
We consider RDEs driven by $B$ which satisfies Condition~\ref{condition on R}.
We can prove $\sup_t\{|J_t|+|J^{-1}_t|\}\in \cap_{p\ge 1}L^p(\Omega)$
by applying the above results in the case where $\rho=1$
to the Milstein approximation
solution $(\Ymh_t, \tilde{J}^{m,1}_t)$ $(t\in D_m)$.
Note that $\Omegam=\Omega_0^{(m)}$ and 
$\liminf_{m\to\infty}\Omega_0^{(m)}=\Omega_0$ hold.
By Theorem~\ref{rough estimate of tzmr},
we see that $\lim_{m\to\infty}\max_{t\in D_m}|\Ymh_t-Y_t|=0$ for all
$\omega\in \Omega_0$.
Let $\hat{J}^{m,1}_t$ and $(\hat{J}^{m,1}_t)^{-1}$ $(t\in [0,1])$
be piecewise linear extensions of
$\tilde{J}^{m,1}_t$ and $(\tilde{J}^{m,1}_t)^{-1}$ 
$(t\in D_m)$ respectively.
Since $\tilde{J}^{m,1}_t$ and $(\tilde{J}^{m,1}_t)^{-1}$ 
are uniform H\"older continuous paths 
on $D_m$ which follow from
Lemmas~\ref{Ist general}, \ref{R estimate for tJmr},
\ref{R estimate for tJmr^-1},
so are $\hat{J}^{m,1}_t$ and $(\hat{J}^{m,1}_t)^{-1}$ on $[0,1]$.
This implies that for any subsequences of 
$\hat{J}^{m,1}_t$ and $(\hat{J}^{m,1}_t)^{-1}$,
there exist subsequences of them which converge uniformly
on $[0,1]$.
By the estimate in Lemmas~\ref{R estimate for tJmr},
\ref{R estimate for tJmr^-1} and the uniqueness of RDEs, any limits
of $\hat{J}^{m,1}_t$ and $(\hat{J}^{m,1}_t)^{-1}$ 
are equal to $J_t$ and $J^{-1}_t$ respectively.
This implies that the limits of themselves without taking subsequences exist and 
the limits $J_t$ and $J_t^{-1}$ also satisfy 
the same estimates as in (\ref{estimate of remainder term of tJmr})
and (\ref{estimate of remainder term of tJmr^-1})
for all $\omega\in\Omega_0$.
\item We can improve the estimate in Theorem~\ref{rough estimate of tzmr} (3)
when the driving process is an fBm as you can see
in Theorem~\ref{main theorem fbm} and Remark~\ref{rem849320482091}.
\end{enumerate}
\end{remark}

\subsection{Estimates of $J_t-\tJm_t$ and
$J_t^{-1}-(\tJm_t)^{-1}$ on $\Omega_0^{(m)}$}\label{4.3}
Throughout this section, $Y_t$ and $J_{t}$ denote
the solutions to \eqref{rde} and \eqref{Jacobian}, respectively.
Recall $\tJm_t=\tilde{J}^{m,0}_t$ is defined by \eqref{Mmr}.
Note that the recurrence relation for $\tJm$ does not contain
the terms $d^m$ and $\hepm$.
Hence we do not need assumptions on $d^m$ and $\hepm$
in this section.
Again, we assume $m$ satisfies (\ref{assumption on m}).
From now on, we will give estimates of
$J_t-\tJm_t$ and
$J_t^{-1}-(\tJm_t)^{-1}$.
We define
$\ep(J)_{\tmkm,\tmk}$ by
\begin{align}
J_{\tmk}
&=
	J_{\tmkm}
	+(D\sigma)(Y_{\tmkm})[J_{\tmkm}]B_{\tmkm,\tmk}
	+(D^2\sigma)(Y_{\tmkm})
		\left[J_{\tmkm},\sigma(Y_{\tmkm})e_{\alpha}\right]e_{\beta}
		\BB^{\alpha,\beta}_{\tmkm,\tmk}
	\nonumber\\
&\quad\qquad
	+
	(D\sigma)(Y_{\tmkm})\left[(D\sigma)(Y_{\tmkm})[J_{\tmkm}]e_{\alpha}\right]e_{\beta}
	\BB^{\alpha,\beta}_{\tmkm,\tmk}
	+(Db)(Y_{\tmkm})[J_{\tmkm}]\Delta_m
	\nonumber\\
&\quad\qquad
	+\ep(J)_{\tmkm,\tmk}.\label{eq of J}
\end{align}

\begin{lemma}\label{tJ-J}
  Let $\omega\in \Omega_0^{(m)}$.
  Let  
  \begin{align*}
    \dmJ_t=-\sum_{i=1}^{2^mt}
      \big(\tJm_{\tmi}\big)^{-1}
      \ep(J)_{\tmim,\tmi},
	\qquad
	t\in\Dm.
  \end{align*}

\begin{enumerate}
  \item It holds that
        \[
          |\ep(J)_{\tmkm,\tmk}|\le C_5\Delta_m^{3H^-},\quad 1\le k\le 2^m,
        \]
        where $C_5$ is the constant in Lemma~\upshape{\ref{R estimate for tJmr}}.
  \item $\{\delta^m(J)_t\}_{t\in\Dm}$ is a 
$\{\Delta_m^{3H^--1}\}$-order nice discrete process 
with the H\"older exponent
$2H^-$ and
        \begin{align*}
			\max_{t\in\Dm}|\delta^m(J)_t|=O(\Delta_m^{3H^--1}).
        \end{align*}
  \item For any natural number $R$, it holds that
        \begin{align}
          \tJm_t
          &=
            J_t
            \Bigg(
              I
              +
              \sum_{r=1}^R
                (\dmJ_t)^r
            \Bigg)
            +
            (\tilde{J}^m_t-J_t)\dmJ_t^R.\label{tJ-J expansion}
        \end{align}
   In particular,
   \begin{align}
    \max_{t\in\Dm}\left|\tilde{J}^m_t-J_t\left(I+\sum_{r=1}^R
  (\dmJ_t)^r\right)\right|=O(\Delta_m^{(3H^--1)(R+1)}).
\label{estimate of tJ-J}
   \end{align}
    \item For any natural numbers $L$ and $R$, it holds that
          \begin{multline*}
            \max_{t\in\Dm}
              \Bigg|
                (\tilde{J}^m_t)^{-1}
                -
	       	      \Bigg\{
                  I
                  +
                  \sum_{l=1}^L
                    \left(
                      -
                      \sum_{r=1}^R
                        (\dmJ_t)^r
                    \right)^l
	              \Bigg\}
				  J_t^{-1}
              \Bigg|\\
              =
                O(\Delta_m^{(3H^--1)(L+1)})+O(\Delta_m^{(3H^--1)(R+1)}).
	      \end{multline*}
   \end{enumerate}
  \end{lemma}

 \begin{proof}
\noindent
(1)~This follows from Lemma~\ref{R estimate for tJmr}
and Remark~\ref{estimate for J and J^-1}.

\noindent
(2) Similarly to $\epm_t$ and $\hepm_t$ (see \eqref{eq89034802914832}),
we set $\ep(J)^m_t=\sum_{i=1}^{2^mt} \ep(J)_{\tmim,\tmi}$ $(t\in \Dm)$.
From assertion~(1), $\ep(J)^m$ is a $\{\Delta_m^{3H^--1}\}$-order 
nice discrete process.
Hence, using the estimate of $\tJm$ and Remark~\ref{discrete young integral},
we see assertion~(2).

\noindent
  (3) 
	From the definition of $\tJm$ and \eqref{eq of J}, we have
\begin{align*}
J_t
&=
\tJm_t+\tJm_t\sum_{i=1}^{2^mt} \left(\tJm_{\tmi}\right)^{-1}
\ep(J)_{\tmim,\tmi}
=
\tJm_t-\tJm_t\dmJ_t
\end{align*}
		Hence
		$
			\tJm_t-J_t
			=
				J_t
				\dmJ_t
				+
				(\tJm_t-J_t)
				\dmJ_t
		$,
		which implies \eqref{tJ-J expansion}.
Noting $\tJm_t-J_t=\tJm_t\delta^m(J)_t$,
we get (\ref{estimate of tJ-J}).
	
\noindent
 (4)
Note that
\begin{align*}
  J_{t}^{-1}-(\tJm_{t})^{-1}
  &=
	-
      (\tJm_{t})^{-1}
      \big(J_{t}-\tJm_{t}\big)
      J_{t}^{-1} \\
  &=
	-
    J_t^{-1}
    \big(J_t-\tJm_t\big)
    J_t^{-1}
	+
	\big(J_{t}^{-1}-(\tJm_{t})^{-1}\big)
    \big(J_t-\tJm_t\big)
    J_t^{-1}.
 \end{align*}
Iterating this $L$ times and using the first identity above, we get
\begin{align*}
  J_t^{-1}-(\tJm_t)^{-1}
    &=
      -
      J_t^{-1}
      \sum_{l=1}^L
        \big[
          \big(J_t-\tJm_t\big)
          J_t^{-1}
        \big]^l
      +
	  (J_t^{-1}-(\tJm_t)^{-1})
      \big[
          \big(J_t-\tJm_t\big)
          J_t^{-1}
        \big]^L\\
    &=
      -
      J_t^{-1}
      \sum_{l=1}^L
        \big[
          \big(J_t-\tJm_t\big)
          J_t^{-1}
        \big]^l
	-
	(\tJm_{t})^{-1}
      \big[
          \big(J_t-\tJm_t\big)
          J_t^{-1}
        \big]^{L+1}.
\end{align*}
and
$
(\tJm_{t})^{-1}
  \big[
    \big(J_t-\tJm_t\big)
    J_t^{-1}
  \big]^{L+1}
  =
  O(\Delta_m^{(3H^--1)(L+1)})
$.
Thus
\begin{multline*}
  J_t^{-1}-(\tJm_t)^{-1}\\
  \begin{aligned}
	&=
		-
		J_t^{-1}
		\sum_{l=1}^L
		\left[
			\left(
			-
			J_t
			\sum_{r=1}^R (\dmJ_t)^r
			+
			O(\Delta_m^{(3H^--1)(R+1)})
			\right)
			J_t^{-1}
		\right]^l 
		+
		O(\Delta_m^{(3H^--1)(L+1)}) \\
	&=
		-
		J_t^{-1}
		\sum_{l=1}^L
		\left[
			\left(-
			J_t\sum_{r=1}^R (\dmJ_t)^r\right)
			J_t^{-1}
		\right]^l
		+O(\Delta_m^{(3H^--1)(L+1)})
		+LO(\Delta_m^{(3H^--1)(R+1)}).
  \end{aligned}
\end{multline*}
Since we have
\begin{align*}
	\left[
		\left(-
		J_t
		\sum_{r=1}^R (\dmJ_t)^r\right)
		J_t^{-1}
	\right]^l
	=
		J_t
		\left(
			-
			\sum_{r=1}^R (\dmJ_t)^r
		\right)^l
		J_t^{-1},
\end{align*}
we arrive at the conclusion.
 \end{proof}

\begin{remark}\label{rem_489011}
Summarizing above, we have the following.
By taking $L=R$ as a positive integer, we have
\begin{align*}
 \tJm_t-J_t&=J_tK^{1,m,R}_t+L^{1,m,R}_t,
 &
(\tJm_t)^{-1}-J_t^{-1}&=K^{2,m,R}_tJ_t^{-1}+L^{2,m,R}_t,
\end{align*}
where $K^{1,m,R}$ and $K^{2,m,R}$ are $\{\Delta_m^{3H^--1}\}$-order nice
discrete processes with the H\"older exponent $2H^-$
and $\max_t\{|L^{1,m,R}_t|+|L^{2,m,R}_t|\}=O(\Delta_m^{(3H^--1)R})$.
\end{remark}

\subsection{Convergence of $\tJmr_t$ and $(\tJmr_t)^{-1}$}\label{4.4}
Here we show convergence of $\tJmr_t$ and $(\tJmr_t)^{-1}$.
To this end we study
$
	N^{m,\rho}_t=(-\tJmr_t)^{-1}\partial_{\rho}\tJmr_t
$.
Note that $N^{m,\rho}_t$ is defined on $\Omegam$ and for large $m$
because $(-\tJmr_t)^{-1}$ can exist under the same condition.
\begin{lemma}\label{Nmr}
Assume that Conditions~\upshape{\ref{condition on dm}} 
and~\upshape{\ref{condition for epsilon}} (1) hold.
	Let $\vep_1$ be the constant given in Condition~\ref{condition on dm}.
   Set $a_m=\max\{\Delta_m^{3H^--1},\Delta_m^{\vep_1}\}$.
  Let $f_1,\dots,f_n$ be the standard basis of $\RR^n$
  and write $\tilde{Z}^{m,\rho,\nu}_t=(\tZmr_t,f_\nu)$ for $\nu=1,\dots,n$.
  Note that $\tilde{Z}^{m,\rho,\nu}$ is a real-valued process.
  \begin{enumerate}
   \item		Let $\omega\in \Omegam$. We have 
	\begin{align*}
		N^{m,\rho}_t
		&=
			\sum_{\nu=1}^n
				\tilde{Z}^{m,\rho,\nu}_t
				I^{m,\rho}(\varphi_\nu)_t
			+
			\sum_{\lambda=0}^3
				I_\lambda(N^{m,\rho})_t
	\end{align*}
	Here, $\varphi_\nu(x,M_1,M_2)$ 
	$\left(x\in \RR^n, M_1, M_2\in \mathcal{L}(\RR^n)\right)$
	is an
	$\mathcal{L}\left(\RR^d, \mathcal{L}(\RR^n)\right)$-valued function
	defined by
	\begin{align*}
		\varphi_\nu(x,M_1,M_2)
		=
			-M_2(D^2\sigma)(x)[M_1f_\nu,M_1]
	\end{align*}
	and $I^{m,\rho}(\varphi_\nu)$ is a discrete rough integral
  defined in Lemma~\upshape{\ref{estimate of discrete rough integral 2}}.
  Explicitly,
  we have, for $t\in\Dm$,
	\begin{multline*}
		\varphi_\nu
			\left(
				\Ymr,\tJmr,(\tJmr)^{-1}
			\right)^{\bcdot}_t[v\otimes w]\\
		\begin{aligned}
			&=
				-
				(\tJmr_t)^{-1}
				(D\sigma)(\Ymr_t)
					\left[
						(D^2\sigma)(\Ymr_t)\left[\tJmr_tf_\nu,\tJmr_t\right]w
					\right]
				v\\
			&\phantom{=}\qquad
				+
				(\tJmr_t)^{-1}
				(D^3\sigma)(\Ymr_t)
					\left[
						\sigma(\Ymr_t)v,
						\tJmr_tf_\nu,
						\tJmr_t
					\right]
				w\\
			&\phantom{=}\qquad
				+
				(\tJmr_t)^{-1}
				(D^2\sigma)(\Ymr_t)
					\left[
						(D\sigma)(\Ymr_t)\left[\tJmr_tf_\nu\right]v,\tJmr_t
					\right]
				w\\
			&\phantom{=}\qquad
				+
				(\tJmr_t)^{-1}
				(D^2\sigma)(\Ymr_t)
					\left[
						\tJmr_tf_\nu,
						(D\sigma)(\Ymr_t)\left[\tJmr_t\right]v
					\right]
				w,
			\qquad
			v,w\in\RR^d.
		\end{aligned}
  \end{multline*}
  Also
\begin{align*}
	I_0(N^{m,\rho})_t
	&=
		-	
		\sum_{\nu=1}^n
		\sum_{j=1}^{2^m t}
			\tilde{Z}^{m,\rho,\nu}_{\tmjm,\tmj}
			I^{m,\rho}(\varphi_\nu)_{\tmj},\\
	I_1(N^{m,\rho})_t
	&=
		\sum_{j=1}^{2^m t}
			(-\tJmr_{\tmjm})^{-1}
			(D^2b)(\Ymr_{\tmjm})[\Zmr_{\tmjm},\tJmr_{\tmjm}]
			\Delta_m,\\
	I_2(N^{m,\rho})_t
	&=
		\sum_{j=1}^{2^m t}
			(-\tJmr_{\tmjm})^{-1}
			\Big\{(Dc)(\Ymr_{\tmjm})[\tJmr_{\tmjm}]\\
	&\qquad\qquad\qquad\qquad\qquad
			+
			\rho(D^2c)(\Ymr_{\tmjm})[\Zmr_{\tmjm},\tJmr_{\tmjm}]
			\Big\}
			d^m_{\tmjm,\tmj},
\end{align*}
	and $I_3(N^{m,\rho})$ is the residual term defined
by 
\begin{align*}
	I_3(N^{m,\rho})_t=N^{m,\rho}_t-
\sum_{\nu=1}^n
				\tilde{Z}^{m,\rho,\nu}_t
				I^{m,\rho}(\varphi_\nu)_t
-\sum_{\lambda=0}^2I_{\lambda}(N^{m,\rho})_t.
\end{align*}
  \item	$I_0(N^{m,\rho})$, $I_1(N^{m,\rho})$, $I_2(N^{m,\rho})$ and
$I_3(N^{m,\rho})$ are 
  $\{a_m\}$-order nice discrete processes 
	with the H\"older exponent $\lambda=\min\{\lambda_1,2H^-\}$.
  In addition, $\sup_{\rho}\|N^{m,\rho}\|_{H^-}=O(a_m)$
  in the sense of Definition~\upshape{\ref{nice discrete process}}~\upshape{(2)}.
  \end{enumerate}
  \end{lemma}

\begin{proof}
	From \eqref{expansion of Mmr2}, we have
	\begin{align*}
		\Nmr_{\tmj}
		=
			\Nmr_{\tmjm}
			+
			(-\tJmr_{\tmj})^{-1}
			\left\{\partial_{\rho}E^{m,\rho}(\Ymr_{\tmjm},\theta^m_{\tmjm}B)\right\}
			\tJmr_{\tmjm}.
	\end{align*}
	Using
	$
		(\tJmr_{\tmj})^{-1}
		=
			(\tJmr_{\tmjm})^{-1}
			\{
				I
				-(D\sigma)(\Ymr_{\tmjm})B_{\tmjm,\tmj}
				+O(\Delta_m^{2H^-})
			\}
	$
	due to Lemma~\ref{expansion of Mmr}
	and the expression of
	$
		\partial_\rho E^{m,\rho}(\Ymr_{\tmjm},\theta^m_{\tmjm}B)
	$,
	we have
	\begin{align}
		\Nmr_{\tmj}-\Nmr_{\tmjm}
		&=
			(-\tJmr_{\tmjm})^{-1}
			\Big\{
			(D^2\sigma)(\Ymr_{\tmjm})[\Zmr_{\tmjm},\tJmr_{\tmjm}]B_{\tmjm,\tmj} \nonumber\\
		&\qquad\qquad\qquad\qquad
				-
				(D\sigma)(\Ymr_{\tmjm})
					\left[
						(D^2\sigma)(\Ymr_{\tmjm})[\Zmr_{\tmjm},\tJmr_{\tmjm}]B_{\tmjm,\tmj}
					\right]
					B_{\tmjm,\tmj} \nonumber \\
		&\qquad\qquad\qquad\qquad
				+
				D^2((D\sigma)[\sigma])(\Ymr_{\tmjm})[\Zmr_{\tmjm},\tJmr_{\tmjm}]\BB_{\tmjm,\tmj}
			\Big\} \nonumber \\
		&\phantom{=}\quad
			+
			(-\tJmr_{\tmjm})^{-1}
			\Big[
				(Dc)(\Ymr_{\tmjm})[\tJmr_{\tmjm}]
				+
				\rho(D^2c)(\Ymr_{\tmjm})[\Zmr_{\tmjm},\tJmr_{\tmjm}]
			\Big]d^m_{\tmjm,\tmj} \nonumber \\
		&\phantom{=}\quad
				+
				(-\tJmr_{\tmjm})^{-1}
				\Big[
				(D^2b)(\Ymr_{\tmjm})[\Zmr_{\tmjm},\tJmr_{\tmjm}]\Delta_m
			\Big]
			+
			O(\Delta_m^{3H^-}).
		\label{eq4349023211}
	\end{align}
	Next we take the sum over $0\leq j\leq 2^mt$.
	Applying $B^\alpha_{s,t}B^\beta_{s,t}-B^{\alpha,\beta}_{s,t}=B^{\beta,\alpha}_{s,t}$
	and substituting $\Zmr_{\tmjm}=\tJmr_{\tmjm}\tZmr_{\tmjm}=\sum_{\nu=1}^n\tilde{Z}^{m,\rho,\nu}_{\tmjm} \tJmr_{\tmjm}f_\nu$,
	we see that the summation of the first term in \eqref{eq4349023211} gives
	\begin{align*}
		\sum_{\nu=1}^n
		\sum_{j=1}^{2^mt}
			\tilde{Z}^{m,\rho,\nu}_{\tmjm}
			I^{m,\rho}(\varphi_\nu)_{\tmj,\tmjm}
		=
			\sum_{\nu=1}^n
				\tilde{Z}^{m,\rho,\nu}_t
				I^{m,\rho}(\varphi_\nu)_t
			-
			\sum_{\nu=1}^n
			\sum_{j=1}^{2^mt}
				\tilde{Z}^{m,\rho,\nu}_{\tmj,\tmjm}
				I^{m,\rho}(\varphi_\nu)_{\tmj}.
	\end{align*}
	The summations of the second and third terms in \eqref{eq4349023211}
	give $I_2(N^{m,\rho})$ and $I_1(N^{m,\rho})$, respectively.
	The summation of the fourth term $O(\Delta_m^{3H^-})$
in \eqref{eq4349023211} is $I_3(N^{m,\rho})$,
	which is an $\{a_m\}$-order nice discrete process.
	This completes the proof of (1).

	We show assertion~(2).
	Recall that the discrete H\"older norm 
$\|I^{m,\rho}(\varphi_\nu)\|_{H^-}$ can be estimated by a constant which
depends on $\sigma, b, c$, $C(B)$ and $\tilde{N}(B)$ polynomially
	(see Lemma~\upshape{\ref{estimate of discrete rough integral 2}})
	and that 
	$\tilde{Z}^{m,\rho,\nu}$ is
	an $\{a_m\}$-order nice discrete process
	(see Theorem~\ref{rough estimate of tzmr}).
	Thus, the discrete version of the estimate of Young integrals
(Remark~\ref{discrete young integral}) implies
	that $I_0(N^{m,\rho})$ is an $\{a_m\}$-order nice discrete process.
	Noting that we have good estimates of
$H^-$-H\"older norm of $\Ymr$, $\tJmr$, $(-\tJmr)^{-1}$ (Lemma~\ref{Ist general}, 
Lemma~\ref{R estimate for tJmr}, Lemma~\ref{R estimate for tJmr^-1}) and
that $\Zmr$ is an $\{a_m\}$-order nice discrete process
(Theorem~\ref{rough estimate of tzmr}),
	we see that 
$I_1(N^{m,\rho})$ is an $\{a_m\}$-order nice discrete process.
Since $d^m$ is an $\{a_m\}$-order nice discrete process, 
$I_2(N^{m,\rho})$ is as well.
As for $I_3^{m,\rho}$, we already proved the assertion.
	Here we used 
	Lemmas~\ref{R estimate for tJmr}, \ref{R estimate for tJmr^-1}, 
	and \ref{estimate of discrete rough integral 2}
	and Theorem~\ref{rough estimate of tzmr}.
	Since $\sup_{\rho}\|\tZmr\|_{H^-}=O(a_m)$
	and other terms are $\{a_m\}$-order nice discrete processes,
	we have $\sup_{\rho}\|N^{m,\rho}\|_{H^-}=O(a_m)$ 
        which completes the proof of
	assertion~(2).
\end{proof}

  \begin{theorem}\label{convergence of J and J-1}
	Assume that Conditions~\upshape{\ref{condition on dm}} 
and~\upshape{\ref{condition for epsilon}} (1) hold.
	Let $\vep_1$ be the constant given in Condition~\ref{condition on dm}.
   Set $a_m=\max\{\Delta_m^{3H^--1},\Delta_m^{\vep_1}\}$.
    Then we have
    \begin{align*}
      \sup_{t,\rho}|\tJmr_t-J_t|
      &=
        O(a_m),
      &
      \sup_{t,\rho}|(\tJmr_t)^{-1}-J_t^{-1}|
      &=
        O(a_m)
    \end{align*}
    in the sense of Definition~\upshape{\ref{nice discrete process}}~\upshape{(2)}.
  \end{theorem}
  \begin{proof}
    Note that
    \begin{align*}
      \tJmr_t-\tJm_t
      &=
        \int_0^\rho
          \partial_{\rho_1}J^{m,\rho_1}_t
          d\rho_1
      =
        \int_0^\rho
          (-J^{m,\rho_1}_t)N^{m,\rho_1}_t
          d\rho_1,\\
      (\tJmr_t)^{-1}-(\tJm_t)^{-1}
      &=
        \int_0^{\rho}
          \partial_{\rho_1}
          (J^{m,\rho_1}_t)^{-1}
          d\rho_1
      =
        \int_0^{\rho}
          N^{m,\rho_1}_t
          (J^{m,\rho_1}_t)^{-1}
          d\rho_1.
    \end{align*}
    From Lemmas~\ref{estimate for tJmr} and \ref{Nmr}, we see that
    $
      \sup_{t,\rho}|\tJmr_t-\tJm_t|=O(a_m)
    $
    and
    $
      \sup_{t,\rho}|(\tJmr_t)^{-1}-(\tJm_t)^{-1}|=O(a_m)
    $.
    This and Remark~\ref{rem_489011} yield the assertion.
  \end{proof}

\section{Proof of main theorem}\label{proof of main theorem}

We prove Theorem~\ref{main theorem} and
Corollary~\ref{corollary for main theorem} 
in Section~\ref{subsection proof of main theorem}.
Section~\ref{lemmas} is a preparation for it.

\subsection{Lemmas}\label{lemmas}

Throughout this section, we assume that 
Conditions~\ref{condition on dm}\,$\sim$\,\ref{condition on K} hold.
Recall that $(\lambda_1,\vep_1, G_1)$ and
$(\lambda_2,\vep_2,G_2)$ are the triples of the two constants
and the random variable specified
in Conditions~\ref{condition on dm} and~\ref{condition on K}, respectively.
Also, set
\begin{align*}
	a_m
	=
		\max\{\Delta_m^{3H^--1},\Delta_m^{4H^- - 2H - \frac{1}{2}},\Delta_m^{\vep_1},\Delta_m^{\vep_2}\}.
\end{align*}
We will give estimate of $\tZmr(\omega)$ for $\omega\in \Omegam$.
Precisely, we prove 

\begin{lemma}\label{estimate of tzmr}
There exists a positive integer $m_0$ such that
for all $p\ge 1$ it holds that
\begin{align*}
 \sup_{m\ge m_0}
\left\|\sup_{0\le \rho\le 1}
\|(2^m)^{2H-\frac{1}{2}}\tZmr\|_{H^-}1_{\Omega^{(m,d^m)}_0}\right\|_{L^p}
<\infty.
\end{align*}
\end{lemma}

We refer the readers to Definition~\ref{definition of tzmr} and \eqref{defIm}
for definition of $\tZmr_t$ and $I^m$.
We decompose as
$\tZmr_t-I^m_t=\sum_{i=1}^5 S^{m,\rho,i}_t$, where
\begin{align*}
	S^{m,\rho,1}_t
	&=
		\sum_{i=1}^{2^mt}
		(\tJmr_{\tmi})^{-1}
		\left(c(\Ymr_{\tmim})-c(Y_{\tmim})\right)
		d^m_{\tmim,\tmi},\\
	S^{m,\rho,2}_t
	&=
		\sum_{i=1}^{2^mt}
			\left((\tJmr_{\tmi})^{-1}-(\tJm_{\tmi})^{-1}\right)
			c(Y_{\tmim})d^m_{\tmim,\tmi},\\
	S^{m,\rho,3}_t
	&=
		\sum_{i=1}^{2^mt}
			\left((\tJm_{\tmi})^{-1}-J_{\tmi}^{-1}\right)
			c(Y_{\tmim})
			d^m_{\tmim,\tmi},
	\quad
	S^{m,\rho,4}_t
	=
		\sum_{i=1}^{2^mt}
			J^{-1}_{\tmim,\tmi}c(Y_{\tmim})d^{m}_{\tmim,\tmi},\\
	S^{m,\rho,5}_t
	&=
		\sum_{i=1}^{2^mt}
			(\tJmr_{\tmi})^{-1}
			\left(\hepm_{\tmim,\tmi}-\epm_{\tmim,\tmi}\right).
\end{align*}

We give estimates for each term $S^{m,\rho,i}$ ($1\le i\le 5$).
First, we consider $S^{m,\rho,1}$.

\begin{lemma}\label{lemma for Smrho1}
Let $\omega\in \Omegam$.
Then we have
 \begin{align*}
  \|(2^m)^{2H-\frac{1}{2}}S^{m,\rho,1}\|_{\lambda_1}\le
a_m C G_1 \sup_{\rho}
\|(2^m)^{2H-\frac{1}{2}}\tZmr\|_{H^-},
 \end{align*}
where $C$ depends only on $\tilde{C}(B)$ and $\tilde{N}(B)$
polynomially.
\end{lemma}

\begin{proof}
	Set $F^{m,\rho}_t=(\tJmr_{t+\Delta_m})^{-1}(c(\Ymr_t)-c(Y_t))$.
 We have
\begin{align*}
 c(\Ymr_t)-c(Y_t)&=
\int_0^{\rho}(Dc)(Y^{m,\rho_1}_t)[Z^{m,\rho_1}_t]d\rho_1
=
\int_0^{\rho}(Dc)(Y^{m,\rho_1}_t)[\tilde{J}^{m,\rho_1}_t
\tilde{Z}^{m,\rho_1}_t]d\rho_1
\end{align*}
and we obtain 
H\"older estimate of the discrete process
$
\| F^{m,\rho}\|_{H^-}
\le
C \sup_{\rho}\|\tZmr\|_{H^-}
$.
Here, $C$ depends on the H\"older norms of $\Ymr$ and $\tJmr$.
By combining the estimate $\|d^m\|_{\lambda_1}\le 2^{-m\vep_1}G_1\leq a_m G_1$
($\omega\in \Omega_0$) and Remark~\ref{discrete young integral},
we complete the proof.
\end{proof}

Next, we consider $S^{m,\rho,4}$ and $S^{m,\rho, 5}$.

\begin{lemma}\label{estimate of remainder term1}
Let $\omega\in \Omegam$.
 We have
\begin{align*}
 J^{-1}_{\tmim,\tmi}c(Y_{\tmim})d^m_{\tmim,\tmi}=
-J^{-1}_{\tmim}(D\sigma)(Y_{\tmim})[c(Y_{\tmim})
d^m_{\tmim,\tmi}]B_{\tmim,\tmi}
+
O(\Delta_m^{4H^-}),
\end{align*}
where the dominated random variable for the term
$O(\Delta_m^{4H^-})$
depends only on $\tilde{C}(B)$ and $\tilde{N}(B)$ polynomially.
\end{lemma}

\begin{proof}
 This follows from Lemma~\ref{R estimate for tJmr^-1} 
and Remark~\ref{estimate for J and J^-1}.
We used $\la_1>H^-$.
\end{proof}

\begin{lemma}\label{estimate of remainder term2}
Let $\omega\in \Omegam$.
There exist $\RR^n$-valued 
bounded Lipschitz functions $\varphi^{\alpha,\beta,\gamma}$,
$\psi_{\alpha}$, $F_{\alpha,\beta,\gamma}$, $F^1_{\alpha}$,
$F^2_{\alpha}$
on $\RR^n$
$(1\le \alpha,\beta,\gamma\le d)$ such that
\begin{multline*}
	(\tJmr_{\tmi})^{-1}\big(\hepm_{\tmim,\tmi}-\epm_{\tmim,\tmi}\big)\\
	\begin{aligned}
		&=
			(\tJmr_{\tmim})^{-1}
			\Big\{
				\sum_{\alpha,\beta,\gamma}
					\varphi_{\alpha,\beta,\gamma}(\hat{Y}_{\tmim})
					B^{\alpha,\beta,\gamma}_{\tmim,\tmi}
				+
				\sum_{\alpha}
					\psi_{\alpha}(\Ymh_{\tmim})
					B^{\alpha}_{\tmim,\tmi}\Delta_m\\
		&\quad\quad
				+
					\sum_{\alpha,\beta,\gamma}
						F_{\alpha,\beta,\gamma}(Y_{\tmim})
						B^{\alpha,\beta,\gamma}_{\tmim,\tmi}
				+\sum_{\alpha}F^1_{\alpha}(Y_{\tmim})B^{0,\alpha}_{\tmim,\tmi}
				+\sum_{\alpha}F^2_{\alpha}(Y_{\tmim})B^{\alpha,0}_{\tmim,\tmi}
			\Big\}\\
		&\quad\quad
			+
			O(\Delta_m^{4H^-}).
	\end{aligned}
\end{multline*}
The dominated random variables for the terms
$O(\Delta_m)^{4H^-}$
depends on $\tilde{C}(B)$ and $\tilde{N}(B)$ polynomially.
\end{lemma}

\begin{proof}
	From \eqref{estimate of Emr^-1}, Condition~\ref{condition for epsilon} (1)
	and Lemma~\ref{lem8345901809321} (1),
	we have
	\begin{align*}
		(\tJmr_{\tmi})^{-1}\big(\hepm_{\tmim,\tmi}-\epm_{\tmim,\tmi}\big)
		=
			(\tJmr_{\tmim})^{-1}\big(\hepm_{\tmim,\tmi}-\epm_{\tmim,\tmi}\big)
			+
			O(\Delta_m^{4H^-}).
	\end{align*}
	Combining this identity with Condition~\ref{condition for epsilon} (2) and 
Lemma~\ref{lem8345901809321} (2) yields the desired estimate.
\end{proof}

As we have shown in the above lemmas, we need estimates
for weighted sum process in Wiener chaos of order $3$
and sum process of $d^{m,\alpha,\beta}_{\tmim,\tmi}B^{\gamma}_{\tmim,\tmi}$.
We refer the readers to \eqref{eqK} for the definition of $\mathcal{K}^3_m$.

\begin{lemma}\label{lem48291143}
	Let $\omega\in \Omegam$. Let $K^{m}\in \mathcal{K}^3_{m}$
	and $\{F^m_t\}_{t\in \Dm}$ be a discrete process satisfying
	$|F^m_0|+\|F^m\|_{H^-}\le C$, where
	$C$ is independent of $m$ and
	depends only on $\tilde{C}(B)$ and $\tilde{N}(B)$ polynomially.
	Let
	$
	I^m(F^m)_t=\sum_{i=1}^{2^mt}F^m_{\tmim}K^{m}_{\tmim,\tmi}
	$
	$(t\in \Dm)$.
	Then it holds that 
	\begin{align*}
	\|(2^m)^{2H-\frac{1}{2}}I^m(F^m)\|_{\lambda_2}
	\le
		a_m
		C
		G_2,
\end{align*}
where $C$ depends on
$\tilde{C}(B)$ and $\tilde{N}(B)$ polynomially.
\end{lemma}

\begin{proof}
	By the assumption on the H\"older norm of
$F^m$ and Condition~\ref{condition on K} and 
using Remark~\ref{discrete young integral}, we have
$
  \|(2^m)^{2H-\frac{1}{2}}I^m(F^m)\|_{\lambda_2}
  \le
	\Delta_m^{\vep_2}
    C
    G_2,
$
which implies the assertion.
\end{proof}

\begin{lemma}\label{lemma for Smrho5}
	Let $\omega\in \Omegam$. 
We have
\begin{align*}
 \|(2^m)^{2H-\frac{1}{2}}S^{m,\rho,4}\|_{\lambda_2}
 +
 \|(2^m)^{2H-\frac{1}{2}}S^{m,\rho,5}\|_{\lambda_2}
 \le  
 a_m
  C
  \{G_2+1\},
\end{align*}
where $C$ depends on
$\tilde{C}(B)$ and $\tilde{N}(B)$ polynomially.
\end{lemma}

\begin{proof}
We use the decompositions in Lemmas~\ref{estimate of remainder term1} and~\ref{estimate of remainder term2}.
First, we consider the sum of $O(\Delta_m^{4H^-})$.
Let $s=\tmk<\tml=t$.
We have
\begin{align*}
	\left|
		(2^m)^{2H-\frac{1}{2}}
		\sum_{i=k}^{l-1}
			O(\Delta_m^{4H^-})
	\right|
\le (2^m)^{2H-\frac{1}{2}}(l-k) C\Delta_m^{4H^-}
=\Delta_m^{4H^--2H-\frac{1}{2}}C(t-s).
\end{align*}
where $C$ depends on $\tilde{C}(B)$ and $\tilde{N}(B)$ polynomially.
This term can be estimated as in the assertion.
As for sum process $K^m_{s,t}=\Delta_m B^{\alpha}_{s,t}$ which defined by the term
$\Delta_m B^{\alpha}_{\tmim,\tmi}$ in Lemma~\ref{estimate of remainder term2}, 
we have similar estimate to
the elements in $\mathcal{K}^3_m$.
See the proof of Lemma~\ref{K satisfies Condition}.
Note that we use Condition~\ref{condition on R} only in that proof.
The remaining main terms 
can be handled by Lemma~\ref{lem48291143} and Condition~\ref{condition on K}.
This completes the proof.
\end{proof}

\begin{remark}
 In the above Lemmas~\ref{lem48291143} and \ref{lemma for Smrho5},
we used the estimate of $K^m_{s,t}$ which is defined as the sum process
of $B^{0,\alpha}_{\tmim,\tmi}$ and $B^{\alpha,0}_{\tmim,\tmi}$
in Condition~\ref{condition on K}.
If we use the estimate $|B^{\alpha,0}_{\tmim,\tmi}|\le C\Delta_m^{1+H^-}$,
which follows form
the H\"older estimate of $B$ only, we obtain a rough estimate
$|(2^m)^{2H-\frac{1}{2}}K^m_{s,t}|\le C\Delta_m^{H^--(2H-\frac{1}{2})}|t-s|$ similarly
to the estimate of $(2^m)^{2H-\frac{1}{2}}\sum O(\Delta_m^{4H^-})$ in the proof of 
Lemma~\ref{lemma for Smrho5}.
However, this estimate will give the estimate
$\vep<\min\{3H^--1,H^--(2H-\frac{1}{2}),\vep_1,\vep_2\}$.
Clearly this estimate gets worse as $H\to \frac{1}{2}$.
\end{remark}

We consider the estimates of $S^{m,\rho,3}$.
To this end, recall definition \eqref{defIm} of $I^m$ and set 
\begin{align}\label{def X}
  X_m
  &=
    \|(2^m)^{2H-\frac{1}{2}}I^m|_{\Dm}\|_{H^-}.
\end{align}
Then from Condition~\ref{minimal moment estimate}, 
we have $\sup_m\|X_m\|_{L^p}<\infty$ for all $p\ge 1$.
\begin{lemma}\label{Smrho3}
  Let $\omega\in \Omegam$.
  We have
  \begin{align*}
    \|(2^m)^{2H-\frac{1}{2}}S^{m,\rho,3}\|_{H^-}
    &\le
		a_m
      C
	  \{
		X_m
      	+G_2
      	+1
	\},
  \end{align*}
  where $C$ depends on $\tilde{C}(B)$ and $\tilde{N}(B)$ polynomially.
  \end{lemma}
  
\begin{proof}
	Let $R$ be a positive integer.
	From Remark~\ref{rem_489011}, we have
	\begin{align}
		(\tJm_{\tmi})^{-1}-J_{\tmi}^{-1}
		&=
			K^{2,m,R}_{\tmi}J_{\tmi}^{-1}
			+L^{2,m,R}_{\tmi} \nonumber \\
		&=
			K^{2,m,R}_{\tmi}J^{-1}_{\tmim}
			+K^{2,m,R}_{\tmi}J^{-1}_{\tmim,\tmi}
			+L^{2,m,R}_{\tmi}.
			\label{eq48930184231}
	\end{align}
	where $K^{2,m,R}$ is an $\{a_m\}$-order nice discrete processes
	and $L^{2,m,R}$ is a small discrete process.
	Hence
  \begin{align*}
    S^{m,\rho,3}_t
    &=
      \sum_{i=1}^{2^m t}
        \left(
            K^{2,m,R}_{\tmi}J^{-1}_{\tmim}
            +K^{2,m,R}_{\tmi}J^{-1}_{\tmim,\tmi}
            +L^{2,m,R}_{\tmi}
        \right)
        c(Y_{\tmim})
        d^m_{\tmim,\tmi}\\
    &=
      S^{m,\rho,3,1}_t
      +S^{m,\rho,3,2}_t
      +S^{m,\rho,3,3}_t,
  \end{align*}
  Then with the help of 
	the summation by parts formula~\eqref{summation by parts formula},
  we have
  \begin{align*}
    S^{m,\rho,3,1}_t
    =
      \sum_{i=1}^{2^m t}
          K^{2,m,R}_{\tmi}
          I^m_{\tmim,\tmi}
    =
		K^{2,m,R}_t
		I^m_t
      -
      \sum_{i=1}^{2^m t}
		K^{2,m,R}_{\tmim,\tmi}
		I^m_{\tmim}.
  \end{align*}
  Recalling that $(2^m)^{2H-\frac{1}{2}}I^m|_{\Dm}$ is discrete $H^-$-H{\"o}lder continuous
  and using Remark \ref{rem_489011},
using $X_m$ defined by \eqref{def X}, we have
  \begin{multline*}
    \|(2^m)^{2H-\frac{1}{2}}S^{m,\rho,3,1}\|_{H^-}\\
	\begin{aligned}
		&\leq
		2
		\|K^{2,m,R}\|_{H^-}
		\|(2^m)^{2H-\frac{1}{2}}I^m|_{\Dm}\|_{H^-}
		+
		\left\|
			\sum_{i=1}^{2^m \cdot}
				K^{2,m,R}_{\tmim,\tmi}
				(2^m)^{2H-\frac{1}{2}}I^m_{\tmim}
		\right\|_{H^-}\\
	&\leq
		C
		\{
			a_m
			\cdot
			X_m
			+
			a_m
			\cdot
			X_m
		\}.
	\end{aligned}
  \end{multline*}
  In a similar way to Lemma \ref{lemma for Smrho5},
  using Lemma~\ref{estimate of remainder term1}, we have
  \begin{align*}
    \|(2^m)^{2H-\frac{1}{2}}
      S^{m,\rho,3,2}
    \|_{\lambda_2}
    \le  
	a_m C\{ G_2+1 \}.
   \end{align*}
   The term $\|(2^m)^{2H-\frac{1}{2}}S^{m,\rho,3,3}\|_{H^-}$ becomes small for large $R$.
   The proof is completed.
\end{proof}

Finally, we estimate $S^{m,\rho,2}$.
To this end, we use 
$
 N^{m,\rho}_t=(-\tJmr_t)^{-1}\partial_{\rho}\tJmr_t
$,
which is introduced in Section~\ref{4.4}.

\begin{lemma}\label{lem483901}
  Let $L$ be a positive integer.
 Then it holds that
  \begin{multline*}
  (\tJmr_t)^{-1}-(\tJm_t)^{-1}
  =
    \sum_{l=1}^{L-1}
      \int_{0<\rho_l<\cdots<\rho_1<\rho}
        d\rho_1
        \cdots
        d\rho_l\,
        N^{m,\rho_1}_t\cdots N^{m,\rho_l}_t
        (\tJm_t)^{-1}\\ 
    +
    \int_{0<\rho_{L}<\cdots<\rho_1<\rho}
      d\rho_1
      \cdots
      d\rho_L\,
      N^{m,\rho_1}_t
      \cdots 
      N^{m,\rho_L}_t
      (\tilde{J}^{m,\rho_L}_t)^{-1}.
 \end{multline*}
 \end{lemma}
 
 \begin{proof}
 Noting $\partial_{\rho}(\tJmr_t)^{-1}=
 -(\tJmr_t)^{-1}\partial_{\rho}\tJmr_t
 (\tJmr_t)^{-1}=N^{m,\rho}_t(\tJmr_t)^{-1}$,
 we have
 \begin{multline*}
  (\tJmr_t)^{-1}-(\tJm_t)^{-1}
  =
    \int_{0<\rho_1<\rho}
      d\rho_1\,
      N^{m,\rho_1}_t
      (\tilde{J}^{m,\rho_1}_t)^{-1}\\
  \begin{aligned}
    &=
      \int_{0<\rho_1<\rho}
        d\rho_1\,
        N^{m,\rho_1}_t
        (\tJm_t)^{-1}
      +
      \int_{0<\rho_1<\rho}
        d\rho_1
        N^{m,\rho_1}_t
        \left\{
			(\tilde{J}^{m,\rho_1}_t)^{-1}-(\tJm_t)^{-1}
        \right\}\\
  &=
    \int_{0<\rho_1<\rho}
      d\rho_1\,
      N^{m,\rho_1}_t
      (\tJm_t)^{-1}
    +
    \int_{0<\rho_1<\rho}
      d\rho_1
      N^{m,\rho_1}_t
      \int_{0<\rho_2<\rho_1}
        d\rho_2\,
        N^{m,\rho_2}_t
        (\tilde{J}^{m,\rho_2}_t)^{-1}.
  \end{aligned}
 \end{multline*}
 Iterating this calculation, we are done.
 \end{proof}

 \begin{lemma}\label{Smrho2}
 For $\omega\in\Omegam$, we have
  \begin{align*}
    \|(2^m)^{2H-\frac{1}{2}}S^{m,\rho,2}\|_{H^-}
    \le
		a_m
      C
      \big\{
        G_1
        \sup_{\rho}\|(2^m)^{2H-\frac{1}{2}}\tZmr\|_{H^-}
        +X_m
        +G_2
        +1
      \big\},
  \end{align*}
 where $C$ depends on $\tilde{C}(B), \tilde{N}(B)$ polynomially.
 \end{lemma}

\begin{proof}
  We use the same notation as in Lemmas~\ref{Nmr} and \ref{lem483901}.
  Set
  \begin{align*}
    \tilde{N}^{m,\rho_1,\ldots,\rho_l}_t
    &=
      \prod_{r=1}^l
        \left\{
          N^{m,\rho_r}_t
          -
          \sum_{\nu=1}^n
            \tilde{Z}^{m,\rho_r,\nu}_t
          I^{m,\rho_r}(\varphi_\nu)_t
        \right\}
	=
		\prod_{r=1}^l
		\sum_{\lambda=0}^3
			I_\lambda(N^{m,\rho_r})_t,\\
    R^{m,\rho_1,\ldots,\rho_l}_t
    &=
      N^{m,\rho_1}_t\cdots N^{m,\rho_l}_t
      -
      \tilde{N}^{m,\rho_1,\ldots,\rho_l}_t.
  \end{align*}
  Note that the product $\prod_{r=1}^l$
  in the above equation should be taken according to the order.
  Then we have $S^{m,\rho,2}_t=S^{m,\rho,2,1}_t+S^{m,\rho,2,2}_t+S^{m,\rho,2,3}_t$,
  where
  \begin{align*}
      S^{m,\rho,2,1}_t
      &=
        \sum_{l=1}^{L-1}
        \int_{0<\rho_l<\cdots<\rho_1<\rho}
          d\rho_1
          \cdots
          d\rho_l\,
          \sum_{i=1}^{2^mt}
            \tilde{N}^{m,\rho_1,\ldots,\rho_l}_{\tmi}
			(\tJm_{\tmi})^{-1}
            c(Y_{\tmim})d^m_{\tmim,\tmi},\\ 
      S^{m,\rho,2,2}_t
      &=
        \sum_{l=1}^{L-1}
        \int_{0<\rho_l<\cdots<\rho_1<\rho}
          d\rho_1
          \cdots
          d\rho_l\,
          \sum_{i=1}^{2^mt}
            R^{m,\rho_1,\ldots,\rho_l}_{\tmi}
            (\tJm_{\tmi})^{-1}
            c(Y_{\tmim})d^m_{\tmim,\tmi},\\ 
      S^{m,\rho,2,3}_t
      &=
        \int_{0<\rho_{L}<\cdots<\rho_1<\rho}
          d\rho_1
          \cdots
          d\rho_L\,
          \sum_{i=1}^{2^mt}
            N^{m,\rho_1}_{\tmi}
            \cdots 
            N^{m,\rho_L}_{\tmi}
            (\tilde{J}^{m,\rho_L}_{\tmi})^{-1}
            c(Y_{\tmim})d^m_{\tmim,\tmi}.
  \end{align*}
  We estimate the terms above.
 
	By the definition, 
	all terms in the expansion of $R^{m,\rho_1,\ldots,\rho_l}_t$
	are given by the product of $l$ terms from 
	$N^{m,\rho_r}_t$ and 
	$\tilde{Z}^{m,\rho_r,\nu}_t I^{m,\rho_r}(\varphi_\nu)_t$
	($1\leq r\leq l$, $1\leq \nu\leq n$)
  	and each term contains at least one $\tilde{Z}^{m,\rho_r,\nu}_t I^{m,\rho_r}(\varphi_\nu)_t$.
  Thus, using Remark~\ref{discrete young integral}, 
  Lemmas~\ref{estimate of discrete rough integral 2} and~\ref{Nmr},
	we have
  \begin{multline*}
	\bigg\|
		(2^m)^{2H-\frac{1}{2}}
		\sum_{i=1}^{2^m \cdot}
			R^{m,\rho_1,\ldots,\rho_l}_{\tmi}
			(\tJm_{\tmi})^{-1}
			c(Y_{\tmim})
			d^m_{\tmim,\tmi}
	\Big\|_{\lambda_1}\\
	\begin{aligned}
		&\leq
			C
			\|
				(2^m)^{2H-\frac{1}{2}}
				R^{m,\rho_1,\ldots,\rho_l}
			\|_{H^-}
			\bigg\|
				\sum_{i=1}^{2^m \cdot}
					(\tJm_{\tmi})^{-1}
					c(Y_{\tmim})
					d^m_{\tmim,\tmi}
			\bigg\|_{\lambda_1}\\
		&\leq
			C
			\|(2^m)^{2H-\frac{1}{2}}\tZmr\|_{H^-}
			\cdot
			a_m 
			C
			G_1,
	\end{aligned}
  \end{multline*}
  from which we obtain an estimate of $S^{m,\rho,2,2}$.
  We next consider $S^{m,\rho,2,1}$.
  Noting \eqref{eq48930184231},
  we have
  \begin{multline*}
    \sum_{i=1}^{2^mt}
      \tilde{N}^{m,\rho_1,\ldots,\rho_l}_{\tmi}
      (\tJm_{\tmi})^{-1}
      c(Y_{\tmim})d^m_{\tmim,\tmi}
      =
      \sum_{i=1}^{2^mt}
      \tilde{N}^{m,\rho_1,\ldots,\rho_l}_{\tmi}
      (I+K^{2,m,R}_{\tmi})
      J_{\tmim}^{-1}
      c(Y_{\tmim})d^m_{\tmim,\tmi}\\
      \begin{aligned}
        &+
        \sum_{i=1}^{2^mt}
        \tilde{N}^{m,\rho_1,\ldots,\rho_l}_{\tmi}
        (I+K^{2,m,R}_{\tmi})
        J_{\tmim,\tmi}^{-1}
        c(Y_{\tmim})
        d^m_{\tmim,\tmi}\\
        &+
        \sum_{i=1}^{2^mt}
        \tilde{N}^{m,\rho_1,\ldots,\rho_l}_{\tmi}
        L^{2,m,R}_{\tmi}
        c(Y_{\tmim})
        d^m_{\tmim,\tmi}.
      \end{aligned}
  \end{multline*}
  All terms can be treated in the similar way as Lemma~\ref{Smrho3}
  because $\tilde{N}^{m,\rho_1,\ldots,\rho_l}$
  is an $\{a_m\}$-order nice discrete
  process independent of $\rho_1,\ldots,\rho_l$
  (see Lemma~\ref{Nmr}).

  Finally, we consider $S^{m,\rho,2,3}$.
  Noting that
  \begin{align*}
    \sup_{\rho_1,\ldots,\rho_L}\|N^{m,\rho_1}\cdots N^{m,\rho_L}\|_{H^-}
    =O(a_m^L),  
  \end{align*}
  we see that this term is small for large $L$.
  This completes the proof.
  \end{proof}

\begin{proof}[Proof of Lemma~\upshape{\ref{estimate of tzmr}}]
We write
$
	f_m
	=
		\sup_{\rho}\|(2^m)^{2H-\frac{1}{2}}\tZmr\|_{H^-}1_{\Omegam}
$.
From the lemmas above, there exist
random variables $\{\varGamma_m\}$ and $\varGamma$ defined on $\Omega_0$
which satisfy $\sup_m\|\varGamma_m\|_{L^p}<\infty$ for all $p\ge 1$
and $\varGamma\in \cap_{p\ge 1}L^p(\Omega_0)$
such that
$
 f_m\le \varGamma_m+a_m \varGamma f_m
$.
Recalling $\tZmr$ is an $\{a_m\}$-order nice discrete process independent of $\rho$
(Theorem~\ref{rough estimate of tzmr}),
there exists $\varGamma'$ such that $f_m\leq (2^m)^{2H-\frac{1}{2}}\varGamma'$
and $\varGamma'\in \cap_{p\ge 1}L^p(\Omega_0)$.
By using this inequality $L$-times and
Theorem~\ref{rough estimate of tzmr}, we get
\begin{align*}
	f_m
	&\le
		\left\{
			\sum_{l=0}^{L-1}
				(a_m \varGamma)^l
		\right\}
		\varGamma_m
		+
		(a_m \varGamma)^{L}
		f_m
	\le
		\left\{
			\sum_{l=0}^{L-1}
				(a_m \varGamma)^l
		\right\}
		\varGamma_m
		+
		(2^m)^{2H-\frac{1}{2}}(a_m \varGamma)^{L}
		\varGamma'.
\end{align*}
By taking $L$ to be sufficiently large, we arrive at the conclusion.
\end{proof}

Finally, using Lemma~\ref{estimate of tzmr}, we prove an 
estimate of $\tZmr-I^m$.

\begin{lemma}\label{lem48901801}
	Let $\vep_1$ and $\vep_2$ be constants 
specified in Conditions~\ref{condition on dm}
	and \ref{condition on K}, respectively.
	Let $0<\vep<\min \{3H^--1,4H^--2H-\frac{1}{2},\vep_1,\vep_2\}$.
	Then, for all $p\ge 1$ it holds that
  \begin{align*}
    \lim_{m\to\infty}
      \left\|
        \sup_{0\le \rho\le 1}
          \|(2^m)^{2H-\frac{1}{2}+\vep}(\tZmr-I^m)\|_{H^-}1_{\Omega^{(m,d^m)}_0}\right\|_{L^p}
    =
      0.
  \end{align*}
\end{lemma}
\begin{proof}
	Write $f_m=\sup_{\rho}\|(2^m)^{2H-\frac{1}{2}}\tZmr\|_{H^-}1_{\Omegam}$.
	Lemmas~\ref{lemma for Smrho1}, \ref{Smrho2} 
	imply
	\begin{align*}
		\|(2^m)^{2H-\frac{1}{2}+\vep}S^{m,\rho,1}\|_{H^-}1_{\Omegam}
		&\leq
			(2^m)^\vep
			\cdot
			a_m
			CG_1
			f_m,\\
		\|(2^m)^{2H-\frac{1}{2}+\vep_1}S^{m,\rho,2}\|_{H^-}1_{\Omegam}
		&\leq
			(2^m)^\vep
			\cdot
			a_m
		  	C
			\big\{
				G_1f_m
				+X_m
				+G_2
				+1
			\big\}.
	\end{align*}
	Lemmas~\ref{lemma for Smrho5} and \ref{Smrho3} gives 
	similar estimates for 
	$\|(2^m)^{2H-\frac{1}{2}+\vep}S^{m,\rho,r}\|_{H^-}1_{\Omegam}$
	for $r=3,4,5$. Combining these estimates and Lemma~\ref{estimate of tzmr},
the proof is finished.
\end{proof}

\subsection{Proofs of Theorem~\ref{main theorem} and Corollary~\ref{corollary for main theorem}}
\label{subsection proof of main theorem}
Here we show Theorem~\ref{main theorem} and Corollary~\ref{corollary for main theorem}.

\begin{proof}[Proof of Theorem~\textup{\ref{main theorem}}]
	Recall that $R^m_t$ $(t\in \Dm)$ is defined by \eqref{remainder term discrete}.
	We will first consider $R^m_t1_{\Omegam}$,
	then $R^m_t1_{(\Omegam)^\complement}$.
	Proposition~\ref{expression of ymh-ym} implies
\begin{align*}
  R^m_t1_{\Omegam}
  =
    (\Ymh_t-Y_t-J_tI^m_t)1_{\Omegam}
  =
    \int_0^1
      \{ \tJmr_t \tZmr_t - J_tI^m_t \}
      1_{\Omegam}\,
      d\rho.
\end{align*}
The integrand scaled by $(2^m)^{2H-\frac{1}{2}+\vep}$ is decomposed into
\begin{multline*}
	(2^m)^{2H-\frac{1}{2}+\vep}
  \{ \tJmr_t \tZmr_t - J_tI^m_t \}
      1_{\Omegam}
  =
    \tJmr_t
    \cdot
	(2^m)^{2H-\frac{1}{2}+\vep}
    \big(
      \tZmr_t-I^m_t
    \big)
    1_{\Omegam}\\
    +
	(2^m)^{\vep}
    \big(\tJmr_t-J_t\big)
    1_{\Omegam}
    \cdot
	(2^m)^{2H-\frac{1}{2}}
    I^m_t.
\end{multline*}
Hence we have
\begin{multline*}
	(2^m)^{2H-\frac{1}{2}+\vep}
  	\max_{t\in\Dm}|R^m_t1_{\Omegam}|\\
	\leq
		\big(
			\max_t | \tJmr_t |
		\big)
		\big(
			\sup_{\rho}
    	    	\|(2^m)^{2H-\frac{1}{2}+\vep}(\tZmr-I^m)\|_{H^-}
		\big)
		1_{\Omega^{(m,d^m)}_0}\\
		+
		(2^m)^{\vep}
		\big(
			\sup_{t,\rho}|\tJmr_t-J_t|
		\big)
		1_{\Omega^{(m,d^m)}_0}
		\cdot
		\|(2^m)^{2H-\frac{1}{2}}I^m|_{\Dm}\|_{H^-}.
\end{multline*}
Here, $I^m|_{\Dm}$ denote the discrete process defined as the restriction of $I^m$ on $\Dm$.
The first term in the right-hand side converges to $0$
due to Lemmas~\ref{L^p estimate for tJmr} and \ref{lem48901801}.
The second term converges to $0$ follows from Theorem~\ref{convergence of J and J-1}
and Condition~\ref{minimal moment estimate}.
From this we have
$
  (2^m)^{2H-\frac{1}{2}+\vep}
  \max_{t\in\Dm}|R^m_t1_{\Omegam}|
$
converges to $0$ in $L^p$.

Next we consider $R^m_t1_{(\Omegam)^\complement}$.
Noting
\begin{align*}
  (2^m)^{2H-\frac{1}{2}+\vep}
  R^m_t1_{(\Omegam)^\complement}
  &=
    \big(\Ymh_t-Y_t\big)
    \cdot
    (2^m)^{2H-\frac{1}{2}+\vep}
    1_{(\Omegam)^\complement}\\
  &\phantom{\leq}\qquad
    -
    J_t
    \cdot
    (2^m)^{2H-\frac{1}{2}}
    I^m_t
    \cdot
    (2^m)^{\vep}
    1_{(\Omegam)^\complement},
\end{align*}
we have
\begin{multline*}
  (2^m)^{2H-\frac{1}{2}+\vep}
  \max_{t\in\Dm} |R^m_t1_{(\Omegam)^\complement}|
  \leq
	\big(
		\max_t | \Ymh_t-Y_t |
	\big)
    \cdot
    (2^m)^{2H-\frac{1}{2}+\vep}
    1_{(\Omegam)^\complement}\\
    +
    \big(
      \max_t | J_t |
    \big)
    \|(2^m)^{2H-\frac{1}{2}}I^m|_{\Dm}\|_{H^-}
    \cdot
    (2^m)^{\vep}
    1_{(\Omegam)^\complement}.
\end{multline*}
Lemma~\ref{Ist general} and Remark~\ref{estimate for J and J^-1}
imply that $\max_t | \Ymh_t - Y_t |$ and $\max_t | J_t |$
are bounded from above by $\cap_{p\ge 1}L^p$ random variable.
By using \eqref{complement of Omegam} and Condition~\ref{minimal moment estimate},
both terms of the right-hand side converge to $0$ in $L^p$.
The proof is completed.
\end{proof}

\begin{proof}[Proof of Corollary~\upshape{\ref{corollary for main theorem}}]
	Recall that $R^m_t$ $(0\leq t\leq 1)$ is defined by \eqref{remainder term conti}.
	Since $R^m_t=R^m_{\tmkm}+(R^m_t-R^m_{\tmkm})$ for $\tmkm\leq t\leq \tmk$, we have
	\begin{align*}
		\max_{0\leq t\leq 1} |R^m_t|
		\leq
			\max_{t\in \Dm} |R^m_t|
			+
			\max_{1\leq k\leq 2^m}
			\max_{\tmkm\leq t\leq \tmk}
				|R^m_t-R^m_{\tmkm}|.
	\end{align*}
	Since the first term is estimated in Theorem~\textup{\ref{main theorem}},
	we give an estimate of the second term.
	Let $\tmkm\leq t\leq \tmk$.
	We decompose $R^m_t-R^m_{\tmkm}$ into two terms;
\begin{align*}
  \Phi^m_1(t)
  &=
    \Ymh_t-\Ymh_{\tmkm}-(Y_t-Y_{\tmkm}),
  &
  \Phi^m_2(t)
  &=
    J_{\tmkm}I^m_{\tmkm}-J_tI^m_t.
\end{align*}
We have
\begin{align*}
  \Phi^m_1(t)&=\big\{\sigma(\Ymh_{\tmkm})-\sigma(Y_{\tmkm})\big\}B_{\tmkm,t}
+\big\{((D\sigma)[\sigma])(\Ymh_{\tmkm})-((D\sigma)[\sigma])(Y_{\tmkm})\big\}
\BB_{\tmkm,t}\\
&\qquad +\big\{b(\Ymh_{\tmkm})-b(Y_{\tmkm})\big\}(t-\tmkm)+
c(\Ymh_{\tmkm})d^{m}_{\tmkm,t}+\big\{\hepm_{\tmkm,t}-\epm_{\tmkm,t}\big\},
\end{align*}
which implies
\begin{align*}
	|\Phi^m_1(t)|
	&\leq
		C
		\big\{
			|\Ymh_{\tmkm}-Y_{\tmkm}|\Delta_m^{H^-}
			+\hat{X}\Delta_m^{2H^-}
			+\hat{X}\Delta_m^{3H^-}
		\big\}\\
	&\leq
		C
		\big\{
			|J_{\tmkm}I^m_{\tmkm}+R^m_{\tmkm}|\Delta_m^{H^-}
			+\hat{X}\Delta_m^{2H^-}
		\big\}.
\end{align*}
Here $C$ is a constant depending on $\sigma$, $b$, $c$ and $C(B)$.
From this we obtain
\begin{align*}
	(2^m)^{2H-\frac{1}{2}+\vep}
	|\Phi^m_1(t)|
	&\leq
		C
		\big(
			1+\|J\|_{H^-}
		\big)
		\|(2^m)^{2H-\frac{1}{2}}I^m|_{\Dm}\|_{H^-}
		\Delta_m^{H^--\vep}\\
	&\qquad
		+
		C
		\big\{
			(2^m)^{2H-\frac{1}{2}+\vep}
			\max_k|R^m_{\tmk}|
		\big\}
		\Delta_m^{H^-}
		+
		C
		\hat{X}
		\Delta_m^{\frac{1}{2}-2H+2H^--\vep}\\
	&=:
		CX_{m,1} \Delta_m^{H^--\vep}
		+CX_{m,2} \Delta_m^{H^-}
		+C
		\hat{X}
		\Delta_m^{\frac{1}{2}-2H+2H^--\vep}.
\end{align*}
We have $I^m_t=I^m_{\tmkm}$ $(\tmkm\le t\le \tmk)$, which implies
\begin{align*}
	(2^m)^{2H-\frac{1}{2}+\vep}|\Phi^m_2(t)|
 	=
		(2^m)^{\vep}
		|J_{\tmkm}-J_t| 
		|(2^m)^{2H-\frac{1}{2}}I^m_{\tmkm}|
	\le
		X_{m,1}
		\Delta_m^{H^--\vep}.
\end{align*}
Noting that the right-hand sides in the two estimates
are independent of $k$, we have
\begin{align*}
	(2^m)^{2H-\frac{1}{2}+\vep}
	\max_{1\leq k\leq 2^m}
	\max_{\tmkm\leq t\leq \tmk}
		|R^m_t-R^m_{\tmkm}|
	&\leq
		(C+1)X_{m,1} \Delta_m^{H^--\vep}
		+CX_{m,2} \Delta_m^{H^-}\\
	&\phantom{\leq}
		\qquad
		+C
		\hat{X}
		\Delta_m^{\frac{1}{2}-2H+2H^--\vep}
\end{align*}
We see that $\sup_m\{\|X_{m,1}\|_{L^p},\|X_{m,2}\|_{L^p},
\|\hat{X}\|_{L^p}\}<\infty$ for all $p\ge 1$ which follows from
Lemma~\ref{L^p estimate for tJmr}, Remark~\ref{estimate for J and J^-1},
Condition~\ref{minimal moment estimate} and Theorem~\ref{main theorem}.
Hence noting $3H^--1\leq H^-$ and $3H^--1\leq \frac{1}{2}-2H+2H^-$,
we complete the proof.
\end{proof}

\noindent
{\bf Acknowledgment}
~~	The authors thank the anonymous referees for their helpful comments
	which improved the quality of this paper.
	This work was supported by JSPS KAKENHI Grant Numbers JP20H01804 and JP22K13932.

\bigskip

\address{
Shigeki Aida\\
Graduate School of Mathematical Sciences,\\
The University of Tokyo,\\
Meguro-ku, Tokyo, 153-8914, Japan}
{aida@ms.u-tokyo.ac.jp}

\address{
Nobuaki Naganuma\\
Faculty of Advanced Science and Technology,\\
Kumamoto University,\\
Kumamoto city, Kumamoto, 860-8555, Japan}
{naganuma@kumamoto-u.ac.jp}

\end{document}